\documentclass[final,leqno]{siamltex}
\usepackage{threeparttable}
\usepackage{dcolumn}
\newcolumntype{d}{D{.}{.}{-1}}
\usepackage{subfigure}
\usepackage[german,english]{babel}
\usepackage{amssymb}
\usepackage{algorithm}
\usepackage{algorithmic}
\usepackage{color}
\usepackage{graphicx}
\usepackage{booktabs}
\usepackage{amsfonts}
\usepackage{soul}
\usepackage{mathrsfs}
\usepackage{hyperref}
\usepackage{amsmath, amssymb}
\usepackage{float}

\usepackage{footnote}
\usepackage{threeparttable}
\usepackage{dcolumn}
\usepackage{fullpage}
\newcolumntype{d}{D{.}{.}{-1}}
\newtheorem{remark}[theorem]{Remark}
\usepackage{enumitem}
\usepackage[small,nohug,heads=vee]{diagrams}
\diagramstyle[labelstyle=\scriptstyle]
\newcommand{\defeq}{{: =}}
\newcommand{\hurwitz}[1]{\mathrm{H}(#1)}
\newcommand{\generalizedham}[1]{\mathrm{GH}(#1)}
\newcommand{\skewsymmetric}[1]{\mathrm{SS}(#1)}
\newcommand{\Sp}{\mathrm{Sp}}
\newcommand{\Vspace}{\mathbb{V}}
\newcommand{\Wspace}{\mathbb{W}}
\newcommand{\SPD}[1]{\mathrm{SPD}(#1)}
\newcommand{\SPSD}[1]{\mathrm{SPSD}(#1)}
\newcommand{\lift}{\phi}
\newcommand{\projection}{\psi}

\renewcommand{\algorithmicrequire}{\textbf{Input:}}
\renewcommand{\algorithmicensure}{\textbf{Output:}}

\newcommand{\RRstar}[1]{\mathbb{R}_*^{#1}}
\newcommand{\RR}[1]{\mathbb{R}^{#1}}

\usepackage{tikz}
\usepackage{tikz-cd}
\usetikzlibrary{%
  matrix,%
  calc,%
  arrows%
}

\title{Structure-preserving  model reduction for\\  marginally stable  LTI
systems
}

\author{Liqian Peng\thanks{Extreme Scale Data Science \& Analytics Department, Sandia National Laboratories, Livermore, CA 94550 (lpeng@sandia.gov).} \and Kevin Carlberg\thanks{Extreme Scale Data Science \& Analytics Department, Sandia National Laboratories, Livermore, CA 94550 (ktcarlb@sandia.gov).}}
\usepackage{xr}
\begin{document}
\setlength{\abovedisplayskip}{3pt}
\setlength{\belowdisplayskip}{3pt}
\setlength{\abovedisplayshortskip}{3pt}
\setlength{\belowdisplayshortskip}{3pt}

\maketitle
\begin{abstract}
This work proposes a structure-preserving model reduction method for marginally stable linear time-invariant (LTI) systems. In contrast to Lyapunov-stability-based approaches---which ensure the poles of the reduced system remain in the open left-half plane---the proposed method preserves marginal stability by reducing the subsystem with poles on the imaginary axis in a manner that ensures those poles remain purely imaginary. In particular, the proposed method decomposes a marginally stable LTI system into (1) an asymptotically stable subsystem with eigenvalues in the open left-half plane and (2) a pure marginally stable subsystem with a purely imaginary spectrum. We propose a method based on inner-product projection and the Lyapunov inequality to reduce the first subsystem while preserving asymptotic stability.  In addition, we demonstrate that the pure marginally stable subsystem is a generalized Hamiltonian system; we then propose a method based on symplectic projection to reduce this subsystem while preserving pure marginal stability.  In addition, we propose both inner-product and symplectic balancing methods that balance the operators associated with two quadratic energy functionals while preserving asymptotic and pure marginal stability, respectively.  We formulate a geometric perspective that enables a unified comparison of the proposed inner-product and symplectic projection methods. Numerical examples illustrate the ability of the method to reduce the dimensionality of marginally stable LTI systems while retaining accuracy and preserving marginal stability; further, the resulting reduced-order model yields a finite infinite-time energy, which arises from the pure marginally
stable subsystem.
\end{abstract}



\section{Introduction}

Reduced-order models (ROMs) are essential for enabling high-fidelity
computational models to be used in many-query and real-time applications such
as control, optimization, and  uncertainty quantification.  Marginally stable
linear time-invariant dynamical (LTI) systems often arise in such
applications; examples include inviscid fluid flow, quantum mechanics,  and
undamped structural dynamics.  An ideal model-reduction approach for such
systems would produce a dynamical-system model that is lower dimensional, is
accurate with respect to the original model, and remains marginally stable,
which is an intrinsic property of the dynamical system (it ensures, e.g., a
finite system response at infinite time).  Unfortunately, most classical
model-reduction methodologies, such as balanced truncation \cite{MooreBC:81a},
Hankel norm approximation  \cite{GloverK:84a}, optimal $\mathcal H_2$
approximation \cite{GugercinS:08a,WilsonDA:07a,MeierL:67a}, and Galerkin
projection exploiting inner-product structure \cite{RowleyCW:04a}, were
originally developed for asymptotically stable LTI systems, i.e., systems with
all poles in the open left half-plane.

Although developed for asymptotically stable systems, balanced truncation and optimal $\mathcal H_2$ approximation 
can be extended to unstable stable systems without
poles on the imaginary axis. In particular,  a reduced-order model can be obtained
by balancing and truncating frequency-domain controllability and observability
Gramians~\cite{PrakashR:90a, ZhouK:99a}.   
By extending the $\mathcal H_2$ norm  to the $\mathcal L_2$-induced
Hilbert-Schmidt norm,  an iteratively corrected rational Krylov algorithm  was
proposed for optimal $\mathcal L_2$  model reduction \cite{GugercinS:10a}.
However,  the methods in Refs.~\cite{PrakashR:90a, ZhouK:99a, GugercinS:10a} cannot
be applied to marginally stable systems, as the frequency-domain
controllability and observability Gramians as well as the $\mathcal L_2$-induced
Hilbert-Schmidt norm are not well defined when there are poles on the imaginary
axis.

Although many well-known model reduction methods can be directly applied to
systems with purely imaginary poles, they do not guarantee stability. These
methods include proper orthogonal decomposition (POD)--Galerkin~\cite{HolmesP:12a}, balanced POD~\cite{RowleyCW:05a},
pseudo balanced POD~\cite{RowleyCW:11a,OrAC:10a}, and moment matching
\cite{BaiZ:02a, FreundRW:03a}.  The shift-reduce-shift-back approach (SRSB)
\cite{YangJ:92a, DoyleJC:02a, SantiagoJM:85a, ZilouchianA:91a, YangJ:93a}
reduces a $\mu$-shifted system $(A-\mu I, B, C)$ by balanced truncation.
However, this approach fails to ensure stability when the balanced reduced
system is shifted back by $\mu$.


In general,  stability-preserving  ROMs fall into roughly two categories.
The first category of methods derives \textit{a priori} a stability-preserving
model reduction framework, often specific to a particular equation set; the
present work falls within this category.
Refs.~\cite{RowleyCW:04a, BaroneMF:09a,KalashnikovaI:10a} construct ROMs  in
an energy-based inner product. Ref.~\cite{SerreG:12a} extends
Ref.~\cite{BaroneMF:09a,KalashnikovaI:10a} by applying the stabilizing
projection to a skew-symmetric system constructed by augmenting a given linear
system with its dual system.  
Refs.~\cite{Marsden:03h, CarlbergK:12a, CarlbergK:14a, vdSchaft:90a,
HartmannC:10a,vdSchaft:10a,vdSchaft:12a,PengL:16a, PengL:16c, HesthavenJS:17a} construct reduced-order
models to preserve the Lagrangian and (port-)Hamiltonian structures of the
original  systems.  However, these methods cannot be applied to general marginally stable LTI systems.

 The second category of methods stabilizes an unstable ROM
through \textit{a posteriori} stabilization step. In particular, Ref.~\cite{KalashnikovaI:14a} stabilizes reduced-order models via optimization-based eigenvalue reassignment. Refs.~\cite{BondBN:08a, AmsallemD:12a, BalajewiczM:16a}  construct reduced basis via minimal subspace rotation on the Stiefel manifold while preserving certain properties of the original system matrix.  Other methods includes to introduce viscosity~\cite{AubryN:88a, PodvinB:88a,DelvilleJ:99a}  or penalty term~\cite{CazemierW:98a}, to enrich basis functions  representing the small and energy dissipation scale~\cite{BalajewiczMJ:13a, NoackBR:03a,BergmannM:09a},  and to calibrate POD coefficients~\cite{CoupletM:05a,KalbVL:07a}.
 In many cases, the stabilization  alters the original unstable ROM and a sacrifice of accuracy is inevitable.

 In this work, we propose a structure-preserving model-reduction method for
 marginally stable systems. The method guarantees marginal-stability
 preservation by executing two steps.  First, the approach decomposes the
 original marginally stable linear system into two subsystems: one with
 eigenvalues in the left-half plane and one with nonzero eigenvalues on the
 imaginary axis.  This is similar to the approach taken in
 Ref.~\cite{MirnateghiN:13a,ZhouK:99a} for performing model reduction of
 unstable systems without poles on the imaginary axis.  Specifically, given a
 marginally stable (autonomous) LTI system $\dot x=A x$, where $A$ is
 invertible and all eigenvalues have a non-positive real part, we apply a
 similarity transformation, which yields $A=T \diag(A_s, A_m) T^{-1} $. Here,
 $A_s$ has eigenvalues in the left-half plane (i.e., is Hurwitz) and $A_m$ has
 purely imaginary eigenvalues.  In this case, the subsystem $\dot x_s =A_s
 x_s$ is asymptotically stable,  while we show that the subsystem $\dot
 x_m=A_m x_m$ is a generalized Hamiltonian system.  Second, the method
 performs structure-preserving model reduction on the subsystems separately;
 namely, inner-product projection based on the Lyapunov inequality is employed to
 reduce the asymptotically stable subsystem, while symplectic projection is
 applied to the pure marginally stable subsystem characterized by purely
 imaginary eigenvalues.

Specific contributions of this work include:

\begin{enumerate} 
\item A novel structure-preserving model reduction method for
marginally stable LTI systems that preserves the asymptotic stability
of the asymptotically stable subsystem via inner-product projection
and the pure marginal stability of pure marginally stable subsystem
via symplectic projection (Algorithm \ref{alg:overall}). 
\item A general inner-product projection framework  (Section
\ref{sec:asystab}), which we demonstrate ensures asymptotic-stability
preservation if the matrix used to define the inner product satisfies the Lyapunov
inequality (Lemma \ref{thm:stabpresadj}).
\item An inner-product balancing approach that enables the operators
associated with any primal or dual quadratic energy functional to be
balanced (Section \ref{sec:innerproductbalancing}). If either of these
satisfies a Lyapunov inequality,
then asymptotic stability is additionally preserved (Corollary
\ref{cor:innerproductbalancestability}). We show that many existing
model-reduction techniques (e.g., POD--Galerkin, balanced truncation, balanced
POD, and SRSB) can be expressed as an inner-product projection and in fact are special
cases of inner-product balancing (Table \ref{tab:comp_inner}).
\item A stabilization approach that produces an asymptotically stable
reduced-order model starting with a subset of the ingredients required for a
stability-preserving inner-product projection, e.g., starting with an arbitrary  trial
basis matrix and a symmetric-positive-definite matrix that satisfies the Lyapunov inequality
(Section \ref{sec:existenceInner}).
\item Analysis that demonstrates that any pure marginally stable system is equivalent to a
generalized Hamiltonian system with marginal stability (Theorem \ref{thm:matrix_similar_Hamil}).
\item A novel symplectic-projection framework (Section \ref{sec:marstab}) that ensures
preservation of pure marginal stability (Theorem \ref{thm:mainstab}).
\item A symplectic balancing approach that enables the operators
associated with any primal or negative dual quadratic energy functional to be
balanced (Section \ref{sec:symplecticbalancing}) and preserve pure marginal stability
(Corollary \ref{cor:symplecticBalanceStable}). In particular, we show that the
generalized Hamiltonians associated with the primal and negative dual systems can be
balanced with this approach.
\item A stabilization approach that produces a pure marginally stable
reduced-order model starting with a subset of the ingredients required for a
symplectic projection (Section \ref{sec:PSD}).
\item A geometric framework that enables a unified analysis and comparison of
inner-product and symplectic projection (Tables \ref{tab:podvspsd} and \ref{tab:balancing}).
	\item Experiments on two model problems that demonstrate that the proposed
method has a small relative error in both the state and total energy (Section
\ref{sec:example}).  Because symplectic model reduction is energy-conserving,
the proposed method ensures that the infinite-time system energy is equal to
the initial energy of the marginally stable subsystem. In contrast, the
infinite-time energy of other reduced models is zero or infinity.
\end{enumerate}
	 
The remainder of the paper is organized as follows.  Section \ref{sec:general}
provides an overall view of the proposed method. Sections \ref{sec:asystab} and
\ref{sec:marstab} present the methodologies to reduce the asymptotically stable
subsystem and marginally stable subsystem, respectively.  Section
\ref{sec:example} illustrates the stability, accuracy, and efficiency of the
proposed method through two numerical examples.  Finally, Section
\ref{sec:conclusion} provides conclusions.

We make extensive use of the following sets in the remainder of the paper:
\begin{itemize} 
\item $\SPD{n}$: the set of
all $n\times n$ symmetric-positive-definite (SPD) matrices.
\item $\SPSD{n}$: the set of
all $n\times n$ symmetric-positive-semidefinite (SPSD) matrices.
\item$\skewsymmetric{n}$: the set of 
$n\times n$ nonsingular, skew-symmetric matrices.
\item $\hurwitz{n}$: the
set of real-valued $n\times n$ matrices whose eigenvalues have strictly
negative real parts (i.e., the set of Hurwitz matrices).
\item $\generalizedham{n}$: the
set of real-valued $n\times n$ diagonalizable matrices with nonzero purely imaginary eigenvalues.
\item $\RRstar{n\times k}$: the set of full-column-rank $n\times k$ matrices
with $k\leq n$ (i.e., the non-compact Stiefel manifold).
\item $O(M,N)$: 
the set of full-column-rank $n\times k$ matrices $V$ with
$k\leq n$ such that $V^\tau MV = N$ with
$M\in\SPD{n}$ and $N\in\SPD{k}$. Note that $O(I_n,I_k)$ represents the
Stiefel manifold.
\item $\Sp(J_\Omega, J_\Pi)$: the set of full-column-rank $2n\times 2k$ matrices $V$ with
$k\leq n$ such that
$V^\tau J_\Omega V = J_\Pi$ with
$J_\Omega\in\skewsymmetric{2n}$ and $J_\Pi\in\skewsymmetric{2k}$.
Note that $\Sp(J_{2n},J_{2k})$ represents the
symplectic Stiefel manifold.
\end{itemize}

\section{Marginally stable LTI systems} \label{sec:general} 
We begin by formulating the full-order model, which is a marginally stable LTI
system (Section \ref{sec:fullmodel}), and subsequently present the formulation
for a general projection-based reduced-order model (Section
\ref{sec:redmodel}).  Then, we present the proposed framework based on system
decomposition (Section \ref{sec:decomp}).

\subsection{Full-order model} \label{sec:fullmodel}
This work considers continuous-time LTI systems of the form
\begin{equation}\label{eq:linear_control_sys}
\begin{aligned} 
\dot x &=A x+Bu\\
y &=Cx
 \end{aligned}
\end{equation}
with $A \in\mathbb{R}^{n\times n}$, $B\in \mathbb{R}^{n \times p}$, and $C
\in\mathbb{R}^{q \times n }$, $x \in\mathbb{R}^{n}$, $u \in\mathbb{R}^{p}$,
and $y \in \mathbb{R}^q$. We denote this system by $(A,B,C)$ and
focus on the particular case where the linear
system is marginally stable. Because stability concerns the spectrum of the
operator $A$, we focus primarily on the corresponding autonomous system
\begin{equation}\label{eq:linearsys}
 \dot x=A x.
 \end{equation}
We now define marginal  stability.
\begin{definition}[Marginal stability]
Linear system \eqref{eq:linear_control_sys}  is 
\emph{marginally stable}, or Lyapunov stable, if for every initial condition
$x(0)=x_0\in \mathbb{R}^n$,   the state response $x(t)$ of the associated
autonomous system \eqref{eq:linearsys} is uniformly bounded.
\end{definition}

The following standard lemmas (e.g.,
Ref.~\cite[pp.~66--70]{HespanhaJP:09a}) provide conditions for marginal stability.
\begin{lemma}\label{lemma:stability0}
The following conditions are equivalent:
\begin{enumerate}[label=(\alph*)]
\setlength{\itemsep}{5pt}
\item The system \eqref{eq:linear_control_sys} is marginally stable.
\item All eigenvalues of $A$ have non-positive real parts and all Jordan
blocks corresponding to eigenvalues with zero real parts are $1 \times 1$.
\end{enumerate}
\end{lemma}
\begin{lemma}\label{lemma:stabilityequiv1}
The system \eqref{eq:linear_control_sys} is marginally stable if one of the following
conditions holds:
\begin{enumerate}[label=(\alph*)]
\setlength{\itemsep}{5pt}
\item There exists 
$\Theta\in\SPD{n}$ that
satisfies the Lyapunov inequality
 \begin{equation}\label{eq:inequiLya}
 A^{\tau} \Theta + \Theta A \preceq 0.
 \end{equation}
\item For every $Q\in\SPSD{n}$, there
exists a unique solution $\Theta\in\SPD{n}$ to the Lyapunov equation
\begin{equation}\label{eq:semiequiLya}
A^{\tau} \Theta + \Theta A =- Q.
\end{equation}
\item There exists $\Theta\in\SPD{n}$ such that the energy $\frac{1}{2}x^\tau \Theta x$ of the corresponding autonomous system is nonincreasing in
time, i.e.,
 \begin{equation} \label{eq:Lyaenergy}
 \frac{d}{dt}\left(\frac{1}{2}x^\tau \Theta x\right) \leq 0,
  \end{equation} 
  for $x\in \mathbb{R}^n$ satisfying \eqref{eq:linearsys}.
 \end{enumerate}
\end{lemma}
\noindent We note that because $\frac{d}{dt}(\frac{1}{2}x^\tau \Theta
x) = \frac{1}{2}
x^\tau (A^\tau \Theta  +  \Theta A)x $, (a) and (c) are  equivalent.
Lemma \ref{lemma:stabilityequiv1} provides sufficient conditions
for marginal stability; not all marginally stable systems have a Lyapunov
matrix $\Theta$ that satisfies \eqref{eq:inequiLya}--\eqref{eq:Lyaenergy}.

\subsection{Reduced-order model}\label{sec:redmodel}
Let $\Psi, \Phi \in \RRstar{n\times k}$ denote  test and trial basis matrices
that are biorthogonal (i.e., $\Psi^{\tau} \Phi= I_k$) and whose columns span $k$-dimensional test and trial subspaces of $\RR{n}$,
respectively.
If the reduced-order model is constructed via Petrov--Galerkin projection performed on
the full-order model, then \eqref{eq:linear_control_sys} reduces to
\begin{equation}\label{eq:reduced_control_sys}
\begin{aligned} 
\dot z &= \tilde A z+\tilde Bu\\
y &=\tilde C z,
 \end{aligned}
\end{equation}
where $\tilde A\defeq\Psi^{\tau}  A \Phi\in \mathbb{R}^{k \times k}$, $\tilde
B\defeq \Psi^{\tau}
B\in \mathbb{R}^{k\times p}$, $\tilde C \defeq C\Phi \in \mathbb{R}^{q\times
k}$, and the state is approximated as $x\approx \Phi z$. We denote this system
by $(\tilde A,\tilde B,\tilde C)$. The 
corresponding autonomous system is
\begin{equation}\label{eq:xrAr}
\dot z =  \tilde A z
\end{equation}
with initial condition $z(0)=\Psi^{\tau} x_0\in \mathbb{R}^k$.

\subsection{System decomposition}\label{sec:decomp}
If the full-order-model system
\eqref{eq:linear_control_sys} is marginally stable and the matrix $A$ has a
full rank,  then
 there exists $T \in \RRstar{n\times n}$ such that
 the similarity transformation satisfies  
  \begin{equation}\label{eq:diagA}
 A= T
 \begin{bmatrix}
A_s & 0  \\
 0 & A_m
 \end{bmatrix} T^{-1},
 \end{equation}
  where $A_s \in \hurwitz{n_s}$, $A_m \in \generalizedham{n_m}$, and
	$n_s+n_m = n$. Let $T=\begin{bmatrix}T_s & T_m\end{bmatrix}$ with $T_s \in \RRstar{n \times n_s}$
	and $T_m \in \RRstar{n \times n_m}$. Then,  $A T_i =T_i A_i$ (for
	$i\in\{s,m\}$),  which
	implies that the columns of $T_i$  span an invariant subspace
	of $A$. Let $x_s\in \mathbb{R}^{n_s}$  and $x_m\in \mathbb{R}^{n_m}$. Substituting $x = T \begin{bmatrix} x_s^{\tau} &  x_m^{\tau} \end{bmatrix}^{\tau}$ into
	\eqref{eq:linear_control_sys} and premultiplying the first set of equations
	by $T^{-1}$ yields a decoupled LTI system 
  \begin{align}\label{eq:relinsys}
  \begin{split}
  \frac{d}{dt}
     \begin{bmatrix}
    x_s \\ x_m 
     \end{bmatrix}
  &=
  \begin{bmatrix}
A_s & 0  \\
 0 & A_m 
 \end{bmatrix} 
      \begin{bmatrix}
    x_s \\ x_m 
     \end{bmatrix} + 
		\begin{bmatrix} 
		B_s\\
		B_m
		\end{bmatrix} 
		u\\
		y &=\begin{bmatrix}
C_s & C_m
		\end{bmatrix}
     \begin{bmatrix}
    x_s \\ x_m 
     \end{bmatrix}
		 ,
  \end{split}
 \end{align}
 where $T^{-1}B = 
		\begin{bmatrix} 
		B_s^{\tau}& B_m^{\tau}
		\end{bmatrix} ^{\tau}$ and 
		$
		CT = 
		\begin{bmatrix}
C_s & C_m
		\end{bmatrix}
		$.
 Here, the subsystem associated with $x_s$ is asymptotically stable, while the
 subsystem associated with $x_m$ is marginally stable.

This decomposition enables each subsystem to be reduced in a manner that
preserves its particular notion of stability. 
In the present context, we can
accomplish this by defining biorthogonal test and trial basis matrices for each
subsystem $\Psi_i\in\RRstar{n_i\times k_i}$, $\Phi_i\in\RRstar{n_i\times
k_i}$, $i\in\{s,m\}$. Applying Petrov--Galerkin projection to
\eqref{eq:relinsys} with test basis matrix $\diag(\Psi_s,\Psi_m)$ and trial basis
matrix $\diag(\Phi_s,\Phi_m)$ yields a decoupled reduced LTI system
 \begin{align}\label{eq:general_relinsysSimple}
 \begin{split}
  \frac{d}{dt}
     \begin{bmatrix}
    z_s \\ z_m
     \end{bmatrix}
  &=
   \begin{bmatrix}
\tilde A_s & 0  \\
 0 & \tilde A_m  \\
 \end{bmatrix}
      \begin{bmatrix}
    z_s \\ z_m 
     \end{bmatrix}+
		\begin{bmatrix} 
		\tilde B_s\\
		\tilde B_m
		\end{bmatrix} 
		u\\
		y &=\begin{bmatrix}
\tilde C_s & \tilde C_m 
		\end{bmatrix}
     \begin{bmatrix}
    z_s \\ z_m 
     \end{bmatrix},
 \end{split}
 \end{align}
where the full state is approximated as 
$$x(t) \approx T\begin{bmatrix}
\Phi_s z_s(t)\\
\Phi_m z_m(t)
\end{bmatrix}.$$
Within this decomposition-based approach, basis matrices $\Psi_s$ and $\Phi_s$ can be
computed to preserve asymptotic stability in the associated reduced subsystem
(e.g., via balanced truncation or other Lyapunov methods). For the marginally
stable subsystem, we will show that the symplectic model reduction method can
be applied to obtain a low-order marginally stable system wherein all
eigenvalues of $\tilde A_m$ are nonzero and purely imaginary.

Algorithm \ref{alg:overall} summarizes the proposed procedure for computing
reduced-order-model operators $(\tilde A_s,\tilde B_s,\tilde C_s )$ and
$(\tilde A_m,\tilde B_m,\tilde C_m )$.  
Here, we have defined
Table \ref{tab:podvspsd} lists the
methods and key properties of each subsystem.  The next two sections explain
Algorithm \ref{alg:overall} and Table \ref{tab:podvspsd} in detail.  

\begin{algorithm}
\caption{Structure-preserving model reduction for marginally stable LTI
systems.
 } 
\begin{algorithmic}[1]
\label{alg:overall}
 \REQUIRE A marginally stable LTI system $(A,B,C)$.\\
\ENSURE Reduced-order-model operators 
$(\tilde A_s,\tilde
B_s,\tilde C_s )$ and
$(\tilde A_m,\tilde
B_m,\tilde C_m )$.
\STATE Compute a matrix $T$ such that $A$ is transformed into block-diagonal
form \eqref{eq:diagA}.
\STATE  Select $M \in \SPD{n_s}$ such that the Lyapunov inequality $A_s^{\tau} M
+ M A_s\prec 0$ is satisfied.\label{step:selectM}
\STATE Construct trial basis matrix $\Phi_s \in O(M, N)$ for some $N \in
\SPD{k_s}$, $k_s < n_s$.
\STATE  Construct test basis matrix $\Psi_s=M \Phi_s N^{-1}$.\label{step:selectPsis}
\STATE Construct the reduced system $\tilde A_s= \Psi_s^{\tau} A_s \Phi_s$, $\tilde
B_s = \Psi_s^{\tau}B_s$, $\tilde C_s = C_s \Phi_s$.
\STATE Select $J_\Omega \in \skewsymmetric{n_m}$ such that $A_m=-J_\Omega^{-1} L$ with $L\in \SPD{n_m}$.
\STATE Construct trial basis matrix $\Phi_m \in \Sp(J_\Omega, J_\Pi)$ for some 
$J_\Pi \in \skewsymmetric{k_m}$, $k_m < n_m$.
\STATE  Construct test basis matrix $\Psi_m=J_\Omega \Phi_m J_\Pi^{-1}$.
\STATE Construct the reduced system $\tilde A_m= \Psi_m^{\tau} A_m \Phi_m$, $\tilde
B_m = \Psi_m^{\tau}B_m$, $\tilde C_m = C_m \Phi_m$.
\end{algorithmic}
\end{algorithm}

Appendix
\ref{sec:decompGeneral} describes how this decomposition approach can be
extended to general unstable LTI systems with $A$ possibly singular.


\begin{table} 
\begin{center}
\caption{Inner-product model reduction v.\ symplectic model reduction.}

 \label{tab:podvspsd}
	{%
\begin{tabular}{|c|c |c|}
\hline
& \begin{tabular}{@{}c@{}} Asymptotically stable \\ subsystem \end{tabular}   &  
\begin{tabular}{@{}c@{}}   Marginally  stable \\ subsystem \end{tabular}  \\
\hline
\hline
 \begin{tabular}{@{}c@{}} Original  space \end{tabular} &
 \begin{tabular}{@{}c@{}}  Inner-product space: \\ $(\mathbb{R}^{n_s}, M)$
 with $M\in\SPD{n_s}$   \end{tabular}   &
 \begin{tabular}{@{}c@{}}  Symplectic space: \\ $(\mathbb{R}^{n_m}, J_\Omega)$
 with $J_\Omega \in \skewsymmetric{n_m}$ \end{tabular}    \\
\hline
 \begin{tabular}{@{}c@{}} System matrix \end{tabular} &
 \begin{tabular}{@{}c@{}}   $A_s   \in \hurwitz{n_s}$ \end{tabular}   &
 \begin{tabular}{@{}c@{}}     $A_m   \in \generalizedham{n_m}$  \end{tabular}    \\
\hline
\begin{tabular}{@{}c@{}}  Autonomous \\ system \end{tabular} &  \begin{tabular}{@{}c@{}}   $\dot x_s =A_s x_s$ \\ with $x_s \in \mathbb{R}^{n_s}$   \end{tabular} 
&  \begin{tabular}{@{}c@{}}   $\dot x_m =A_m x_m$ \\ with $x_m \in \mathbb{R}^{n_m}$   \end{tabular} \\
 \hline
 \begin{tabular}{@{}c@{}}  Key property \\ of full system \end{tabular} 
 &  \begin{tabular}{@{}c@{}}  Lyapunov inequality: \\ $A_s^{\tau} M+M A_s \prec 0$ \end{tabular}  &  
 \begin{tabular}{@{}c@{}}  Generalized Hamiltonian property: \\ $ A_m^{\tau} J_\Omega+ J_\Omega A_m=0$  \end{tabular}  \\
 \hline
\begin{tabular}{@{}c@{}} 
 Energy property \\
 of full system
 \end{tabular}
 & $\frac{d}{dt}\left(\frac{1}{2}x_s^\tau M x_s\right) < 0$ & 
$\frac{d}{dt}\left(\frac{1}{2}x_s^\tau L x_s\right) = 0$ with $A_m =
-J_\Omega^{-1}L$\\
 \hline
 \begin{tabular}{@{}c@{}}  Canonical  form \end{tabular} 
 &  \begin{tabular}{@{}c@{}}  $M=I_n$ \\ $A_s^{\tau}+ A_s \prec 0$ \end{tabular}  &  
 \begin{tabular}{@{}c@{}}   $J_\Omega=J_{2n}$  \\ $ A_m^{\tau} J_{2n}+ J_{2n}
 A_m=0$ \end{tabular}  \\
 \hline
 \begin{tabular}{@{}c@{}} Reduced  space \end{tabular} &
 \begin{tabular}{@{}c@{}}  Inner-product space: \\ $(\mathbb{R}^{k_s}, N)$
 with $N\in \SPD{k_s}$  \end{tabular}   &
 \begin{tabular}{@{}c@{}}  Symplectic space: \\ $(\mathbb{R}^{k_m}, J_\Pi)$
 with $J_\Pi \in \skewsymmetric{k_m}$ \end{tabular}    \\
 \hline
 Projection & Inner-product projection & Symplectic projection\\
 \hline
 \begin{tabular}{@{}c@{}}  Trial basis matrix\end{tabular}  &
   \begin{tabular}{@{}c@{}}
	 $\Phi_s\in O(M,N)$
   \end{tabular} &
   \begin{tabular}{@{}c@{}}
	 $\Phi_m \in \Sp(J_\Omega, J_\Pi)$
   \end{tabular} \\
   \hline
   \begin{tabular}{@{}c@{}} Test basis  matrix \end{tabular} &   $\Psi_s =M
	 \Phi_s N^{-1} \in \RRstar{n_s\times k_s}$ &  $\Psi_m= J_\Omega \Phi_m
	 J_\Pi^{-1} \in \RRstar{n_m\times k_m}$\\
\hline
   \begin{tabular}{@{}c@{}} Reduced-system matrix \end{tabular} &
	 \begin{tabular}{@{}c@{}} $\tilde A_s=\Psi_s^{\tau} A_s \Phi_s \in \hurwitz{k_s}$ \end{tabular}   &
  \begin{tabular}{@{}c@{}}$\tilde A_m =\Psi_m^{\tau} A_m \Phi_m \in \generalizedham{k_m}$
	\end{tabular} \\
\hline
   \begin{tabular}{@{}c@{}} Reduced autonomous \\ system \end{tabular} &  \begin{tabular}{@{}c@{}}   $\dot z_s = \tilde A_s z_s$ \end{tabular} &
 {\begin{tabular}{@{}c@{}}  $\dot z_m = \tilde A_m z_m$  \end{tabular} }     \\
  \hline
   \begin{tabular}{@{}c@{}}  Key property \\ of reduced system \end{tabular} 
		&  \begin{tabular}{@{}c@{}}  Lyapunov inequality: \\ $\tilde A_s^{\tau}
		N+N \tilde A_s \prec 0 $ \end{tabular}  &  
 \begin{tabular}{@{}c@{}}  Generalized Hamiltonian property: \\ $ \tilde A_m ^{\tau} J_\Pi+ J_\Pi \tilde A_m=0$  \end{tabular}   \\
 \hline
\begin{tabular}{@{}c@{}} 
 Energy property \\
 of reduced system
 \end{tabular}
 & $\frac{d}{dt}\left(\frac{1}{2}z_s^\tau N z_s\right) < 0$ & 
$\frac{d}{dt}\left(\frac{1}{2}z_s^\tau \tilde L z_s\right) = 0$ with $\tilde A_m =
-J_\Pi^{-1}\tilde L$\\
 \hline
   \begin{tabular}{@{}c@{}}  Approximate \\ solution \end{tabular} &  $x_s(t) \approx \Phi_s z_s(t)$ &  $x_m(t) \approx \Phi_m z_m(t)$  \\
 \hline
\end{tabular}}
\end{center}
\end{table}

\section{Reduction of asymptotically stable subsystems}\label{sec:asystab} 
This section focuses on reducing the asymptotically stable subsystem $\dot x_s
=A_s x_x$. Section  \ref{sec: inner-product_space} introduces inner-projection
projection, Section \ref{sec: inner-product_dyn} demonstrates that a
model-reduction method based on inner-projection projection preserves
asymptotic stability, Section \ref{sec:innerproductbalancing} presents the
inner-product-balancing framework, and Section \ref{sec:existenceInner}
describes methods for constructing the basis matrices that lead to a
inner-product projection given a subset of the required ingredients.
For notational simplicity, we omit the subscript $s$ throughout this section. 

\subsection{Asymptotically stable systems}\label{sec:asymp_stable_systems}
We begin by defining asymptotic stability.
\begin{definition}[Asymptotic stability]
Linear system \eqref{eq:linear_control_sys}  is \emph{asymptotically stable}
if, in addition to being marginally stable, $x(t)\to 0$  as $t\to \infty$ for
	every initial condition $x(0)=x_0\in \mathbb{R}^n$.
\end{definition}

\noindent In analogue to Lemmas \ref{lemma:stability0}--\ref{lemma:stabilityequiv1}, we
now provide conditions for asymptotic stability.
\begin{lemma}\label{lemma:stabilityequiv}
The following conditions are equivalent:
\begin{enumerate}[label=(\alph*)]
\setlength{\itemsep}{5pt}
\item The system \eqref{eq:linear_control_sys} is asymptotically stable.
\item $A\in\hurwitz{n}$.
\item There exists $\Theta\in\SPD{n}$ that satisfies the Lyapunov inequality
\begin{equation}\label{eq:equiLya4}
A^{\tau} \Theta + \Theta A \prec 0.
\end{equation}
\item For every $Q\in\SPD{n}$, there exists a unique Lyapunov matrix
$\Theta\in\SPD{n}$ that
satisfies \eqref{eq:semiequiLya}.
\item  There exists $\Theta\in\SPD{n}$ such that the energy $\frac{1}{2}x^\tau \Theta x$ of the corresponding autonomous system is strictly decreasing in
time, i.e.,
 \begin{equation} 
 \frac{d}{dt}\left(\frac{1}{2}x^\tau \Theta x\right) < 0,
  \end{equation} 
  for any $x\ne 0\in \mathbb{R}^n$ satisfying \eqref{eq:linearsys}.
\end{enumerate}
\end{lemma}
\noindent We note that 
$A\in\hurwitz{n}$ does not necessarily imply that the
symmetric part of $A$ is negative definite.
However, $A\in\hurwitz{n}$ if and only if it can be transformed into a matrix with negative
symmetric part by similarity transformation with a real matrix; see Lemma
\ref{lemma:spd_can} in Appendix \ref{sec:can_lya} for details. 


\noindent We now connect asymptotic stability of the primal system to that of
its dual.
\begin{lemma}[Dual version of Lemma \ref{lemma:stabilityequiv}]\label{lemma:stabilityequivDual}
If any condition of Lemma \ref{lemma:stabilityequiv} holds, then the following
conditions hold:
\begin{enumerate}[label=(\alph*)]
\setlength{\itemsep}{5pt}
\item The dual system $(A^\tau, C^\tau, B^\tau)$ is asymptotically stable.
\item $A^\tau\in\hurwitz{n}$.
\item There exists $\Theta'\in\SPD{n}$ that satisfies the dual Lyapunov inequality
\begin{equation}\label{eq:dualinequiLya}
A \Theta' + \Theta' A^\tau \prec 0.
\end{equation}
\item For every $Q'\in\SPD{n}$, there exists a unique Lyapunov matrix $\Theta'\in\SPD{n}$ that satisfies 
\begin{equation}\label{eq:semiequiLyaDual}
A \Theta' + \Theta' A^\tau =- Q'.
\end{equation}
\item  There exists $\Theta' \in\SPD{n}$ such that the energy $\frac{1}{2}x^\tau \Theta' x$ of the corresponding autonomous dual system is strictly decreasing in time, i.e.,
 \begin{equation} 
 \frac{d}{dt}\left(\frac{1}{2}x^\tau \Theta' x\right) < 0,
  \end{equation} 
  for any $x\ne 0 \in \mathbb{R}^n$ satisfying $\dot x=A^\tau x$.
\end{enumerate}
\end{lemma}
\begin{proof}
Because the eigenvalues of $A$ are identical to the eigenvalues of $A^\tau$,
$A\in\hurwitz{n}$ if and only if $ A^\tau \in\hurwitz{n}$. Thus the 
LTI system associated with $A^\tau$ is asymptotically stable
and satisfies the corresponding conditions of Lemma \ref{lemma:stabilityequiv}.
\hfill
\end{proof}

\begin{remark}[Relationship with negative dual system: asymptotic stability]\label{rmk:dual}
Thus, any method proposed in this work for ensuring asymptotic stability of a
given (sub)system also ensures asymptotic  stability of the associated dual
(sub)system.  However, because the trial basis $\Phi$ associated with $(\tilde
A,\tilde B,\tilde C)$ corresponds to the test basis of $(\tilde A^\tau, \tilde
C^\tau, \tilde B^\tau)$ (i.e., $\tilde A^\tau = \Phi^\tau A^\tau \Psi$), the
proposed methods for constructing a trial basis matrix $\Phi$ should be
applied to the dual system as a test basis matrix. Similarly, the proposed
methods for constructing a test basis matrix $\Psi$ should be applied to the
dual system as a trial basis matrix.
\end{remark}


\subsection{Inner-product projection of spaces} \label{sec: inner-product_space}
Let $\mathbb{V}\cong \RR{n}$ 
and $\mathbb{W}\cong \RR{k}$  with $k \leq n$
denote vector spaces equipped
with inner products
 $\left\langle \cdot, \cdot \right\rangle_{\mathbb{V}}: \mathbb{V}\times
 \mathbb{V} \to \mathbb{R}$ and
$\left\langle \cdot, \cdot \right\rangle_{\mathbb{W}}: \mathbb{W}\times
 \mathbb{W} \to \mathbb{R}$ 
 respectively.
These inner products  can be represented by matrices
$M\in\SPD{n}$ and $N\in\SPD{k}$, respectively, i.e.,
\begin{align*}
\left\langle  \hat x_1, \hat x_2 \right\rangle_{\mathbb{V}} &\equiv x_1^{\tau}
M x_2,\quad
\forall x_1,x_2 \in \RR{n}\\
\left\langle  \hat z_1, \hat z_2 \right\rangle_{\mathbb{W}} &\equiv z_1^{\tau}
N z_2,\quad
\forall z_1,z_2 \in \RR{k},
\end{align*}
where the operator $\hat\cdot$ provides the representation of an element of a
vector space from its coordinates, i.e., $\hat x \in\mathbb {V}$, $\forall
x\in\RR{n}$ and $\hat z\in\mathbb{W}$, $\forall z\in\RR{k}$.  We represent
these inner-product spaces $\Vspace$ and $\Wspace$ by $(\mathbb{R}^{n}, M)$
and $(\mathbb{R}^{k}, N)$ respectively.

\begin{definition}[Inner-product lift]
 An {\emph{inner-product lift}} is a linear mapping  $\lift: \Wspace \to
 \Vspace$ that preserves inner-product structure:
\begin{equation} \label{eq:innerpropres}
\left\langle \hat z_1, \hat z_2 \right\rangle_{\mathbb{W}}= \left\langle
\lift(\hat z_1), \lift(\hat z_2) \right\rangle_{\mathbb{V}}, \quad \forall \hat z_1, \hat z_2 \in \mathbb{W}.
\end{equation}
\end{definition}

\begin{definition}[Inner-product projection]
Let $\lift: \Wspace \to \Vspace$ be an inner-product lift. The adjoint of $\lift$
is the linear mapping $\projection: \Vspace  \to \Wspace$ satisfying
\begin{equation}\label{eq:innerproject}
\left\langle  \projection(\hat x) ,\hat z  \right\rangle_{\mathbb{W}} = \left\langle
 \hat x ,\lift(\hat z)\right\rangle_{\mathbb{V}}, \quad \forall \hat z\in \mathbb{W}, \ \hat x\in \mathbb{V}.
\end{equation}
 We say $\projection$ is the {\emph{inner-product projection}} induced by $\lift$.
\end{definition}

In coordinate space, this inner-product lift and projection can be expressed equivalently as
\begin{align*}
\lift(\hat z) &\equiv \Phi z,\quad \forall z \in \RR{k}\\
\projection(\hat x)&\equiv\Psi^{\tau} x,\quad \forall x\in\RR{n},
\end{align*}
respectively, where \eqref{eq:innerpropres}--\eqref{eq:innerproject} imply
that $\Phi \in \RRstar{n\times k}$  and $\Psi \in
\RRstar{n\times k}$ satisfy
\begin{align}\label{eq:innerpropres1}
\Phi^{\tau} M \Phi &=N\\
 \label{eq:innerpropres2}\Psi N &=M\Phi ,
\end{align}
from which it follows that
\begin{equation}\label{eq:defpsi}
\Psi =  M \Phi N^{-1}.
\end{equation}
For convenience, we write $\Phi\in O(M,N)$. Although $\Psi^{\tau}$ is not in general equal to the Moore--Penrose pseudoinverse $(\Phi^{\tau}
\Phi)^{-1}\Phi^{\tau}$, it can be verified that it is indeed a left inverse of $\Phi$, 
 which implies that $\projection \circ \lift$ is the identity map on $\mathbb{W}$.

\subsection{Inner-product projection of dynamics}\label{sec: inner-product_dyn}

This section describes the connection between inner-product 
projection and asymptotic-stability preservation in model reduction. Namely,
we show that if inner-product projection is employed to construct the
reduced-order model with $M$ corresponding to a Lyapunov matrix of the
original system, then the reduced-order model inherits asymptotic stability.

 \begin{definition}[Model reduction via inner-product projection]\label{def:inner-product-proj}
A reduced-order model 
$(\tilde A, \tilde B, \tilde C)$ with $\tilde A=\Psi^{\tau}  A \Phi$, $\tilde B=\Psi^\tau B$, and $\tilde C=C\Phi$
is constructed by an inner-product projection if
$\Phi \in O(M,N)$, $\Psi= M \Phi N^{-1}$, where $M\in\SPD{n}$ and
$N\in\SPD{k}$.
 \end{definition}

\begin{lemma}[Inner-product projection preserves asymptotic stability]\label{thm:stabpresadj}
If the original LTI system $(A, B, C)$ has a Lyapunov matrix $\Theta$
satisfying
\eqref{eq:equiLya4} and the
reduced-order model is constructed by inner-product projection with $M=\Theta$, then the
reduced-order model $(\tilde A, \tilde B, \tilde C)$ is asymptotically
stable with Lyapunov matrix $N$.
  \end{lemma}
 \begin{proof}
 Left- and right-multiplying  inequality
\eqref{eq:equiLya4} (with $\Theta=M$)
by $\Phi^{\tau}$ and $\Phi$, respectively,
yields
 \begin{equation}\label{eq:redLya}
\Phi^{\tau} A^{\tau} M \Phi+ \Phi^{\tau} M A \Phi\prec 0.
\end{equation}
 Substituting \eqref{eq:innerpropres2} 
and $\tilde A=\Psi^{\tau} A \Phi$ 
 in \eqref{eq:redLya}  yields
 \begin{equation}\label{eq:ANl0}
\tilde A^{\tau} N+ N \tilde A \prec 0,
\end{equation}
which implies that the reduced system is asymptotically stable by Lemma
\ref{lemma:stabilityequiv}.
\hfill
 \end{proof}


We note that Lemma \ref{thm:stabpresadj} is a generalization of the
stability-preservation property in Ref.~\cite{RowleyCW:04a}, which required the
reduced space to be Euclidean (i.e., $N=I_k$ in the present notation). Lemma
\ref{thm:stabpresadj} considers a more general form where the reduced space can
be any inner-product space, i.e., $N\in\SPD{k}$ but otherwise arbitrary.
\medskip



\subsection{Inner-product balancing}\label{sec:innerproductbalancing}

We now describe an inner-product-balancing approach that leverages
inner-product structure. Table \ref{tab:balancing} compares this
approach with a novel symplectic-balancing approach, which will be described
in Section \ref{sec:symplecticbalancing}.
\begin{table} 
\begin{center}
\caption{Inner-product balancing v.~symplectic balancing. Both methods require
inputs $\Xi,\Xi'
\in\SPD{n}$ and employ
decompositions $\Xi = RR^\tau$, $\Xi' = SS^\tau$, and $R^\tau S = U\Sigma V^\tau$ and
$\bar \Phi = SV_1\Sigma_1^{-1/2}$, $\bar \Psi = R U_1\Sigma_1^{-1/2}$.}
 \label{tab:balancing}
	%
\begin{tabular}{|c|c |c|}
\hline
 & Inner-product balancing & Symplectic balancing\\
\hline
\hline
\begin{tabular}{@{}c@{}}Primal energy $\frac{1}{2}x^\tau M x$ 
\end{tabular}& \begin{tabular}{@{}c@{}} $M=\Xi$ \end{tabular} &
\begin{tabular}{@{}c@{}} $M = G^{-\tau}\diag(\Xi,\Xi')G^{-1}$ 
\end{tabular}
\\
\hline
\begin{tabular}{@{}c@{}}Dual energy $\frac{1}{2}(x')^\tau M' x'$ 
\end{tabular} & \begin{tabular}{@{}c@{}} $M'=\Xi'$ \end{tabular} &
\begin{tabular}{@{}c@{}} $M'=G\diag(\Xi',\Xi)G^{\tau}$ 
\\  \end{tabular}
\\
\hline
{Autonomous dual system} & $\dot x' = A^\tau x'$ & $\dot x' = -A^{\tau}x'$\\
\hline
Trial basis $\Phi$ & 
$\Phi = \bar \Phi$ & 
\begin{tabular}{@{}c@{}}
$\Phi =
G\diag(\bar \Phi,\bar \Psi)$
\end{tabular}\\
\hline
Test basis $\Psi$ & 
$\Psi = \bar \Psi$ & $\Psi = G^{-\tau}\diag(\bar \Psi,\bar \Phi)$\\
\hline
Structure & 
\begin{tabular}{@{}c@{}}
Inner-product projection with \\
$M = \Xi$, $M' = \Xi'$,\\ $N=N' = \Sigma_1$
\end{tabular}
& 
\begin{tabular}{@{}c@{}}
Symplectic projection with \\
$J_{\Omega}= -J_{\Omega'}^{-1}=G^{-\tau} J_{2n} G^{-1}$,\\
$J_{\Pi}= J_{\Pi'} = J_{2k}$
\end{tabular}
\\
\hline
Balancing property & 
\begin{tabular}{@{}c@{}}
$\Phi\in O(M,\Sigma_1)$\\ $\Psi\in O(M',\Sigma_1)$ 
\end{tabular}
&
\begin{tabular}{@{}c@{}}
$\Phi\in O(M,\diag(\Sigma_1,\Sigma_1))$\\ $\Psi\in
O(M',\diag(\Sigma_1,\Sigma_1))$
\end{tabular}
\\
\hline
\begin{tabular}{@{}c@{}}
Canonical choice for \\
energies $\Xi$ and $\Xi'$ 
\end{tabular}
& $\Xi = W_o$, $\Xi' = W_c$ & $\Xi = \Xi' = \beta$\\
\hline
Stability preserved & 
\begin{tabular}{@{}c@{}}
Asymptotic stability if \\
$A^\tau \Xi + \Xi A \prec 0$ or\\
$A \Xi' + \Xi' A^\tau \prec 0$
\end{tabular}
& \begin{tabular}{@{}c@{}}
Pure marginal stability if\\
$A=JL$ with $J=-J_\Omega^{-1}$\\
 and $L\in \SPD{2n}$
\end{tabular}\\
\hline
\end{tabular}
\end{center}
\end{table}

\begin{definition}[Inner-product balancing]\label{def:innerprodbalancing}
Given any $\Xi\in\SPD{n}$ and $\Xi' \in\SPD{n}$, the trial and test basis matrices
characterizing an inner-product balancing 
correspond to 
\begin{equation}\label{eq:basesGenBT}
\Phi = SV_1\Sigma_1^{-1/2}\quad\text{and}\quad \Psi = R U_1\Sigma_1^{-1/2},
\end{equation}
respectively, where $\Xi = RR^\tau$, $\Xi' = S S^\tau$, and $R^\tau S = U\Sigma V^\tau$
is the singular value decomposition.  Here, we have defined $U =
[U_1\ U_2]$, $\Sigma = \diag(\Sigma_1,\Sigma_2)$, and $V = [V_1\
V_2]$, where $U_1,V_1\in O(I_n, I_k)$ and $\Sigma_1 =
\diag(\sigma_1,\ldots,\sigma_k)$ contains the $k$ largest singular values of $R^\tau S$.
\end{definition}
We now show that inner-product balancing leads to an inner-product projection.

\begin{lemma}\label{lem:GBTinnerProd}
An  inner-product balancing characterized by the test and trial basis matrices $(\Psi, \Phi)$  
 with $\Xi, \Xi' \in \SPD{n}$ has the following properties:
\begin{enumerate}[label=(\alph*)]
\setlength{\itemsep}{5pt}
\item The basis matrices $(\Psi, \Phi)$ correspond to an inner-product
projection performed on an LTI  system $(A,B,C)$ with $M=\Xi$ and $N = \Sigma_1$.
\item  The  basis matrices $(\Phi, \Psi)$ correspond to an inner-product projection performed on the dual system $(A^\tau, C^\tau, B^\tau)$ with $M'=\Xi'$ and $N' = \Sigma_1$.
\item The  basis matrices $(\Psi, \Phi)$ balance $\Xi$ and $\Xi'$, i.e., $\Phi \in O(\Xi,\Sigma_1)$ and $\Psi \in O(\Xi',\Sigma_1)$.
\end{enumerate}
\end{lemma}

\begin{proof}
To prove (a), we verify that $\Phi\in O(\Xi, \Sigma_1)$ and $\Psi=\Xi \Phi \Sigma_1^{-1}$, as 
\begin{gather*}
\Phi^\tau\Xi\Phi =   \left(  SV_1\Sigma_1^{-1/2}\right)^\tau \left (RR^\tau \right) \left( SV_1\Sigma_1^{-1/2} \right )=
 \left( \left (R^\tau S \right )V_1\Sigma_1^{-1/2}\right)^\tau  \left( \left(  R^\tau S\right )V_1\Sigma_1^{-1/2} \right)=\Sigma_1,\\
 \Xi\Phi\Sigma_1^{-1} = \left( R R^\tau\right) \left( SV_1\Sigma_1^{-1/2} \right) \Sigma_1^{-1} = R U_1\Sigma_1^{-1/2}=\Psi.
\end{gather*}
 Thus, the conditions for an inner-product projection are satisfied; note that $\Psi^\tau\Phi = I_k$.
To prove (b), recall from Remark \ref{rmk:dual} that the test basis of the dual system
corresponds to
$\Phi$, while the trial basis corresponds to $\Psi$. Thus, we aim to verify
that $\Psi\in O(\Xi',\Sigma_1)$ and $\Phi = \Xi'\Psi\Sigma_1^{-1}$, which can
be done similarly to the steps above.
Finally, (c) holds because $\Phi \in O(\Xi,\Sigma_1)$ and $\Psi \in O(\Xi',\Sigma_1)$.
\hfill
\end{proof}

\noindent We now show that any inner-product projection corresponds to a
particular balancing.
\begin{lemma}\label{lem:innerProdGBT}
If test and trial basis matrices $(\Psi, \Phi)$ characterize  an inner-product
projection with $M\in\SPD{n}$ and $N\in\SPD{n}$, then there exists
$M'\in \SPD{n}$  such that 
$\Phi \in O(M,N)$ and $\Psi\in O(M',N)$.
Further, there exists a realization of that reduced-order model that
corresponds to an inner-product balancing with
$\Xi =M$, $\Xi' = \hat M'$, and $\Sigma_1 = \diag(\lambda_1,\ldots,\lambda_k)$
characterized by basis matrices 
$(\hat \Psi, \hat\Phi)$
that satisfy
$\hat\Phi\in O(M,\Lambda)$, $\text{Ran}(\Phi) = \text{Ran}(\hat \Phi)$, $\hat
\Psi \in O(\hat M', \Lambda)$,
$\text{Ran}(\Psi) = \text{Ran}(\hat \Psi)$, and $\Psi^\tau\Phi =
\hat\Psi^\tau\hat\Phi = I_k$.
\end{lemma}
\begin{proof}
From the definition of an inner-product projection, we have
$\Phi\in O(M,N)$ and $\Psi = M\Phi N^{-1}$. Using
$M'=\Psi_\mathrm{ext}^{-\tau}N_\mathrm{ext}\Psi_\mathrm{ext}^{-1}$, we have 
$
[\Psi\ \bar\Psi]^\tau M' [\Psi\ \bar\Psi] = 
\diag(N,\bar N) 
$, whose (1,1) block gives
$\Psi\in O(M',N)$.
Similarly, $\Phi^\tau M\Phi = N$ implies $U^\tau\Phi^\tau M\Phi U = \Lambda$,
thus $\Phi\in O(M,N)$ implies $\hat\Phi\in O(M,\Lambda)$
with 
$\hat \Phi = \Phi U$.
Using
$\hat M'=\hat \Psi_\mathrm{ext}^{-\tau}\hat N_\mathrm{ext}\hat \Psi_\mathrm{ext}^{-1}$, we have $
[\hat \Psi\ \bar\Psi]^\tau \hat M' [\hat \Psi\ \bar\Psi] = 
\diag(\Lambda,\bar N) 
$, whose (1,1) block gives
$\hat \Psi\in O(M',\Lambda)$. Now, to ensure $(\hat\Psi,\hat\Phi)$ corresponds
to an inner-product projection with $M$ and $\Lambda$, we set $\hat \Psi = M\hat
\Phi\Lambda^{-1} = M\Phi
U\Lambda^{-1}=M\Phi N^{-1}N U\Lambda^{-1} = \Psi N U\Lambda^{-1}$.
Noting that $\text{Ran}(\Phi) = \text{Ran}(\hat \Phi)$ and 
$\text{Ran}(\Psi) = \text{Ran}(\hat \Psi)$ as well as $\Psi^\tau\Phi =
\hat\Psi^\tau\hat\Phi = I_k$, we conclude that basis
matrices $(\Psi,\Phi)$ and $(\hat \Psi,\hat \Phi)$ yield different
realizations of the same reduced-order model.
\hfill
\end{proof}

\medskip

\begin{corollary}\label{cor:innerproductbalancestability}
A inner-product-balancing reduced-order model preserves
asymptotic stability if $\Theta=\Xi$ satisfies the Lyapunov inequality \eqref{eq:equiLya4},
or if $\Theta'=\Xi'$ satisfies the  dual Lyapunov inequality \eqref{eq:dualinequiLya}.
\end{corollary}
\begin{proof}
The result follows directly from Lemmas \ref{lem:GBTinnerProd} (a) and 
 \ref{thm:stabpresadj}, as Lemma \ref{lem:GBTinnerProd} (b)  and 
 \ref{thm:stabpresadj} for the dual system.
\hfill
\end{proof}

We note that many existing model-reduction methods correspond to an
inner-product balancing; these methods are reported in Table
\ref{tab:comp_inner}. Appendix \ref{sec:geoexisting} discusses these model-reduction
methods in more detail.

\begin{table}
\begin{center}
\caption{Comparison of different model-reduction methods with
inner-product-balancing structure defined by $\Xi$ and $\Xi'$. In all cases,
$N = \Sigma_1$ is defined by the inner-product balancing. The remaining
quantities are defined in Appendix \ref{sec:geoexisting}.
}
\label{tab:comp_inner}
 {
\begin{tabular}{|c|c|c|c|c|c|}
\hline
& POD--Galerkin  & \begin{tabular}{@{}c@{}}   Balanced\\ truncation
\end{tabular}  & \begin{tabular}{@{}c@{}}   Balanced\\ POD \end{tabular} & SRSB & 
\begin{tabular}{@{}c@{}}  Proposed\\
inner-product\\ projection   \end{tabular}
 \\
\hline
\hline
$M = \Xi$  & $XX^\tau$ & $ W_o$ & $ \hat W_o$ & $W_o^\mu$ & $ \Theta$ satisfying \eqref{eq:equiLya4}
\\
\hline
$M' = \Xi'$ & $XX^\tau$ & 
$W_c$
&
$\hat W_c$
& 
$W_c^\mu$
&
$\Theta'$
satisfying \eqref{eq:dualinequiLya}\\
\hline
 \begin{tabular}{@{}c@{}}   Stability preservation?
 \end{tabular}    &   No  &  Yes\footnotemark  & No & No &
  Yes \\
\hline
 \end{tabular}}
\end{center}
\end{table}
\footnotetext{With the current framework, we can
only show that balanced truncation preserves marginal stability, as the
right-hand-side matrices in the Lyapunov equations are symmetric positive
\textit{semidefinite}. However,
additional analyses based on controllability and observability demonstrate that
balanced truncation does preserve asymptotic stability (e.g.,
Ref.~\cite[pp.~213--215]{AntoulasAC:05a}).}

\subsection{Construction of basis matrices given subset of ingredients} \label{sec:existenceInner}
Lemma \ref{thm:stabpresadj} demonstrated that a ROM will preserve asymptotic
stability if it is constructed via inner-product projection with $M=\Theta$
a Lyapunov matrix satisfying
\eqref{eq:equiLya4}.  Unfortunately, as reported in Table
\ref{tab:comp_inner}, while
many typical model-reduction techniques associate with an inner-product
projection (and an inner-product balancing),
the associated operator $M$ does not often satisfy the Lyapunov inequality, which
precludes assurances of stability preservation (e.g., POD--Galerkin,  Balanced POD, and SRSB).

We propose three methods (including inner-product balancing)  for constructing
a stability-preserving inner-product projection satisfying the conditions of
Lemma \ref{thm:stabpresadj}.  Table \ref{tab:two_case} summarizes these
methods; this corresponds to Steps 6--8 in Algorithm \ref{alg:overall}.  
Methods 2 and 3 assume that we are given a subset of the required ingredients,
which can be computed by any technique.  For example, the trial basis can be
computed by POD, balanced POD, or rational approximation; the metric $\Theta$
can be obtained by solving Lyapunov equation \eqref{eq:semiequiLya} with
$Q\in\SPD{n}$ but otherwise arbitrary. Thus, these methods can be viewed as
stabilization techniques applied to the provided inputs.

\begin{table}
\begin{center}
\caption{Algorithms for constructing an inner-product projection that ensure the
conditions of Definition \ref{def:inner-product-proj}. } 
\label{tab:two_case}
 {%
\begin{tabular}{|c|c|c|c|}
\hline
&
Method 1 (inner-product balancing)
&Method 2 & Method 3\\
\hline
\hline
 {Input} 
 &
 \begin{tabular} {@{}l@{}}   
     $\Xi,\Xi'\in\SPD{n}$ with \\
		$\Xi = \Theta$ satisfying \eqref{eq:equiLya4} or\\
		$\Xi' = \Theta'$ satisfying \eqref{eq:dualinequiLya}
     \end{tabular} 
 &
 \begin{tabular} {@{}l@{}}   
    $\Phi\in\RRstar{n\times k}$,\\
     $\Theta$ satisfying \eqref{eq:equiLya4}
     \end{tabular} 
    & 
    \begin{tabular} {@{}l@{}}   
    $\Phi_0\in O(M_0,  N_0)$, \\
    $N_0, N\in \SPD{k}$,\\
      $M_0 \in \SPD{n}$,\\
      $\Theta$ satisfying \eqref{eq:equiLya4}
     \end{tabular}   
        \\
\hline
 {Output} &  \begin{tabular} {@{}l@{}}   $M\in \SPD{n}$, $N\in \SPD{k}$, \\ $\Phi\in O(M,N)$,   $\Psi \in O(M', N)$  \end{tabular}  
 &\begin{tabular} {@{}l@{}} $M\in\SPD{n}$, $N\in \SPD{k}$,\\ $\Psi\in
 \RRstar{n\times k}$ \end{tabular} &
 \begin{tabular} {@{}l@{}}
 $M\in\SPD{n}$,\\
 $\Phi\in O(M, N)$, $\Psi \in \RRstar{n\times k}$   
\end{tabular} 
 \\
\hline
{Algorithm}  
&
  \begin{tabular} {@{}l@{}}   
    1.  Compute symmetric factorization\\
     \quad $\Xi=R R^{\tau}$, $\Xi'=S S^{\tau}$ \\
    2. Compute SVD $R^\tau S =  U \Sigma V^\tau$ \\
    3. $\bar \Phi = S V_1 \Sigma_1 ^{-1/2}$ \\
    4. $\bar \Psi = R U_1\Sigma_1 ^{-1/2}$\\
    5. $M = \Xi$, $M' = \Xi'$, $N = \Sigma_1$
     \end{tabular}
&
    \begin{tabular} {@{}l@{}}   
    1. $M=\Theta$ \\
    2. $N= \Phi^{\tau} M \Phi$ \\
    3. $\Psi =  M \Phi N^{-1}$
     \end{tabular} 
     &
     \begin{tabular} {@{}l@{}}   
    1. Set $M=\Theta$ \\
    2. Construct  $G  \in  O(M, M_0)$\\
    3. Construct  $\tilde G  \in  O(N, N_0)$ \\
    4. $\Phi =G \Phi_0 \tilde G^{-1}$ \\
    5. $\Psi =  M \Phi N^{-1}$
     \end{tabular}\\
\hline
\end{tabular}}
\end{center}
\end{table}

Method 2 constructs a stability-preserving inner-product projection starting
with any arbitrarily chosen trial basis matrix $\Phi\in \RRstar{n\times
k}$ and a Lyapunov matrix $\Theta$ satisfying
\eqref{eq:equiLya4}.  
Method 3 constructs a stability-preserving inner-product
projection starting with a basis $\Phi_0\in O(M_0, N_0)$, where $M_0$ might not
satisfy Lyapunov inequality \eqref{eq:equiLya4}, 
and a Lyapunov matrix $\Theta$ satisfying
\eqref{eq:equiLya4}. For simplicity, we can choose 
$G = M^{-1/2}M_0^{1/2}$ and  $\tilde G = N^{-1/2}N_0^{1/2}$. 
Lemma \ref{lem:sigma_innerpro} demonstrates that we can compute the trial basis matrix
in this context as $\Phi= G \Phi_0^{\tau} \tilde G^{-1}\in O(M,N)$, which constitutes
step 4 of the algorithm.

\begin{lemma}\label{lem:sigma_innerpro}
Let $\lift: (\mathbb{R}^{k}, N)\to (\mathbb{R}^{n}, M),x\mapsto \Phi x$ denote an
inner-product lift with $\Phi \in  O(M, N)$.
 Let $g: (\mathbb{R}^{n}, M_0) \to (\mathbb{R}^{n}, M),x\mapsto G x$ and
 $\tilde g: (\mathbb{R}^{k}, N_0) \to (\mathbb{R}^{k}, N),x\mapsto \tilde G x$
 represent (invertible) inner-product transformations, i.e., $G \in  O(M, M_0)\subseteq \RRstar{n\times n}$ and
 $\tilde G \in  O(N, N_0)\subseteq \RRstar{k\times k}$, respectively. Then, there exists a unique
 inner-product lift $\lift_0:  (\mathbb{R}^{k}, N_0)\to (\mathbb{R}^{n},
 M_0),x\mapsto\Phi_0 x$ with $\Phi_0 \in  O(M_0, N_0)$, such that the following diagram commutes:
\begin{center}
\begin{tikzpicture}[>=triangle 45]
\matrix(m)[matrix of math nodes,column sep={90pt,between origins},row
sep={50pt,between origins},nodes={asymmetrical rectangle}] (s)
{
|[name=M]| (\mathbb{R}^{n}, M_0)   &|[name=M1]| (\mathbb{R}^{n}, M) \\
|[name=N]| (\mathbb{R}^{k}, N_0)   &|[name=N1]| (\mathbb{R}^{k},  N)\\
};
\path[->] (M) edge node[auto] {\(g\)} (M1)
          (N) edge node[left] {\(\lift_0\)} (M)
          (N1) edge node[auto] {\(\tilde g^{-1}\)} (N)
          (N1) edge node[right] {\(\lift\)} (M1);
\end{tikzpicture}
\end{center}
Equivalently,  for all $z \in (\mathbb{R}^k, N)$,
\begin{equation}
\lift(z)=g(\lift_0(\tilde g^{-1}(z)))
\end{equation}
and $\Phi =G \Phi_0 \tilde G^{-1}$ in matrix representation.
\end{lemma}
\begin{proof}
Because  $G\in O(M, M_0)$,  we have $G^{\tau} M G=M_0$. It follows that $G^{-\tau}M_0 G^{-1}=M$.  By the same argument, $\tilde G\in O(N, N_0)$ implies that $\tilde G^{\tau} N \tilde G= N_0$.
Because $\Phi \in O(M, N)$, we have $\Phi^{\tau} M \Phi=N$. Because $g$ is invertible, we can define $\lift_0:  (\mathbb{R}^{k}, N_0)\to (\mathbb{R}^{n},  M_0)$  by $g^{-1} \circ \lift \circ \tilde g$ with matrix representation $\Phi_0= G^{-1} \Phi \tilde G$.  It follows that
$\Phi_0^{\tau} M_0 \Phi_0 =   \tilde G^{\tau} \Phi^{\tau} (G^{-\tau}  M_0
G^{-1}) \Phi \tilde G= \tilde G^{\tau} (\Phi^{\tau} M \Phi) \tilde G=\tilde
G^{\tau} N \tilde G=N_0$.
The last equation implies that $\Phi_0\in O(M_0, N_0)$. Finally, if $\lift_0$
satisfies $\lift =g\circ \lift_0  \circ \tilde g^{-1}$, $\lift_0$ is uniquely
determined by $\lift_0=g^{-1} \circ \lift \circ \tilde g$.
\hfill
\end{proof}

\noindent We now show that if the original trial basis matrix $\Phi_0$
exhibits a POD-like optimality property, then $\Phi$ computed by Method 3 in
Table \ref{tab:two_case} will inherit a related
optimality property. Given a set of snapshots $\{x_i\}_{i=1}^N$ with 
$x_i \in
(\mathbb{R}^{n}, M)$, we define the
projection error of the ensemble in the $M$-induced norm by
$\sum_{i=1}^N \left \| {x_i-  \Phi \Psi^\tau x_i } \right \|_{M} ^2 = \sum_{i=1}^N \left \|
{x_i-  \Phi N^{-1} \Phi^\tau M x_i } \right \|_{M} ^2$, where we have used
$\Psi= M \Phi N^{-1}$.

\begin{theorem} Let $M_0$, $M$, $N_0$, $N$, $G$, $\tilde G$, $\Phi_0$ and $\Phi$ be as defined in Lemma \ref{lem:sigma_innerpro}.
If $\Phi_0$ minimizes the projection of the snapshot ensemble $\{y_i\}_{i=1}^N$  with $y_i \in (\mathbb{R}^{n}, M_0)$, i.e.,
\begin{equation}\label{eq:PODoptimality}
\Phi_0=\underset{V_0\in O(M_0,N_0)}{\arg\min}\sum_{i=1}^N\| y_i -
V_0 N_0^{-1}V_0^\tau M_0 y_i\|_{M_0}^2,
\end{equation}
then $\Phi=G \Phi_0 \tilde G^{-1}$ minimizes the projection of the snapshot
ensemble $\{x_i\}_{i=1}^N$  with $x_i= G y_i \in (\mathbb{R}^{n}, M)$, i.e.,
\begin{equation}\label{eq:PODoptimality_phi}
\Phi=\underset{V\in O(M,N)}{\arg\min}\sum_{i=1}^N\| x_i - VN^{-1}V^\tau M x_i\|_{M}^2.
\end{equation}
Moreover, the cost function in \eqref{eq:PODoptimality} and \eqref{eq:PODoptimality_phi} achieves the same minimal value.
\end{theorem}
\begin{proof}
By Lemma \ref{lem:sigma_innerpro}, for any $V\in O(M,N)$, there exists a unique $V_0\in O(M_0, N_0)$ such that $V=G V_0 \tilde G^{-1}$. For any $i\in \{1, \ldots, N\}$, we have
\begin{align*}
\|x_i - VN^{-1}V^\tau M x_i\|_{M}^2 &= \| x_i -  (G V_0 \tilde G^{-1}) N^{-1} (G V_0 \tilde G^{-1})^\tau M x_i\|_{M}^2  \quad \quad  \quad \quad (V=G V_0 \tilde G^{-1})\\
&=\|G y_i  -  G V_0 (\tilde G^{-1} N^{-1} G^{-\tau}) V_0^\tau  (G^{\tau} M G) y_i  \|_{M}^2   \        \quad \quad (x_i = G y_i )   \\ 
&= \| G y_i  -  G V_0  N_0^{-1} V_0^\tau  M_0 y_i  \|_{M}^2    \ \quad \quad \quad \quad \quad \quad   \quad \quad   \quad (\tilde G ^\tau N \tilde G =N_0, \ G^{\tau} M G=M_0)\\
&= \|  y_i  -   V_0  N_0^{-1} V_0^\tau  M_0 y_i  \|_{M_0}^2.   \ \  \quad \quad \quad \quad \quad \quad \quad   \quad \quad   \quad (G^{\tau} M G=M_0)
\end{align*}
Then, the  cost function in \eqref{eq:PODoptimality} and
\eqref{eq:PODoptimality_phi} have the same  value when $V=G V_0 \tilde
G^{-1}$. Thus, if $\Phi_0$ is given by \eqref{eq:PODoptimality}, then $\Phi =G
\Phi_0 \tilde G^{-1}$ is the optimal value in \eqref{eq:PODoptimality_phi}.
Moreover, two cost functions achieve the same minimal value. 
\hfill
\end{proof}

\noindent
We note that (typical) POD satisfies optimality
property \eqref{eq:PODoptimality} with $M_0 = I_n$, $N_0 = I_k$, and
$\{y_i\}_{i=1}^N$ corresponding to snapshots of the system state, while
balanced POD  \cite{willcox2002bmr,RowleyCW:05a} satisfies this
property with $M_0=\hat W_o$ and $N_0 = \Sigma_1$, and $y_i$,
$i\in \{1,\ldots, N\}$ corresponding to snapshots arising from an impulse response. Because $\Phi$  constructed by Method 3 in
Table \ref{tab:two_case} satisfies $\Phi=G \Phi_0 \tilde G^{-1}$, Theorem
\ref{eq:PODoptimality} implies that $\Phi$ inherits the optimality to
minimize the projection error.

\section{Reduction of pure marginally stable subsystems}\label{sec:marstab}

This section focuses on reducing the pure marginally stable subsystem $\dot
x_m =A_m x_m$.  While the inner-product-projection approach could be applied
to the marginally stable subsystem if it has Lyapunov structure (i.e., if
\eqref{eq:inequiLya}--\eqref{eq:semiequiLya} hold), not all marginally stable
systems exhibit this structure; further, such a reduction would not guarantee
the poles remain nonzero and purely imaginary. Instead, we pursue an approach
that is valid for all pure marginally stable subsystems. It is based on the
key observation that all pure marginally stable systems are equivalent to a
generalized Hamiltonian system.

Section \ref{sec:background} introduces LTI Hamiltonian systems and
demonstrates that the marginally stable  subsystem has symplectic structure
(Theorem \ref{thm:matrix_similar_Hamil}). 
Subsequently, Section
\ref{sec:symproj_spa} introduces 
symplectic projection, Section
\ref{sec:sym_proj_dyn} demonstrates that a model-reduction method based on
symplectic projection preserves symplectic structure of generalized LTI
Hamiltonian systems and thus preserves pure marginal stability,
Section \ref{sec:symplecticbalancing} presents the symplectic-balancing framework,
and Section \ref{sec:PSD} describes methods for constructing the basis matrices
that lead to a symplectic projection given a subset of the required
ingredients. For notational simplicity, we omit the subscript $m$
throughout this section.

\subsection{Pure marginally stable systems}\label{sec:background}
We begin by defining pure marginal stability.
\begin{definition}[Pure marginal stability]
Linear system \eqref{eq:linear_control_sys}  is 
\emph{pure marginally stable},  if the system matrix $A$ is nonsingular and  diagonalizable, and has  a purely imaginary
spectrum.
\end{definition}

\noindent
If $A$ is a $2n\times 2n$ matrix, pure marginal stability means $A\in \generalizedham{2n}$.

We next  introduce the concept of symplectic spaces, and subsequently
introduce the LTI Hamiltonian and generalized LTI Hamiltonian equations.
Then, Theorem \ref{thm:matrix_similar_Hamil} proves the key result: any pure
marginally stable system  is a
generalized Hamiltonian system.

Let $\mathbb{V}\cong \mathbb{R}^{2n}$ denote a vector space. A symplectic form
$\Omega: \mathbb{V}\times \mathbb{V} \to \mathbb{R}$   is  a skew-symmetric,
nondegenerate,  bilinear function on the vector space $\mathbb{V}$.
The pair $(\mathbb{V}, \Omega)$ is called a symplectic vector space. Assigning
a symplectic form $\Omega$ to $\mathbb{V}$ is referred to as equipping
$\mathbb{V}$ with symplectic structure.

By choosing canonical coordinates on $\mathbb{V}$, the symplectic vector space
can be represented by $(\mathbb{R}^{2n}, J_{2n})$, where $J_{2n} \in
\{0, \pm1\}^{2n\times 2n}$ is a Poisson matrix defined as
\begin{equation*}
 {J_{2n}} \defeq \begin{bmatrix}
  0_n & {I_n} \\
 {-I_n} & 0_n \\
 \end{bmatrix}
\end{equation*}
that satisfies $J_{2n}J_{2n}^{\tau}=J_{2n}^{\tau}J_{2n}=I_{2n}$, and $J_{2n}J_{2n}=J_{2n}^{\tau}J_{2n}^{\tau}=-I_{2n}$. 
 The symplectic form
 $\Omega$ can be
 represented by the Poisson matrix as
\begin{equation*}
\Omega(\hat x_1, \hat x_2)=x_1^{\tau} J_{2n} x_2, \quad
\forall x_1,x_2 \in \RR{2n},
\end{equation*}
 where (as before) the operator $\hat\cdot$ provides the representation of an element of a
vector space from its coordinates, i.e., $\hat x \in\mathbb {V}$, $\forall
x\in\RR{2n}$.

 \begin{definition}[LTI Hamiltonian system]
An LTI system $(A,B,C)$ is an LTI Hamiltonian system if its corresponding
autonomous
system is given by
\begin{equation}\label{eq:LTI_Ham_can}
\dot x=J_{2n} \nabla_x H_0(x)= J_{2n} L_0x,
\end{equation}
where $L_0\in \RRstar{2n\times 2n}$ is symmetric and 
defines the 
(quadratic)
Hamiltonian 
\begin{equation}\label{def:Linear_Ham}
H_0:\RR{2n}\rightarrow \RR{},\ x\mapsto {\frac{1} {2}}{x^{\tau}}L_0x.
\end{equation}
 \end{definition}
 \begin{definition}[Hamiltonian matrix]
A Hamiltonian matrix is given by 
\begin{equation}\label{def:Hamilmatrix_can}
A_0= J_{2n} L_0\in \RRstar{2n \times 2n},
\end{equation} 
where $L_0\in \RRstar{2n\times 2n}$ is symmetric.
 \end{definition}

\noindent Thus, the $A$ matrix characterizing an LTI Hamiltonian system
$(A,B,C)$ is a Hamiltonian matrix.
\medskip
\begin{lemma}\label{lemma:hamil}
$A_0 \in \RRstar{2n \times 2n}$ is a Hamiltonian matrix if and only if it
satisfies
\begin{equation}\label{eq:sym_eqn}
 A_0^{\tau} J_{2n}+ J_{2n} A_0=0.
 \end{equation}
 \vspace*{-5mm}
\end{lemma}
\begin{proof}
Suppose the matrix $A_0 \in \RRstar{2n \times 2n}$ is Hamiltonian. Substituting $A_0$ with $J_{2n} L_0$, 
\eqref{eq:sym_eqn} holds for any symmetric $L_0$. Conversely, suppose
\eqref{eq:sym_eqn} holds. This implies that 
$ J_{2n}^\tau A_0$
is symmetric.  Let  $L_0=J_{2n}^\tau A_0$. Then, $A_0=J_{2n} L_0$, which  is a Hamiltonian matrix. 
\hfill
\end{proof}

\noindent We note that the Hamiltonian is constant in time, as  
\begin{equation*}
\frac{d}{dt} H_0(x(t))= \frac{d}{dt}  \left ( \frac{1}{2} x^{\tau} L_0 x \right
)=x^{\tau} L_0 \dot x= x^{\tau}L_0 A_0x= x^{\tau} L_0 J_{2n} L_0 x=0.
\end{equation*}

More generally, if non-canonical coordinates are chosen,  the symplectic
vector space $(\mathbb{V}, \Omega)$ can be represented by $(\mathbb{R}^{2n},
J_\Omega)$, where $J_\Omega \in \skewsymmetric{2n}$. Then, the symplectic form can be represented
as
\begin{equation*}
\Omega(\hat x_1, \hat x_2)=x_1^{\tau} J_\Omega x_2, \quad
\forall x_1,x_2 \in \RR{2n}.
\end{equation*}
\begin{definition}[Generalized LTI Hamiltonian system]
An LTI system $(A,B,C)$ is a generalized LTI Hamiltonian system if its
corresponding autonomous
system is given by
\begin{equation}\label{eq:LTI_Ham}
\dot x= J \nabla_x H(x)= J L x,
\end{equation}
where $J\in \skewsymmetric{2n}$ and $L \in \RRstar{2n\times 2n}$ is symmetric. The matrix $L$ defines the (quadratic)
generalized Hamiltonian
\begin{equation}\label{def:Linear_HamGen}
H:\RR{2n}\rightarrow \RR{},\ x\mapsto {\frac{1} {2}}{x^{\tau}}Lx.
\end{equation}
\end{definition}

 \begin{definition}[Generalized Hamiltonian matrix]
A generalized Hamiltonian matrix is given by 
\begin{equation}\label{def:A_matrix_gen_Ham}
A=J L\in\RRstar{2n\times 2n},
\end{equation} 
where $J \in\skewsymmetric{2n}$ and $L\in \RRstar{2n\times 2n}$ is symmetric.
 \end{definition}

\noindent Thus, the $A$ matrix characterizing a generalized LTI Hamiltonian system
$(A,B,C)$ is a generalized Hamiltonian matrix.

\begin{lemma}\label{lemma:hamil_generalized}
$A \in \RRstar{2n \times 2n}$ is a generalized Hamiltonian matrix if and only if it
satisfies
\begin{equation}\label{eq:condition_gen_Ham}
 A^{\tau} J_{\Omega}+ J_{\Omega} A=0
\end{equation} 
for some $J_\Omega\in\skewsymmetric{2n}$.
\end{lemma}

\begin{proof} 
Suppose $A$ is a generalized Hamiltonian matrix with $A=J L$, where $J \in \skewsymmetric{2n}$ and  $L\in \RRstar{2n \times 2n}$ is symmetric. Let $J_\Omega=-J^{-1}$. Then, we obtain $J_\Omega \in \skewsymmetric{2n}$ and $A^\tau J_\Omega + J_\Omega A = (JL)^\tau (-J^{-1})+ (-J^{-1}) (JL)=0$. Conversely, suppose \eqref{eq:condition_gen_Ham} holds with $J_\Omega\in \skewsymmetric{2n}$. Then, $-J_{\Omega} A$ is symmetric. Let  $L=-J_{\Omega} A$.  Because both $-J_{\Omega} $ and $A$ are nonsingular, so is $L$. Moreover, 
because $J_{\Omega}\in \skewsymmetric{2n}$, $-J_\Omega^{-1} \in \skewsymmetric{2n}$.
 Thus, $A=-J_\Omega^{-1} L$ is a generalized Hamiltonian matrix. 
\hfill 
\end{proof}


As with the Hamiltonian, the generalized Hamiltonian is constant in time, as  
\begin{equation*}
\frac{d}{dt} H(x(t))= \frac{d}{dt}  \left ( \frac{1}{2} x^{\tau} L x \right
)=x^{\tau} L \dot x= x^{\tau}L Ax= x^{\tau} L J L x=0
\end{equation*}
due to skew-symmetry of $J$.


We now derive the transformation between the coordinates defining the
Hamiltonian and generalized Hamiltonian systems.
\begin{lemma}\label{lem:skewsymTransform}
$J_\Omega\in \skewsymmetric{2n}$ if and only
if there exists a $G\in \RRstar{2n\times 2n}$ such that 
\begin{equation}\label{eq:transform_can}
G^{\tau} J_{\Omega} G=J_{2n}.
\end{equation}
\end{lemma}

\noindent The proof is provided in  \cite[Corollary 5.4.4]{EvesH:80a}.
Note that $\Sp(J_{\Omega}, J_{2n})$ is not empty for $J_\Omega\in \skewsymmetric{2n}$ by Lemma \ref{lem:skewsymTransform}. This set
represents the set of  (invertible, linear) \textit{symplectic transformations}
  $g: (\mathbb{R}^{2n}, J_{2n}) \to (\mathbb{R}^{2n},
 J_\Omega), \ y \mapsto G y $ with $G\in\Sp(J_{\Omega}, J_{2n})$
  such that 
\begin{equation*}
x_1^{\tau} J_{2n} x_2 = (G x_1)^{\tau} J_\Omega (G x_2),\quad \forall
x_1, x_2 \in  (\mathbb{R}^{2n},
 J_{2n})
.
\end{equation*}

\medskip

\begin{lemma}\label{lemma:hami_gen_can}
$A$ is a generalized
Hamiltonian matrix
if and only if it can be transformed into a Hamiltonian matrix $A_0$ by a
similarity transformation with a matrix $G\in\RRstar{2n\times 2n}$, i.e.,
\begin{equation}\label{eq:TAT}
A_0=G^{-1}A G.
\end{equation}
\end{lemma}
\begin{proof}
Assume $A=J L$ is a generalized Hamiltonian matrix, where $J \in \skewsymmetric{2n}$ and  $L\in \RRstar{2n\times 2n}$ is symmetric. Because $-J^{-1} \in \skewsymmetric{2n}$, Lemma \ref{lem:skewsymTransform} implies that there exists $G\in\Sp(J_{\Omega}, J_{2n})$ such that $G^\tau (-J^{-1}) G=J_{2n}$. It follows that
\begin{equation*}
G^{-1} AG=G^{-1} (J L) G = G^{-1} (G J_{2n} G^{\tau}) L G = J_{2n} (G^\tau L G) .
\end{equation*}
By setting $L_0=G^\tau L G$, $L_0$ is symmetric and nonsingular. Letting $A_0= J_{2n}  L_0$, the last expression implies that 
$G^{-1}A G =A_0$, which is a Hamiltonian matrix.

Conversely, suppose that $A_0 = G^{-1}AG$, where $A_0 = J_{2n} L_0$ is a Hamiltonian matrix and $G$ is nonsingular.
Substituting $A_0 = G^{-1} A G$ into \eqref{eq:sym_eqn} yields
\begin{equation*} 
(G^{-1} A G)^{\tau} J_{2n}+ J_{2n} (G^{-1} A G)=0.
 \end{equation*}
Left-multiplying by $G^{-\tau}$ and right-multiplying by $G^{-1}$ yields
\begin{equation*} 
A^{\tau} G^{-\tau} J_{2n} G^{-1}+ G^{-\tau} J_{2n}  G^{-1} A =0.
 \end{equation*}
Letting  $J_{\Omega} =G^{-\tau} J_{2n}G^{-1}$, we have $J_\Omega\in \skewsymmetric{2n}$. The above equation is equivalent
to \eqref{eq:condition_gen_Ham}, which implies that $A$ is a generalized
Hamiltonian matrix.
\hfill
\end{proof}
\medskip

By Lemma \ref{lemma:hami_gen_can}, an autonomous LTI system $\dot x =
Ax$ can be transformed into an LTI Hamiltonian system  \eqref{eq:LTI_Ham_can}
if and only if the autonomous LTI system is a generalized LTI Hamiltonian system
\eqref{eq:LTI_Ham},
 as substituting $x=G y $ in 
\eqref{eq:LTI_Ham} yields
\begin{equation*}
\dot y= G^{-1} \dot x= G^{-1} A x = G^{-1} A G y =A_0 y.
\end{equation*}

The next theorem shows that  any pure marginally stable system $\dot x= A
x$ with $A\in\generalizedham{2n}$ is a generalized LTI Hamiltonian system,
where $\generalizedham{n}$ denotes the set of real-valued $n\times n$
diagonalizable matrices with
nonzero purely imaginary eigenvalues.

\medskip

 \begin{theorem}\label{thm:matrix_similar_Hamil}
	The following conditions are equivalent:
	\begin{enumerate}[label=(\alph*)]
\setlength{\itemsep}{5pt}
\item $A\in \generalizedham{2n}$.
\item $A$  is a generalized  Hamiltonian matrix whose
corresponding generalized LTI Hamiltonian system is marginally stable.
\item There exists $G\in \RRstar{2n\times 2n}$ such that $G^{-1} A G= J_{2n}
L_0$, where $L_0=\diag(\beta, \beta)$, $\beta=\diag(\beta_1, \ldots,
\beta_n)$, and $\beta_1\ge \ldots \ge \beta_n>0$.  Further, we have
\begin{gather} 
\label{eq:generalizedHamPrimal}A=J L, \quad J =G J_{2n} G^{\tau},\quad  L=G^{-\tau}
L_0 G^{-1},
\end{gather} 
where $J\in \skewsymmetric{2n}$ and $L\in\SPD{2n}$.
\end{enumerate}
\end{theorem}
\begin{proof}
 $(b) \Rightarrow (a)$.
Recall from Lemma \ref{lemma:stability0} that the marginal-stability assumption is equivalent to assuming
that 
eigenvalues of $A$ have non-positive real parts
	and all Jordan
blocks corresponding to eigenvalues with zero real parts are $1 \times 1$.
Now, assume that $A$ is a generalized Hamiltonian matrix whose eigenvalues
have non-positive real parts. By Lemma \ref{lemma:hamil_generalized}, a generalized Hamiltonian matrix
$A$ satisfies $A^{\tau} J_\Omega + J_\Omega A =0$ for some
$J_\Omega\in\skewsymmetric{2n}$. It follows that
\begin{equation}\label{eq:sym_similarity}
A= J_\Omega^{-1}(-A^{\tau}) J_\Omega,
\end{equation}
i.e., $A$ is similar to $-A^{\tau}$. So, if $\lambda$ is an
eigenvalue of $A$, $\lambda$ is an eigenvalue of $-A^{\tau}$ and thus an eigenvalue
of $-A$. This implies that $-\lambda$ is also an eigenvalue
of $A$. Thus, the eigenvalues of $A$ would have positive real parts unless the
real part of $\lambda$ is zero, i.e., the eigenvalues of $A$ are purely imaginary. Due to
the marginal-stability assumption on 
$A$, every Jordan block for purely imaginary eigenvalues
must has dimension $1\times 1$. Therefore,  $A$ is diagonalizable and  has only
nonzero purely imaginary eigenvalues, i.e., $A\in \generalizedham{2n}$.

 $(a) \Rightarrow (c)$.
Assume $A\in \generalizedham{2n}$. Let $\lambda$ be an eigenvalue of $A$. Then $\lambda$ is a root of the
characteristic polynomial $\det(\lambda I_{2n}- A)=0$.  Because the matrix $A$
is a real matrix, the characteristic polynomial only contains real
coefficients of $\lambda$. Thus, if $i \beta_0$ with $\beta_0 \in \mathbb{R}$ is a
root of $\det(\lambda I_{2n}- A)=0$, so is $-i \beta_0$. Moreover, $i \beta_0$ and $-i
\beta_0$ must have the same algebraic multiplicity. It follows that $A$ contains
eigenvalues of the form $\{\pm i \beta_1,  \ldots, \pm i \beta_n \}$, where
$\beta_1\ge \ldots \ge \beta_n > 0$. Because the system matrix $A$ is assumed to
be diagonalizable, 
there exists a matrix $P_1\in \mathbb{C}_*^{2n \times 2n}$ such that 
\begin{equation}\label{eq:diagAP1}
A=P_1 \diag(i \beta_1, -i \beta_1, \ldots, i \beta_n, -i \beta_n) P_1^{-1}.
\end{equation}
Let $\beta=\diag(\beta_1, \ldots, \beta_n)$, it is straightforward to verify that the  matrix $J_{2n} \diag(\beta, \beta) $ also contains eigenvalues $\{\pm i \beta_1,
\ldots, \pm i \beta_n \}$ and is diagonalizable. Thus, there exists a matrix $P_2\in
\mathbb{C}_*^{2n \times 2n}$ such that 
\begin{equation} \label{eq:diagP2} 
J_{2n} \diag(\beta, \beta) =P_2 \diag(i \beta_1, -i \beta_1, \ldots, i \beta_{n}, -i \beta_{n})P_2^{-1}.
\end{equation}
With $P_3=P_1 P_2^{-1}\in \mathbb{C}_*^{2n \times 2n}$, we have
\begin{equation}\label{eq:P3trans}
P_3^{-1}A P_3= P_2 (P_1^{-1}  A P_1) P_2^{-1} 
=P_2 \diag(i \beta_1, -i \beta_1,  \ldots, i \beta_n -i  \beta_n) P_2^{-1}
= J_{2n} L_0 .
\end{equation}
where $L_0=\diag(\beta, \beta) \in\RRstar{2n}$ is symmetric . 
Equation \eqref{eq:P3trans}  implies that $A$ is similar to the Hamiltonian
matrix $J_{2n} L_0$ via a complex matrix $P_3$. Let $A_0=J_{2n} L_0$,
we can also rewrite \eqref{eq:P3trans} as 
\begin{equation*}
A P_3 = P_3 A_0.
\end{equation*}
Decomposing this matrix as $P_3=P_4 + i P_5$ with $P_4, P_5 \in \RR{2n\times
2n}$ and noting 
that both $A$ and $A_0$ are real matrices, the above equation implies that
 \begin{equation*}
 A P_4= P_4 A_0, \quad  A P_5= P_5 A_0.
\end{equation*}
 This implies that $A$ is similar to $A_0$ via a real matrix $P_4+\alpha P_5$
 for any $\alpha \in \mathbb{R}$. Because $\text{det}(P_3)=\text{det}(P_4+i P_5)\ne 0$,  $P(\alpha):=\text{det}(P_4+\alpha P_5)$ is a nonzero polynomial of $\alpha$ with degree no great than $2n$. Thus, the  equation $P(\alpha)=0$ contains $2n$ roots in $\mathbb{C}$ at most. Thus, we can choose $\alpha_0\in \mathbb{R}$ such that  $G=P_4+\alpha_0 P_5$  is invertible.	
 Thus, we obtain $G^{-1} A G =J_{2n}L_0$ with  $G \in  \RRstar{2n \times 2n}$ and $L_0=\diag(\beta, \beta)$. Setting $J=G J_{2n} G^{\tau}$ and $L=G^{-\tau} L_0 G^{-1}$, we have $J\in \skewsymmetric{2n}$ and $L\in \SPD{2n}$. It follows that $$A=G \left (J_{2n}L_0 \right)G^{-1}=(G J_{2n} G^\tau )(G^{-\tau} L_0 G^{-1})=J L.$$  
	 
 $(c) \Rightarrow (b)$. Suppose $G^{-1}AG=J_{2n}L_0 \in
\mathbb{R}^{2n \times 2n}$ with $G \in  \RRstar{2n \times 2n}$ and $L_0=\diag(\beta, \beta)$.
  Lemma \ref{lemma:hami_gen_can}, implies that $A$ is a generalized Hamiltonian matrix. Because   $J_{2n} \diag(\beta, \beta) $ contains eigenvalues  $\{\pm i \beta_1, \ldots, \pm i \beta_n \}$ and is diagonalizable,  so is $A$. Therefore, the corresponding system of $A$ is marginally stable.
  \hfill
\end{proof}

The part  $(a) \Rightarrow (c)$ is a constructive proof.   Algorithm \ref{alg:constr_cannon} lists the detailed procedure.

\begin{algorithm}
\caption{Transform $A\in \generalizedham{2n}$ into a canonical Hamiltonian matrix.} \label{alg:constr_cannon}
\begin{algorithmic}[1]
 \REQUIRE
$A\in \generalizedham{2n}$.
\ENSURE $G\in \RRstar{2n\times 2n}$  satisfying
$G^{-1} A G= J_{2n} L_0$, where $L_0=\diag(\beta, \beta)$ and 
$\beta=\diag(\beta_1, \ldots, \beta_n)$.
\STATE  Compute the eigenvalue decomposition \eqref{eq:diagAP1} of $A$ to obtain the eigenvalues $\{\pm i \beta_1, \ldots, \pm i \beta_n \}$ and the transformation matrix $P_1\in \mathbb{C}_*^{2n \times 2n}$. 
\STATE Construct the matrix $A_0=J_{2n} L_0$, where $L_0=\diag(\beta, \beta)$ with $\beta=\diag(\beta_1, \ldots, \beta_n)$.
\STATE Compute the eigenvalue decomposition \eqref{eq:diagP2} of $A_0$ to obtain the transformation matrix $P_2 \in \mathbb{C}_*^{2n \times 2n}$. 
\STATE Compute $P_3=P_1 P_2^{-1} \in \mathbb{C}_*^{2n \times 2n}$. 
\STATE Decompose $P_3=P_4 + i P_5$ with $P_4, P_5 \in \mathbb{R}^{2n \times
2n}$ and 
define $P(\alpha)\defeq	P_4 + \alpha P_5$.
\STATE $\alpha\leftarrow 0$
\WHILE{$\det(P(\alpha))=0$ and $\alpha<2n$}
\STATE $\alpha\leftarrow \alpha+1$.
\STATE Update  $P(\alpha)=P_4 + \alpha P_5$.
\ENDWHILE
\STATE $G=P(\alpha)$.
\end{algorithmic}
\end{algorithm}

Theorem \ref{thm:matrix_similar_Hamil} implies that performing model reduction
in a manner that preserves generalized Hamiltonian structure and marginal
stability  will ensure that the reduced-order
model retains pure marginal stability. We will accomplish this via symplectic
projection.

\begin{corollary}[Dual version of Theorem \ref{thm:matrix_similar_Hamil}]\label{cor:decomp_A_Adual}
If  any condition of Theorem \ref{thm:matrix_similar_Hamil} holds, then the
following conditions hold:
\begin{enumerate}[label=(\alph*)]
\setlength{\itemsep}{5pt}
\item $-A^\tau \in \generalizedham{2n}$.
\item $-A^\tau$  is a  generalized  Hamiltonian matrix whose
corresponding generalized LTI Hamiltonian system is marginally stable.
\item With $G$ and $L_0$ defined in Theorem \ref{thm:matrix_similar_Hamil}, we
have $G^\tau (-A^\tau) G^{-\tau}= J_{2n} L_0$. Moreover, we can have
\begin{gather}\label{eq:generalizedHamDual}-A^\tau=J' L',\quad J'=G^{-\tau} J_{2n}
G^{-1}= -J^{-1},\quad  L'=G L_0 G^{\tau},
\end{gather}
where $J'\in \skewsymmetric{2n}$ and $L'\in\SPD{2n}$.
\end{enumerate}
\end{corollary}

\begin{proof} 
Suppose $A\in  \generalizedham{2n}$. Then,  \eqref{eq:sym_similarity} holds for some $J_\Omega\in \skewsymmetric{2n}$, which implies that $A$ is similar to $-A^{\tau}$.
 Thus, $-A^\tau \in  \generalizedham{2n}$. By Theorem   \ref{thm:matrix_similar_Hamil}, (a) and (b) in this  corollary are equivalent. 

Using (c) in Theorem   \ref{thm:matrix_similar_Hamil}, i.e., $G^{-1} A G=J_{2n} L_0$ with $L_0=\diag(\beta, \beta) $, we have 
$$
G^{-1} A G=
\begin{bmatrix}
0 & \beta \\
-\beta & 0
\end{bmatrix}.
$$
It follows that 
$$
G^{\tau} (-A^\tau) G^{-\tau}=
\begin{bmatrix}
0 & \beta \\
-\beta & 0
\end{bmatrix} = J_{2n}L_0.
$$
Defining $J'=G^{-\tau} J_{2n}G^{-1}$ and $ L'=G L_0 G^{\tau}$, we have $J'\in \skewsymmetric{2n}$ and $L'\in\SPD{2n}$.
Moreover, the above equation yields $-A^\tau=G^{-\tau} (J_{2n}L_0) G^\tau=J' L'$. Finally, with  $J =G J_{2n} G^{\tau}$, we obtain $J'=-J^{-1}$.
\hfill
\end{proof}

\begin{remark}[Relationship with dual system: pure marginal
stability]\label{rmk:dualHam}
\noindent Any method proposed in this work for ensuring pure marginal stability of a
given (sub)system also ensures pure marginal stability of the associated
\textit{negative} dual
(sub)system $(-A^\tau, C^\tau,  B^\tau)$.  However, as before,  the
proposed methods for constructing a trial basis matrix $\Phi$ should be
applied to the negative dual system as a test basis matrix. Similarly, the proposed
methods for constructing a test basis matrix $\Psi$ should be applied to the
negative dual system as a trial basis matrix.
\end{remark}

\subsection{Symplectic projection of spaces}\label{sec:symproj_spa} Let
$(\mathbb{V}, \Omega)$ and $(\mathbb{W}, \Pi)$ be two symplectic vector spaces
with coordinate representations $(\mathbb{R}^{2n}, J_\Omega)$ and $(\mathbb{R}^{2k}, J_\Pi)$,
respectively, 
$\dim(\mathbb{V})=2n$, $\dim(\mathbb{W})=2k$, and $k \le n$.

\begin{definition}[Symplectic lift]
A {\emph{symplectic lift}} is a linear mapping $\lift: (\mathbb{W}, \Pi) \to (\mathbb{V}, \Omega)$ that preserves symplectic structure:
 \begin{equation} \label{def:symp_lift}
 \Pi(\hat z_1, \hat z_2)=\Omega(\lift (\hat z_1), \lift (\hat z_2)), \quad \forall \hat z_1, \hat z_2 \in \mathbb{W}.
 \end{equation}
\end{definition}

\begin{definition}[Symplectic projection]
Let $\lift: (\mathbb{W}, \Pi) \to (\mathbb{V}, \Omega)$ be a symplectic lift.  The adjoint of $\lift$ is the linear mapping $\projection:  (\mathbb{V}, \Omega) \to (\mathbb{W}, \Pi)$ satisfying
 \begin{equation}\label{def:sym_pro}
\Pi( \projection(\hat x), \hat z)=\Omega( \hat x,\lift(\hat z)), \quad \forall
\hat z\in \mathbb{W}, \ \hat x\in \mathbb{V}.
\end{equation}
 We say $\projection$ is the {\emph{symplectic projection}} induced by $\lift$.
\end{definition}
\medskip

As in the case of the inner-product lift and projection, the symplectic lift
and projection can be expressed in coordinate space as
\begin{align*}
\lift(\hat z) &\equiv \Phi z,\ \forall z \in \RR{2k},\quad\quad
\projection(\hat x)\equiv\Psi^\tau x,\ \forall x\in\RR{2n},
\end{align*}
respectively, where \eqref{def:symp_lift}--\eqref{def:sym_pro} imply
that $\Phi \in \RRstar{2n\times 2k}$  and $\Psi \in
\RRstar{2n\times 2k}$ satisfy
\begin{align}\label{eq:AJA}
\Phi^{\tau} J_\Omega \Phi &=J_\Pi\\
 \label{eq:JPhi_Trans}
\Psi J_\Pi  &=J_\Omega \Phi,
\end{align}
from which it follows that
\begin{equation}\label{def:Phi+}
\Psi= J_\Omega  \Phi J_\Pi^{-1}.
\end{equation}

When \eqref{eq:AJA} holds, we say $\Phi$  is a \emph{symplectic matrix}  with
respect to $J_\Omega$ and $J_\Pi$, which we denote by $\Phi\in
\Sp(J_\Omega,J_\Pi)$.
As in the inner-product projection case, it can be verified that
$\Psi^\tau$ is a left inverse of $\Phi$, as 
 \begin{equation}
 \Psi^{\tau} \Phi=(J_\Omega  \Phi J_\Pi^{-1})^{\tau} \Phi=J_\Pi^{-1} (\Phi^{\tau} J_\Omega \Phi )=J_\Pi^{-1}  J_\Pi= I_{2k},
 \end{equation}
which implies that $\projection \circ
\lift$ is the identity map on $\mathbb{W}$.


\subsection{Symplectic projection of dynamics}\label{sec:sym_proj_dyn} 
This section first defines symplectic projection of dynamics. We show that if the original system is a generalized Hamiltonian LTI system, then the reduced system constructed by symplectic projection is also a  generalized Hamiltonian LTI system.

\begin{definition}[Model reduction via symplectic projection]\label{def:sym_proj_dyn} 
A reduced-order model 
$(\tilde A, \tilde B, \tilde C)$ with $\tilde A=\Psi^{\tau}  A \Phi$, $\tilde B=\Psi^\tau B$, and $\tilde C=C\Phi$
is constructed by a symplectic projection if  $\Phi \in \Sp(J_\Omega, J_\Pi)$ and $\Psi= J_\Omega  \Phi J_\Pi^{-1}$, where $J_\Omega\in \skewsymmetric{2n}$ and $J_\Pi\in\skewsymmetric{2k}$. 
 \end{definition}
 
\medskip
\begin{lemma}\label{lem:symplec_preserv}
If the original LTI system $(A, B, C)$ is a generalized LTI Hamiltonian
system---i.e., $A = J L$ with $J \in\skewsymmetric{2n}$ and
$L\in\RRstar{2n\times 2n}$ is symmetric---and the reduced-order model
is constructed by symplectic projection with $J_\Omega=-J^{-1}$, then the reduced-order model $(\tilde A,
\tilde B, \tilde C)$ remains a generalized LTI Hamiltonian system, i.e., 
 $\tilde A=-J_{\Pi}^{-1}\tilde L$, where $J_\Pi\in \skewsymmetric{2k}$ and
 $\tilde L = \Phi^\tau L\Phi\in \RRstar{k\times k}$ is symmetric.
\end{lemma}

\begin{proof}
Because $A= -J_\Omega^{-1} L$ and $\Phi \in \Sp(J_\Omega,
J_\Pi)$, we have from  \eqref{def:Phi+} that
$$
\tilde A= \Psi^{\tau}A \Phi  = (J_\Pi^{-1} \Phi^{\tau} J_\Omega) (-J_\Omega^{-1} L) \Phi = -J_\Pi^{-1} (\Phi^{\tau} L  \Phi).
$$
Because $J_\Pi\in\skewsymmetric{2k}$, $-J_\Pi^{-1}\in \skewsymmetric{2k}$. Define $\tilde L=\Phi^{\tau} L \Phi\in\RRstar{2k\times 2k}$. Because $L$ is symmetric and nonsingular, so is $\tilde L$. Thus, $\tilde A=-J_{\Pi}^{-1}\tilde L$ is  a generalized Hamiltonian matrix.
\hfill
\end{proof}
\medskip


Recall that if $A\in \generalizedham{2n}$,  Theorem \ref{thm:matrix_similar_Hamil}  (c) implies that there exists $G\in \RRstar{2n \times 2n}$ such that 
$A=J L$,  $J=G J_{2n} G^{\tau} \in \skewsymmetric{2n}$, and $L=G^{-\tau} L_0 G^{-1}\in \SPD{2n}$, where  $L_0= \diag(\beta_1, \ldots, \beta_n, \beta_1, \ldots,  \beta_n)$ with $\beta_i > 0$.

\medskip
\begin{theorem} [{\rm  Preservation of pure marginal stability}] \label{thm:mainstab}
Suppose the original system $(A, B, C)$ is pure marginally stable, i.e.,  $A=JL\in
\generalizedham{2n}$ with $J\in \skewsymmetric{2n}$ and $L\in \SPD{2n}$.   Then the  reduced system  $(\tilde A, \tilde B, \tilde
C)$ constructed by  symplectic projection with $J_\Omega=-J^{-1}$ and any $J_\Pi \in \skewsymmetric{2k}$ remains pure marginally stable,  i.e., $\tilde A \in \generalizedham{2k}$.
\end{theorem}
\medskip
\begin{proof}
 Lemma \ref{lem:symplec_preserv}, the reduced system matrix $\tilde A$ constructed by  symplectic projection with $J_\Omega=-J^{-1}$ can be written as $\tilde A = -J_\Pi ^{-1} \tilde L$ with $\tilde L=\Phi^{\tau} L  \Phi$ and $\Phi\in \Sp{(J_\Omega, J_\Pi)}$. Because $L\in\SPD{2n}$, we have $\tilde L\in\SPD{2k}$.  

Let $\tilde H:z\mapsto\frac{1}{2} z^{\tau} \tilde L z$ denote  the Hamiltonian
function of the reduced system $\dot z =\tilde A z$. Because $\tilde
L\in\SPD{2k}$, there exists a $\delta>0$ such that $\tilde H(z)> \tilde
H(z_0)$ for all $\|z\|=\delta$, where  $z_0$ is the initial condition. Because
the reduced system is a generalized Hamiltonian system, the Hamiltonian
function satisfies $\tilde H(z(t))=\tilde H(z_0)$ for all $t\geq 0 $. Thus,
$\|z(t)\| < \delta$ for all $t\geq 0$, i.e., the reduced-order-model solution
is uniformly bounded.  Because the reduced system  is also linear, it is
marginally stable. 

Finally, since $\tilde A$ is a generalized  Hamiltonian matrix with marginal stability,  Theorem \ref{thm:matrix_similar_Hamil} implies $\tilde A\in \generalizedham{2k}$.
\hfill
\end{proof}

\begin{corollary}
Suppose $J_\Omega\in \skewsymmetric{2n}$ and the original system $(A, B, C)$ satisfies 
\begin{equation}\label{eq:sym_stable_condition}
-J_\Omega A \in \SPD{2n}.
\end{equation}
 Then the  reduced system  $(\tilde A, \tilde B, \tilde
C)$ constructed by  symplectic projection with $J_\Omega$ and any $J_\Pi \in \skewsymmetric{2k}$  is pure marginally stable,  i.e., $\tilde A \in \generalizedham{2k}$.
\end{corollary}
\begin{proof}
Define $L=-J_\Omega A$. Then, $L\in \SPD{2n}$ and $A=-J_\Omega^{-1}L$.  Theorem \ref{thm:mainstab} implies that $\tilde A \in \generalizedham{2k}$.
\hfill
\end{proof}


\subsection{Symplectic balancing}\label{sec:symplecticbalancing}

In analogue to Section \ref{sec:innerproductbalancing}, we now discuss a
symplectic-balancing approach that leverages symplectic structure. Recall that Table
\ref{tab:balancing} compares the proposed symplectic-balancing approach with
inner-product balancing.


\begin{definition}[Symplectic balancing]
Given any $\Xi, \Xi' \in\SPD{n}$, $J_\Omega \in \skewsymmetric{2n}$, and $G\in \Sp(J_{\Omega}, J_{2n})$, the trial
and test bases
characterizing a symplectic balancing
correspond to 
\begin{equation}\label{eq:subsetSBT}
\Phi = G\diag(\bar\Phi,\bar\Psi)\quad\text{and}\quad \Psi =
G^{-\tau}\diag(\bar\Psi,\bar\Phi),
\end{equation}
where basis matrices $(\bar\Psi,\bar\Phi)$ characterize an inner-product balancing
on matrices $\Xi$ and $\Xi'$, i.e.,
\begin{equation}\label{eq:basesGenSBT}
\bar\Phi = SV_1\Sigma_1^{-1/2}\quad \text{and}\quad\bar\Psi = R U_1\Sigma_1^{-1/2},
\end{equation}
where  quantities ($R, S, U_1, \Sigma_1, V_1$) are defined in Definition \ref{def:innerprodbalancing}.
\end{definition}

\begin{lemma}\label{lem:SBTinnerProd}
A symplectic balancing characterized by the test and trial basis matrices $(\Psi, \Phi)$  
 with $\Xi, \Xi' \in \SPD{n}$ and $J_\Omega \in \skewsymmetric{2n}$ has the following properties:
\begin{enumerate}[label=(\alph*)]
\setlength{\itemsep}{5pt}
\item The test and trial subsystem basis matrices $(\bar\Psi,\bar\Phi)$ balance $\Xi$ and $\Xi'$, i.e., $\bar\Phi\in O(\Xi,\Sigma_1)$ and $\bar\Psi\in O( \Xi',\Sigma_1)$. 
\item The test and trial (full-system) basis matrices $(\Psi,\Phi)$ balance $M=G^{-\tau} \diag(\Xi,\Xi')G^{-1}$ and $M'=G\diag(\Xi',\Xi)G^\tau$, i.e., $\Phi\in O
(M,\diag(\Sigma_1,\Sigma_1))$ and $\Psi\in O (M',\diag(\Sigma_1,\Sigma_1))$.
\item The basis matrices $(\Psi,\Phi)$  correspond to  a symplectic  projection with   $J_\Omega$ and $J_{2k}.$
\item The  basis matrices $(\Phi, \Psi)$  correspond to a symplectic projection with $J_{\Omega'}$ and $J_{2k}$, where $J_{\Omega'}=-J_\Omega^{-1}$.
\end{enumerate}
\end{lemma}
\begin{proof}
The conclusion (a) directly follows from  Lemma \ref{lem:GBTinnerProd} (c). Thus, we have $(\bar \Phi)^\tau \Xi \bar \Psi = (\bar \Psi )^\tau \Xi' \bar \Psi  = \Sigma_1$ and $(\bar \Psi)^\tau \bar \Phi = (\bar \Phi)^\tau \bar \Psi=I_k$.

To prove (b), we verify that $\Phi\in O(M,\diag(\Sigma_1,\Sigma_1))$ and $\Psi\in O(M',\diag(\Sigma_1,\Sigma_1))$. Using $(\bar \Phi)^\tau \Xi \bar \Psi = (\bar \Psi )^\tau \Xi' \bar \Psi  = \Sigma_1$, we have
\begin{gather*}
\Phi^\tau M \Phi=
\left (G 
\begin{bmatrix}
\bar \Phi & 0 \\
0 & \bar \Psi 
\end{bmatrix}
\right )^\tau
\left (
 G^{-\tau}
\begin{bmatrix}
\Xi & 0 \\
0 & \Xi'
\end{bmatrix}
G^{-1}
\right )
\left ( G
\begin{bmatrix}
\bar \Phi & 0 \\
0 & \bar \Psi 
\end{bmatrix}
\right )
 = 
\begin{bmatrix}
(\bar\Phi)^\tau\Xi\bar\Phi & 0  \\
  0 & (\bar\Psi)^\tau\Xi'\bar\Psi
 \end{bmatrix} = 
\begin{bmatrix}
\Sigma_1 & 0  \\
  0 & \Sigma_1
 \end{bmatrix}.
\end{gather*}
Similarly, we can obtain $\Psi^\tau M' \Psi =\diag(\Sigma_1, \Sigma_1)$.

To prove (c),
we verify that $\Phi \in\Sp(J_\Omega,J_{2k})$ and $\Psi = J_\Omega\Phi J_{2k}^{-1}$.
Using $G^\tau J_\Omega G =J_{2n}$ and $(\bar \Psi)^\tau \bar \Phi = (\bar \Phi)^\tau \bar \Psi=I_k$, we obtain
\begin{gather*}
\Phi^\tau J_{\Omega}\Phi = 
\begin{bmatrix}
(\bar\Phi)^\tau & 0  \\
  0& (\bar\Psi)^\tau
 \end{bmatrix}  \left( G^\tau J_{\Omega} G \right )
\begin{bmatrix}
\bar\Phi & 0  \\
  0& \bar\Psi
 \end{bmatrix}
= 
\begin{bmatrix}
(\bar\Phi)^\tau & 0  \\
  0& (\bar\Psi)^\tau
 \end{bmatrix}  J_{2n}
\begin{bmatrix}
\bar\Phi & 0  \\
  0& \bar\Psi
 \end{bmatrix}
=
\begin{bmatrix}
0 & (\bar\Phi)^\tau\bar\Psi  \\
 -(\bar\Psi)^\tau\bar\Phi & 0
 \end{bmatrix} = J_{2k},
\\
J_\Omega \Phi J_{2k}^{-1} =  \left (G^{-\tau}J_{2n}G^{-1} \right) 
\left (G
\begin{bmatrix}
\bar\Phi & 0  \\
  0& \bar\Psi
 \end{bmatrix}\right )
 \left( - J_{2k} \right )
  =
G^{-\tau}
\begin{bmatrix}
\bar\Psi & 0  \\
  0 & \bar\Phi
 \end{bmatrix} = \Psi.
\end{gather*}
 
This proves (d);  we can verify $\Psi \in\Sp(-J_\Omega^{-1},J_{2k})$ and $\Phi
=\left( -J_\Omega^{-1} \right ) \Psi J_{2k}^{-1}$ in a similar manner. 
\hfill
\end{proof}
\medskip


\begin{corollary}\label{cor:symplecticBalanceStable}
Suppose the original system $(A, B, C)$ is pure marginally stable, i.e.,  $A=JL\in
\generalizedham{2n}$ with $J\in \skewsymmetric{2n}$ and $L\in \SPD{2n}$.  If  $J_\Omega=-J^{-1}$,  then 
\begin{enumerate}[label=(\alph*)]
\setlength{\itemsep}{5pt}
\item The  reduced system  $(\tilde A, \tilde B, \tilde C)$ constructed by  symplectic balancing characterized by  $(\Psi, \Phi)$   remains pure marginally stable.
\item The  reduced system  $(-\tilde A^\tau, \tilde C^\tau, \tilde B^\tau)$ constructed by  symplectic balancing characterized by  $(\Phi, \Psi)$ remains pure marginally stable. 
\end{enumerate}
\end{corollary}
\begin{proof}
The conclusion (a)  directly follows from   Lemma \ref{lem:SBTinnerProd} (c) and Theorem \ref{thm:mainstab}.

If  $A\in \generalizedham{2n}$,  Corollary \ref{cor:decomp_A_Adual} (c) implies that $-A^\tau=J' L'$ with $J'=-J^{-1} \in \skewsymmetric{2n}$ and $L'\in \SPD{2n}$.  With the dual relationships $J_{\Omega'}=-J_\Omega^{-1}$ and $J'=-J^{-1}$,  the condition $J_\Omega=-J^{-1}$ yields $J_{\Omega'}=-(J')^{-1}$.
 Then, by Lemma \ref{lem:SBTinnerProd} (d) and Theorem \ref{thm:mainstab}, the conclusion (b) holds.
 \hfill
\end{proof}
\medskip

\begin{corollary}
Let $L, L' \in \SPD{2n}$ defined in \eqref{eq:generalizedHamPrimal} and \eqref{eq:generalizedHamDual},
performing symplectic balancing with $\Xi = \Xi' = \beta =
\diag(\beta_1,\ldots,\beta_n)$ preserves pure marginal stability and
balances the Hamiltonians of the primal and negative dual systems, i.e.,
$\Phi\in O(L,\Sigma_1)$, $\Psi\in O(L',\Sigma_1)$.
\end{corollary}
\begin{proof}
By \eqref{eq:generalizedHamPrimal}, $L=G^{-\tau} L_0 G^{-1}=G^{-\tau}\diag(\Xi,\Xi')G^{-1}$. By \eqref{eq:generalizedHamDual}, $L'=G L_0 G^{\tau}=G\diag(\Xi',\Xi)G^{\tau}$.
Then, the result follows trivially from Corollary \eqref{cor:symplecticBalanceStable} and Lemma \ref{lem:SBTinnerProd} (b).
\hfill
\end{proof}

Thus, in analogue to inner-product balancing for asymptotically stable systems,
performing symplectic balancing with $L$ and $L'$ (i.e., $\Xi = \Xi' = \beta =
\diag(\beta_1,\ldots,\beta_n)$ with $\beta_i>0$) not only preserves stability in the
appropriate (i.e., pure marginal) sense, it also balances the quadratic energy
functionals that characterize the primal and the dual systems. 

In analogue to Section  \ref{sec:existenceInner}, the next section presents
three algorithms (including symplectic balancing) for constructing basis
matrices that ensure symplectic projection given a subset of the required
ingredients.

\subsection{ Construction of  basis matrices given a subset of ingredients}\label{sec:PSD} 
Theorem \ref{thm:mainstab}  demonstrated that a ROM will preserve marginally
stability if it is constructed via symplectic  projection when the original
system has a symplectic structure.  Unfortunately,  most other model reduction
methods, such as POD--Galerkin and balanced POD, does not preserves the
symplectic structure, consequently the reduced model can be unstable.

We propose three methods (including symplectic balancing)  for constructing a
stability-preserving symplectic projection satisfying the conditions of
Definition \ref{def:sym_proj_dyn}. 
Table \ref{tab:comp_inner} summarizes these methods; this 
corresponds to Steps 6--8 in Algorithm
\ref{alg:overall}.

\begin{table}
\begin{center}
\caption{Algorithms for constructing a symplectic projection that ensure the
conditions of Definition \ref{def:sym_proj_dyn}.}
\label{tab:two_case_symplectic}
{
\begin{tabular}{|c|c|c|c|}
\hline
&Method 1 (symplectic balancing)&Method 2 &Method 3\\
\hline
\hline
 {Input} &
  \begin{tabular} {@{}l@{}}   
    $\Xi, \Xi'\in \SPD{n}$,\\
		$J_\Omega \in \skewsymmetric{2n}$ satisfying
		\eqref{eq:sym_stable_condition},\\
		$G$ satisfying \eqref{eq:generalizedHamPrimal}
     \end{tabular} 
 &
 \begin{tabular} {@{}l@{}}   
    $\Phi \in \Sp{(J_\Omega, J_\Pi)}$,\\
     $J_\Pi \in \skewsymmetric{2k}$, \\
     $J_\Omega \in \skewsymmetric{2n}$ \\satisfying \eqref{eq:sym_stable_condition}
     \end{tabular} 
    & 
    \begin{tabular} {@{}l@{}}   
    $\Phi_0\in \Sp(J_{2n},  J_{2k})$, \\
    $J_\Pi \in \skewsymmetric{2k}$, \\
    $J_\Omega \in \skewsymmetric{2n}$ \\satisfying \eqref{eq:sym_stable_condition}
     \end{tabular}   
        \\
\hline
 {Output} & 
 \begin{tabular} {@{}l@{}}
 $J_{\Pi}\in\skewsymmetric{2k}$,\\$\Phi\in \Sp{(J_{\Omega}, J_{2k})}$, $\Psi\in \Sp{(J_{\Omega'}, J_{2k})}$
\end{tabular}
 &
  $\Psi\in \RRstar{2n \times 2k}$  & $\Phi\in \Sp{(J_\Omega, J_\Pi)}$, $\Psi\in \RRstar{2n \times 2k}$   \\
\hline
{Algorithm} & 
 \begin{tabular} {@{}l@{}}   
    1.  Compute symmetric factorization\\
     \quad $\Xi=R R^{\tau}$, $\Xi'=S S^{\tau}$ \\
    2. Compute SVD
		$R^\tau S =  U \Sigma V^\tau$ \\
    3. $\bar \Phi = S V_1 \Sigma_1 ^{-1/2}$,
     $\bar \Psi = R U_1\Sigma_1 ^{-1/2}$\\
    4. $\Phi=G\diag(\bar \Phi, \bar \Psi)$,
    $\Psi=G^{-\tau}\diag(\bar \Psi, \bar \Phi)$ \\
    5. $J_\Pi = J_{2k}$
     \end{tabular}
&
    \begin{tabular} {@{}l@{}}   
    1. $\Psi =  J_\Omega \Phi J_\Pi^{-1}$
     \end{tabular} 
     &
     \begin{tabular} {@{}l@{}}   
    1.  Compute $G\in \Sp(J_\Omega, J_{2n})$\\
		\quad via Algorithm \ref{alg:constr_cannon}\\
    2.  Compute $\tilde G  \in  \Sp(J_\Pi, J_{2k})$ \\
		\quad via Algorithm \ref{alg:constr_cannon}\\
    3. $\Phi =G \Phi_0 \tilde G^{-1}$ \\
    4. $\Psi =  J_\Omega \Phi J_\Pi^{-1}$
     \end{tabular}
     \\
\hline
\end{tabular}}
\end{center}
\end{table}

Method 2 constructs a stability-preserving symplectic projection
starting with any trial basis matrix satisfying $\Phi \in Sp(J_\Omega,
J_\Pi)$,
$J_\Pi\in\skewsymmetric{2k}$, and
 $J_\Omega \in \skewsymmetric{2n}$ satisfying \eqref{eq:sym_stable_condition}.
Method 3 constructs a stability-preserving symplectic
projection starting with a basis $\Phi_0\in O(J_{2n}, J_{2n})$,
$J_\Pi\in\skewsymmetric{2k}$, and
 $J_\Omega \in \skewsymmetric{2n}$ satisfying \eqref{eq:sym_stable_condition}.
Lemma \ref{lem:sigma} demonstrates that we can compute the trial basis matrix
in this context as $\Phi= G \Phi_0^{\tau} \tilde G^{-1}\in Sp(J_\Omega, J_\Pi)$, which constitutes
step 3 of the algorithm.

\begin{lemma}\label{lem:sigma}
Let $\lift: (\mathbb{R}^{2k}, J_\Pi)\to (\mathbb{R}^{2n}, J_\Omega)$ denote a  symplectic lift with matrix presentation  $\Phi \in \Sp(J_\Omega, J_\Pi)$.
 Let $g: (\mathbb{R}^{2n}, J_{2n}) \to (\mathbb{R}^{2n}, J_\Omega)$ and $\tilde g: (\mathbb{R}^{2k}, J_{2k}) \to (\mathbb{R}^{2k}, J_{\Pi})$ represent symplectic transformations, represented by $G \in  \Sp(J_\Omega, J_{2n})$ and $\tilde G \in  \Sp(J_\Pi, J_{2k})$ respectively. Then, there exists a unique symplectic lift  $\lift_0:  (\mathbb{R}^{2k}, J_{2k})\to (\mathbb{R}^{2n}, J_{2n})$, represented by $\Phi_0\in \Sp(J_{2n}, J_{2k})$, such that the following diagram commutes:
\begin{center}
\begin{tikzpicture}[>=triangle 45]
\matrix(m)[matrix of math nodes,column sep={90pt,between origins},row
sep={50pt,between origins},nodes={asymmetrical rectangle}] (s)
{
|[name=M]| (\mathbb{R}^{n}, J_{2n})   &|[name=M1]| (\mathbb{R}^{2n}, J_\Omega)\\
|[name=N]| (\mathbb{R}^{k}, J_{2k})   &|[name=N1]|  (\mathbb{R}^{2k}, J_\Pi)\\
};
\path[->] (M) edge node[auto] {\(g\)} (M1)
          (N) edge node[left] {\(\lift_0\)} (M)
          (N1) edge node[auto] {\(\tilde g^{-1}\)} (N)
          (N1) edge node[right] {\(\lift\)} (M1);
\end{tikzpicture}
\end{center}
Equivalently,  for all $z \in (\mathbb{R}^{2k}, J_\Pi)$,
\begin{equation}
\lift(z)=g(\lift_0(\tilde g^{-1}(z))),
\end{equation}
and $\Phi =G \Phi_0 \tilde G^{-1}$ in matrix representation.
\end{lemma}
\begin{proof}
Because  $G \in \Sp(J_\Omega, J_{2n})$, we have $G^{\tau} J_\Omega G=J_{2n}$.
It follows that $G^{-\tau} J_{2n} G^{-1}= J_\Omega$.  By the same argument,
$\tilde G\in \Sp(J_\Pi, J_{2k})$ implies that $\tilde G^{\tau} J_\Pi  \tilde
G=J_{2k}$ and $\tilde G^{-\tau} J_{2k} \tilde G^{-1}=J_\Pi$.  Because $\Phi
\in \Sp(J_\Omega, J_\Pi)$, we have $\Phi^{\tau}J_\Omega \Phi=J_\Pi$. Because $g$ is  invertible, we can define $\lift_0:  (\mathbb{R}^{2k}, J_{2k})\to (\mathbb{R}^{2n}, J_{2n})$ by $g^{-1} \circ \lift \circ \tilde g$ with matrix representation $\Phi_0=G^{-1} \Phi \tilde G$. 
It follows that
\begin{equation*}
\Phi_0^{\tau} J_{2n} \Phi_0= \tilde G^{\tau} \Phi^{\tau} (G^{-\tau} J_{2n} G^{-1}) \Phi \tilde G=\tilde G^{\tau} (\Phi^{\tau}J_\Omega \Phi )\tilde G=\tilde G^{\tau} J_\Pi \tilde G
=J_{2k}.
\end{equation*}
The last equation implies that $\Phi_0\in \Sp(J_{2n}, J_{2k})$.  Finally, if $\lift_0$
satisfies $\lift =g\circ \lift_0  \circ \tilde g^{-1}$, $\lift_0$ is uniquely
determined by $\lift_0=g^{-1} \circ \lift \circ \tilde g$.
\hfill
\end{proof}

Apart from the symplectic-balancing approach we propose, there is no standard
method to construct a trial basis matrix satisfying $\Phi \in \Sp(J_\Omega,
J_\Pi)$. However, Ref.~\cite{PengL:16a} proposed several empirical methods to
construct $\Phi_0\in \Sp(J_{2n}, J_{2k})$, including cotangent lift (reviewed
in Appendix \ref{sec:cotangent_lift}), the complex SVD, and nonlinear
optimization.  Alternatively, we can also use a greedy
algorithm~\cite{HesthavenJS:17a} to construct $\Phi_0\in \Sp(J_{2n},
J_{2k})$ from empirical data. 


\section{Numerical examples}\label{sec:example}
This section illustrates the performance of the proposed structure-preserving
method (SP) using two numerical examples.  We  compare the full-order model
with reduced-order models constructed by POD--Galerkin (POD) (Appendix
\ref{sec:POD}), shift-reduce-shift-back method (SRSB) (Appendix
\ref{sec:SRSB}), balanced POD (BPOD) (Appendix \ref{sec:BPOD}), as well as
the proposed structure-preserving (SP) method. For reference, Table \ref{tab:four_methods} reports the algorithms for the existing
model-reduction methods. For simplicity, we focus on
autonomous systems $\dot x =Ax$ and employ the analytical solution $x(t)=
\exp(At) x_0$ as the `truth' solution. When applying BPOD and SRSB---which
require a full $(A,B,C)$ description---we set $B=C^\tau=x_0$.
For POD (see Appendix \ref{sec:POD}), we employ $N$ snapshots $\{x(i\Delta
t)\}_{i=0}^{N-1}$ with $\Delta t$ the specified snapshot interval. For balanced POD, we compute the primal
and dual snapshots according to \eqref{eq:defU} and \eqref{eq:defV},
respectively, with the same snapshot interval $\Delta t$.
For SRSB, we must define only the shift margin $\mu$. For each example, we compare two different SP methods: a POD-like method and a balancing method.

For time discretization, we define a uniform grid $\{t_i\}_{i=0}^{i=T-1}$ with 
$t_0=0$ and $t_{T-1}=t_f$, which employs a
uniform time step
$\delta t$ such that $t_i =  t_{i-1} + \delta t$. 
We apply the midpoint rule $x(t_{i+1}) = x(t_i)+ \frac{\delta t}{2}  A( x(t_i)
+x(t_{i+1} ))$ for performing time integration of both the full-order and
reduced-order models. When $A$ is a Hamiltonian matrix, this scheme
corresponds to a symplectic integrator; this ensures that
the time-discrete system will inherit any Hamiltonian structure that exists
in the time-continuous system.


To assess the accuracy of each method, we define the relative state-space error as 
\begin{equation}\label{relerror}
\eta =\frac{ \left (\sum_{i=0}^{T-1}  \|x(t_i)-\hat x(t_i)\|_2 ^2 \right )^{1/2}}{ \left (\sum_{i=0}^{T-1} \|x(t_i)\|_2^2 \right )^{1/2}
},
\end{equation} 
where $x(t_i)$ and $\hat x(t_i)$ denote the benchmark and approximate
solutions computed at time instance $t_i$.
 We also consider the relative system-energy error  as
\begin{equation}\label{relerrorenergy}
\eta_E =\frac{\left ( \sum_{i=0}^{T-1}  \left(E(t_i)-\hat  E(t_i))\right )^2 \right)^{1/2} }{\left( \sum_{i=0}^{T-1}
E(t_i)^2 \right) ^{1/2} },
\end{equation} 
where $E(t_i)$ and $\hat E(t_i)$ denote the benchmark and approximate system energies at time instance $t_i$.

\subsection{A 1D example}\label{sec:1Dexample}
To provide a simple illustration of the merits of the proposed technique, we first
consider a simple linear system with $n=8$, where
 \begin{equation}A=
\begin{bmatrix}
-8   & -29  & -72 & -139 & -192 & -171 & -128 &  -60 \\
     1  &   0  &    0 &    0 &    0 &    0 &    0 &    0 \\
     0  &   1  &    0 &    0 &    0 &    0 &    0 &    0 \\
     0  &   0  &    1 &    0 &    0 &    0 &    0 &    0 \\ 
     0  &   0   &    0 &    1 &    0 &    0 &    0 &    0\\ 
     0  &   0    &   0 &    0 &    1 &    0 &    0 &    0 \\
     0  &   0   &    0 &    0  &   0 &    1 &    0 &    0 \\
     0 &     0  &  0   &  0   &  0    & 0     & 1    & 0
 \end{bmatrix}
  \end{equation}
The eigenvalue decomposition of $A$ gives $A = P_1 \Lambda P_1^{-1}$, where $\Lambda=\diag(-3, -2+ i, -2-i, -1, 2i, -2i, i, -i)$ and $P_1 \in \mathbb{C}_*^{8\times 8}$.
Thus, the original system is marginally stable. Then we construct $T=\begin{bmatrix} T_s \ T_m \end{bmatrix}\in \RRstar{n\times n}$ such that $A T_i= T_i A_s$ (for $i\in \{s, m\}$), where 
the spectrum of $A_s$ and $A_m$ are given by $\{-3, -2\pm i, -1\}$ and $\{\pm 2i, \pm i\}$, respectively. In particular, we obtain 
\begin{equation} 
 A_s =
\begin{bmatrix}
   -3   &  0  &     0 &    0  \\
     0  &   -2 &    1 &    0 \\
     0  &   -1  &   -2 &    0  \\
     0  &   0  &    0 &    -1 
 \end{bmatrix},\quad
  A_m =
 \begin{bmatrix}
       0 &    0 &    2 &    0\\ 
       0 &    0 &    0 &    1 \\
      -2 &    0 &    0 &    0 \\
      0    &   -1     &  0    & 0
      \end{bmatrix} .
\end{equation} 
We  choose
$ M=\Theta = \diag \left (   \frac{1}{6}, \frac{1}{4}, \frac{1}{4}, \frac{1}{2}   \right )$
(which satisfies the Lyapunov equation \eqref{eq:semiequiLya} with $Q = I_4$)  
 and $J_\Omega=J_4$ (which satisfies $A_m=J_4 L_0$ with $L_0=\diag(2,1, 2,1)$). 
 
We test two SP methods; both of them reduce $A_s$ and $A_m$ from dimension $n/2$ to dimension $k/2$. The first SP method, SP1, 
is a POD-like method: SP1 applies Method 2 in Table \ref{tab:two_case} for the
asymptotically subsystem, where $\Phi$ is computed via POD with snapshots  $\{x_s(i\Delta t)\}_{i=0}^{N-1}$;
SP1  applies Method 2 in Table
\ref{tab:two_case_symplectic} for the pure marginally stable subsystem, where $J_\Pi$ is a Poisson matrix and  $\Phi$ is constructed via cotangent
lift (see Algorithm \ref{alg:lift} of Appendix \ref{sec:cotangent_lift}) with
snapshots $\{x_m(i\Delta
t)\}_{i=0}^{N-1}$. The second SP method, SP2,  is a balancing method:  SP2
applies the balanced truncation for the asymptotically stable subsystem and
symplectic balancing with $\Xi = \Xi' = \mathrm{diag}(\beta,\beta)$ for the
pure marginally stable subsystem; this approach balances the primal and
negative dual Hamiltonians; this approach balances the primal and negative
dual Hamiltonians.

We set the initial condition to the first canonical unit vector, i.e.,
$x_0=e_1$. For the purpose of constructing basis matrices, we collect $N = 11$ snapshots from the time domain $[0, 5]$ with snapshot interval $\Delta t=0.5$.
For SRSB, we set the shift margin to $\mu=0.01$.

\begin{figure}
\begin{center}
\subfigure[The evolution of the state-space error  $\|e(t)\|=\|x(t)-\hat x(t)\|$
]{
\includegraphics[width=0.45\textwidth]{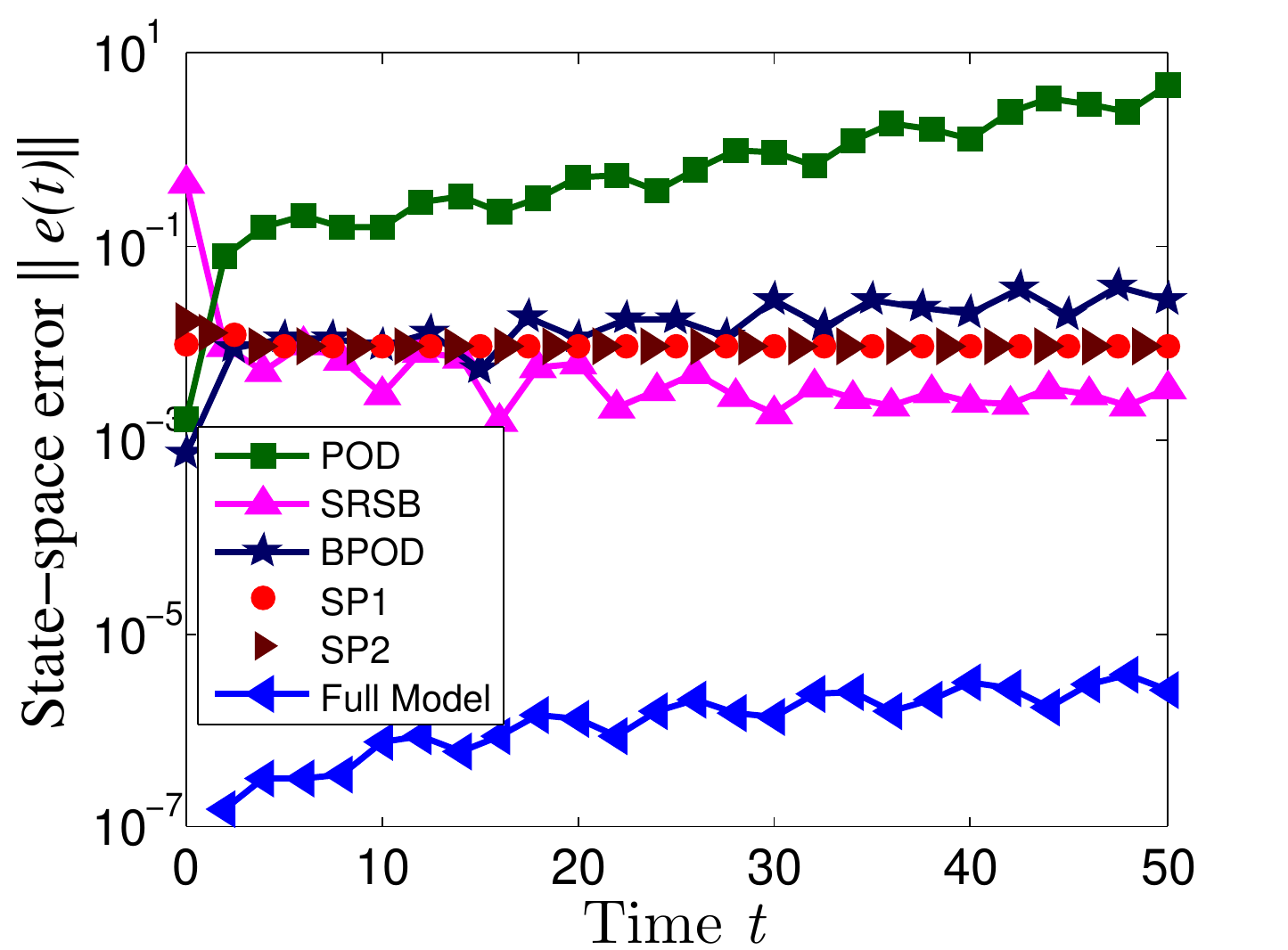}
	}
\subfigure[The evolution of the system energy $E(t)$
]{
\includegraphics[width=0.45\textwidth]{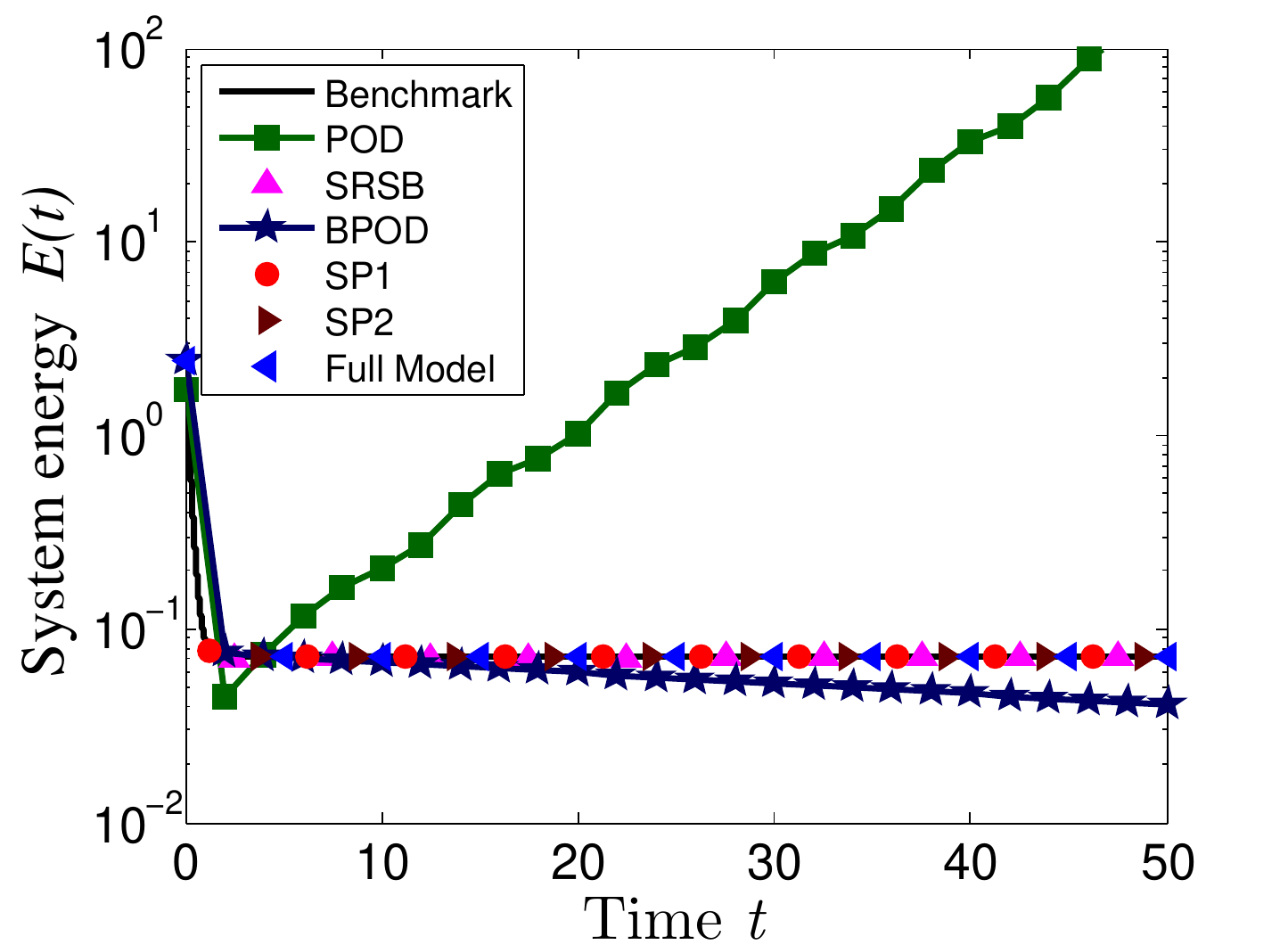}
	}
\caption{\textit{1D example}.
The evolution of the state-space error and energy for all tested
methods and reduced dimension $k=4$.} \label{fig:err_1D}
\vspace{-3mm}
\end{center}
\end{figure}

All experiments consider reduced-order models of dimension $k=4$.
Table \ref{tab:test_eg1} compares the performance of  different methods for
this example; we compute the infinite-time energy via eigenvalue analysis. We
choose a longer time interval characterized by $t_f=50$ (and time step $\delta
t =0.001$) to compute the errors $\eta$ and $\eta_E$; thus, $T=50001$
in \eqref{relerror} and \eqref{relerrorenergy}.  Figure \ref{fig:err_1D}  (a) plots
the evolution of the $\ell^2$-norm of the state-space error $e(t)\defeq x(t)-\hat{x}(t)$.
 Figure \ref{fig:err_1D}  (b) plots
the evolution of the system energy $E(t)$ , which is defined in \eqref{ref:energy} of Appendix
\ref{sec:sysenergy}.


\begin{table}
\begin{center}
\caption{\textit{1D example}. Comparison of different model-reduction methods
for reduced dimension $k=4$.}
\label{tab:test_eg1}
 \resizebox{\textwidth}{!}
 {%
\begin{tabular}{|c|c|c|c|c|c|c|}
\hline
& POD & SRSB & BPOD & SP1 & SP2  & \begin{tabular}{@{}c@{}}  Full-order \\ model   \end{tabular}   \\
\hline
\hline
  Eigenvalues $\lambda$     & \begin{tabular}{@{}c@{}}  -6.9457 \\ -0.4456 \\  $0.0828\pm1.9679i$  \end{tabular} & \begin{tabular}{@{}c@{}} $-0.0008\pm0.9999i$ \\ $0.0002\pm 2.0002i$  \end{tabular}  & \begin{tabular}{@{}c@{}} $-2.8850 \pm 0.5713i$ \\ $-0.0063 \pm 2.0060i$  \end{tabular}  & \begin{tabular}{@{}c@{}}  $-2.3590\pm 0.3684i$  \\  $\pm 1.9998i$  \end{tabular}   & \begin{tabular}{@{}c@{}}  $-2.8663\pm 1.8442i$  \\  $\pm 2.0000 i$  \end{tabular}   & 
  \begin{tabular}{@{}c@{}}  $-3$  \\ $-2\pm i$\\ $-1$\\ $\pm i$ \\ $\pm 2i$  \end{tabular}  \\
 \hline
 \begin{tabular}{@{}c@{}}  Marginal-stability \\ preservation \end{tabular}   &  No & No  & Yes & Yes & Yes & Yes\\
 \hline
 \begin{tabular}{@{}c@{}}  Relative  state-space \\ error $\eta$ \end{tabular}    &  16.2082 &  0.2245   &  0.1936 & 0.0870  & 0.0774 & $1.9269 \times 10^{-5}$  \\
\hline
\begin{tabular}{@{}c@{}}  Relative system-energy\\ error $\eta_E$ \end{tabular}   &  302.2681  & 0.7900  & 0.1401 & 0.0055 & 0.0852 & $2.8761 \times 10^{-7}$ \\
\hline
\begin{tabular}{@{}c@{}}  Infinite-time\\ energy\end{tabular}  & $+\infty$ & $+\infty$  & 0 & 0.07179  & 0.07181 & 0.07216 \\
\hline
\end{tabular}}
\end{center}
\end{table}

First, note that POD yields the largest state and system-energy errors, and
its energy grows rapidly, even within the considered time interval.  This can
be attributed to its eigenvalues ($\lambda = 0.0828 \pm 1.9679 i$), which
correspond to unstable modes. Because the POD reduced-order model is unstable
(its matrix $\tilde A$ has eigenvalues with positive real parts), its
infinite-time energy will be unbounded. While SRSB yields lower errors $\eta$
and $\eta_E$ than POD, it has larger errors over the first part of the time
interval.  Even though SRSB
does not preserve marginal stability, its instability margin is  only $2\times
10^{-4}$, which is relatively small and precludes instabilities from becoming
apparent over the finite time interval considered; however, its infinite-time
energy is unbounded.  BPOD has smaller average errors than both POD and SRSB;
however, the associated reduced-order model is asymptotically stable, which
implies that its infinite-time energy is zero; thus, the reduced-order model
does not have a pure marginal subsystem.  Not only do the proposed SP methods
produce the smallest average errors over all reduced-order models, they are
also the only methods that preserve marginal stability, including the pure
marginally stable subsystem. As a result, two proposed SP methods yield a
finite infinite-time energy; in fact, this energy incurs a sub-1\% error with
respect to the infinite-time energy of the full-order model.  Critically, note
that extreme pure imaginary eigenvalues are exact ($\pm 2i$) in the case of SP2; this
results from the fact that it balances the Hamiltonians directly. 

\subsection{2D mass--spring system}
We now consider a 2D mass--spring system. Each mass is located on a grid point
of an $(\bar n+2)\times (\bar n+2)$ grid with $\bar n = 49$.
The governing equations associated with mass $(i,j)$, $i,j =1,\ldots \bar n$ is given by
\begin{equation}\label{eq:2dmspring}
  \begin{aligned}
  m \ddot u_{i,j}&=k_x(u_{i+1,j}+u_{i-1,j}-2 u_{i,j}) -2b \dot u_{i,j},\\
  m \ddot v_{i,j}&=k_y(v_{i,j+1}+v_{i,j-1}-2 v_{i,j}),
    \end{aligned}
\end{equation}
where $u_{i,j}$ and $v_{i,j}$ are state variables representing the $x$- and
$y$-displacements of mass $(i,j)$, $m=1$ denotes the mass, $k_x$
and $k_y$ denote spring constants with $k_x=k_y=2500$, and $b=1$ denotes the damping coefficient in
the $x$-direction. We apply homogeneous Dirichlet boundary conditions
$u_{0,j}=u_{\bar n+1, j}=v_{i,0}=v_{i,\bar n+1}=0$.

Define canonical coordinates: $q_{i,j}=u_{i,j}$, $p_{i,j}=m \dot q_{i,j}$,
$r_{i,j}=v_{i,j}$, and $s_{i,j}=m \dot r_{i,j}$.  The system Hamiltonian is
given by $H(q_{i,j}, p_{i,j}, r_{i,j}, s_{i,j})=H_x(q_{i,j}, p_{i,j})+H_y(r_{i,j}, s_{i,j})$, where     
\begin{equation}
  \begin{aligned}
& H_x(q_{i,j}, p_{i,j})=\frac{1}{2m} \sum_{i,j=1}^{\bar n} p_{i,j}^2 + \frac{k_x}{2}  \sum_{i,j=0}^{\bar n}  (q_{i+1,j}-q_{i,j})^2 ,\\
&H_y(r_{i,j}, s_{i,j})=\frac{1}{2m} \sum_{i,j=1}^{\bar n} s_{i,j}^2 + \frac{k_y}{2} \sum_{i,j=0}^{\bar n}  (r_{i,j+1}- r_{i,j})^2.
    \end{aligned}
\end{equation}

Now, the  original system \eqref{eq:2dmspring} can represented by dissipative
Hamiltonian ordinary differential equations,
\begin{equation}\label{eq:Hamilmspring}
  \begin{aligned}
\dot q_{i,j}&=\frac{\partial H}{\partial p_{i,j}},  \quad \quad \quad \dot p_{i,j}=-\frac{\partial H}{\partial q_{i,j}}- \frac{2b}{m} p_{i,j},\\
\dot r_{i,j}&=\frac{\partial H}{\partial s_{i,j}}, \quad \quad \quad \dot s_{i,j}=-\frac{\partial H}{\partial r_{i,j}}.\\
  \end{aligned}
\end{equation}
Let 
$q=[q_{1,1}\ \cdots\ q_{1,\bar n}\ \cdots\ q_{\bar n,1}\ \cdots\ q_{\bar n, \bar n}]^{\tau}$
and $r=[r_{1,1}\ \cdots\ r_{1, \bar n}\ \cdots\ r_{\bar n,1}\ \cdots\ r_{\bar n,\bar n}]^{\tau}$ denote the generalized coordinates. Let   
$p=[p_{1,1}\ \cdots\ p_{1,\bar n}\ \cdots\ p_{\bar n,1}\ \cdots\ p_{\bar n, \bar n}]^{\tau}$
and 
$s=[s_{1,1}\ \cdots\ s_{1,\bar n}\ \cdots\ s_{\bar n,1}\ \cdots\ s_{\bar n, \bar n}]^{\tau}$
denote the generalized momenta. 
With $\psi_x=[q^{\tau}\ p^{\tau}]^{\tau}\in \mathbb{R}^{2\bar n^2}$ and
$\psi_y=[r^{\tau}\ s^{\tau}]^{\tau} \in \mathbb{R}^{2\bar n^2}$, the above equation can be written as matrix form, i.e.,
  \begin{equation}
  \frac{d}{dt}
     \begin{bmatrix}
    \psi_x \\ \psi_y
     \end{bmatrix}
  =
  \begin{bmatrix}
A_s & 0  \\
 0 & A_m  \\
 \end{bmatrix} 
      \begin{bmatrix}
    \psi_x \\ \psi_y
     \end{bmatrix},
 \end{equation}
where $\dot \psi_x=A_s \psi_x$ represents an asymptotically stable system and
$\dot \psi_y=A_m \psi_y$ represents a (pure marginally stable) Hamiltonian system.
Thus, the dimension of the full-order model is $n = 4\bar n^2$.

Because the original system is neither controllable
nor observable, SRSB cannot be directly used, as it requires solvability of
Lyapunov equations \eqref{eq:reWo} and \eqref{eq:reWc}. Instead, we compute the Gramians
 $M_o^\mu$ and  $M_c^\mu$ by
solving the modified Lyapunov equations  $(A-\mu I)^{\tau} M_o^\mu + M_o^\mu
(A-\mu I) =-(C^{\tau}C+\varepsilon I)$ and $(A-\mu I) M_c^\mu + M_c^\mu (A-\mu
I)^{\tau} =-(BB^{\tau}+\varepsilon I)$, and let $\varepsilon=10^{-4}$. We set the shift margin to $\mu =1$.

We again test two SP methods; both of them reduce
$A_s$ and $A_m$ from dimension $n/2$ to dimension $k/2$. The first SP method,
SP1, is identical to the SP1 method employed in the previous
example.  The second SP method, SP2,  applies a different balancing approach. Because the
asymptotically stable  subsystem is neither controllable nor observable,
balanced truncation cannot be directly used for this subsystem as well. Thus,
we also compute Gramians  by solving the modified Lyapunov equations 
 $A^{\tau} M_o + M_o
A=-(C^{\tau}C+\varepsilon I)$
and
$A M_c +
M_c A^{\tau} =-(BB^{\tau}+\varepsilon I)$. For the pure marginally stable subsystem, we
collect snapshot ensemble  $\{\psi_y(i\Delta t)\}_{i=0}^{N-1}$ and construct two
snapshot matrices $R=\begin{bmatrix}r_0\ \cdots\ r_{N-1} \end{bmatrix}$ and
$S=\begin{bmatrix}s_0\ \cdots\ s_{N-1} \end{bmatrix}$ in $\mathbb{R}^{\bar n^2 \times N}$, where $\psi_y(i \Delta
t)=\begin{bmatrix}r_i ^\tau\ s_i^\tau \end{bmatrix}^\tau$. Then, the symplectic balancing method (the
first method in Table \ref{tab:two_case_symplectic}) is employed with $\Xi =
RR^\tau$ and $\Xi' =SS^\tau$.  Since the pure marginally stable
subsystem in this example is a standard Hamiltonian,  we have $J_\Omega=J_{n/2}$ and $G=I_{n/2}$.

\begin{figure}
\begin{center}
\subfigure[$q(0)$]{
\includegraphics[width=0.3\linewidth]{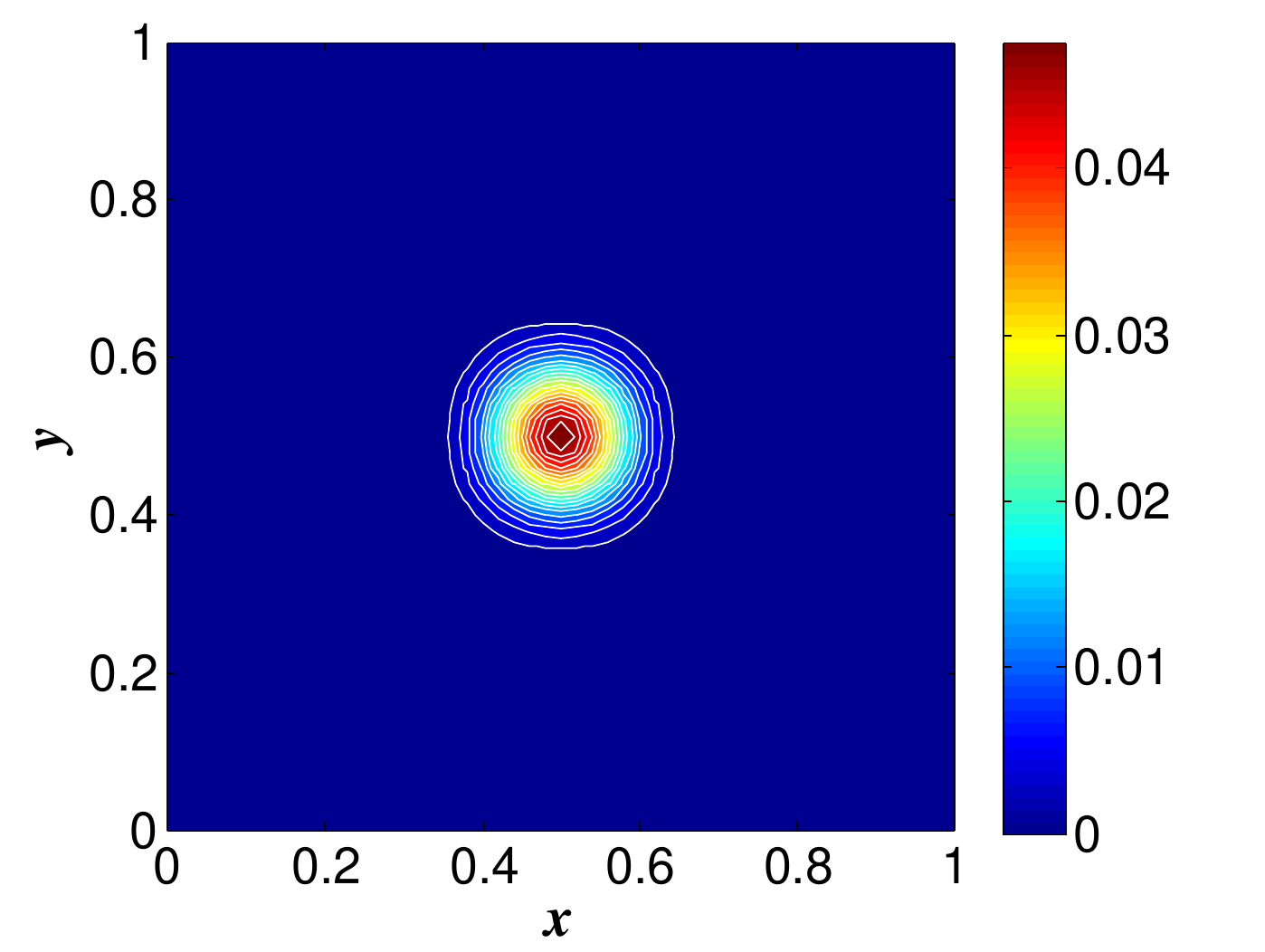}
}
\subfigure[$r(0)$]{
\includegraphics[width=0.3\linewidth]{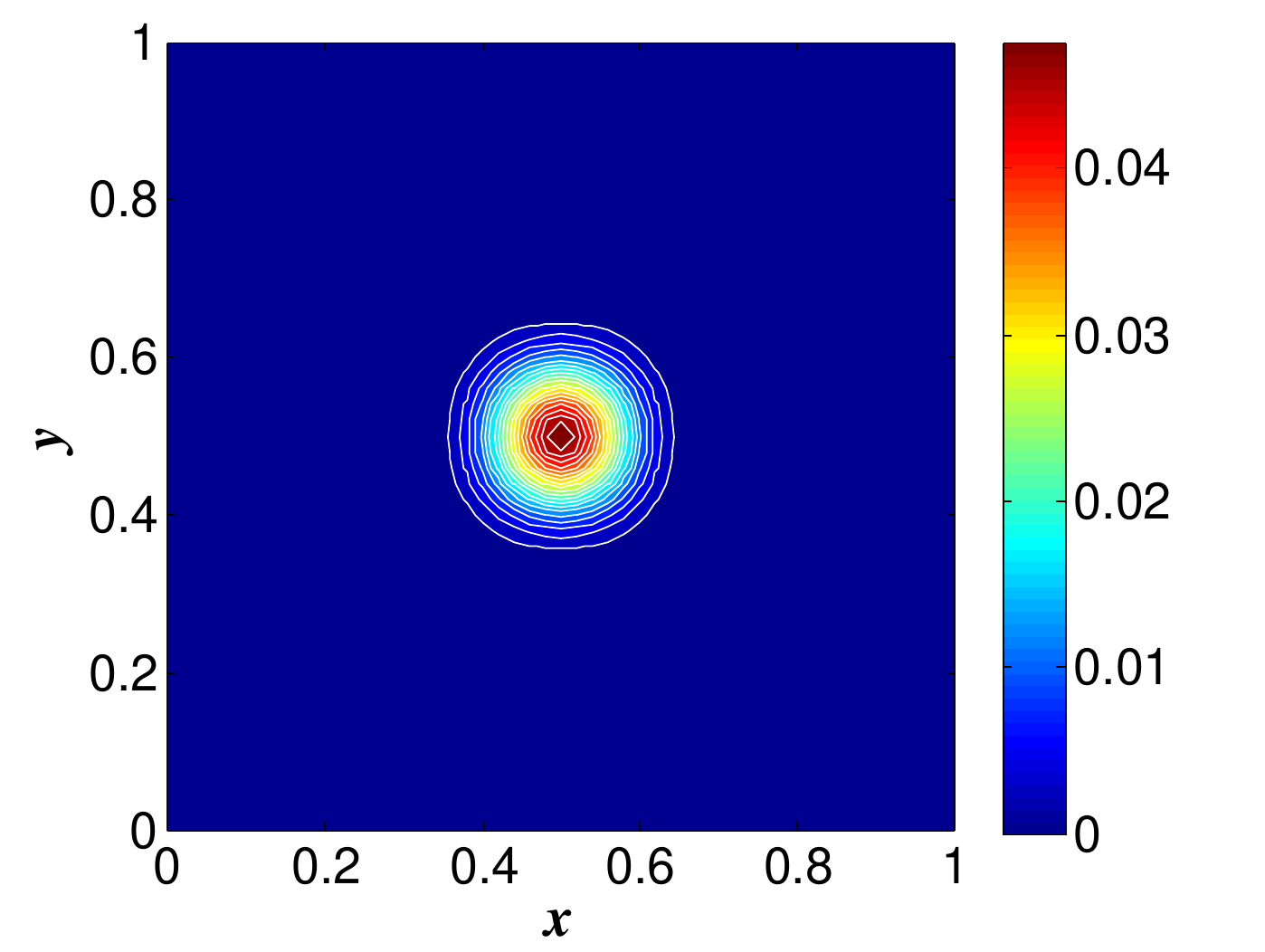}
}
\subfigure[Mass $(i,j)$ at $t=0$, $i,j=0, \ldots, 50$]{
\includegraphics[width=0.3\linewidth]{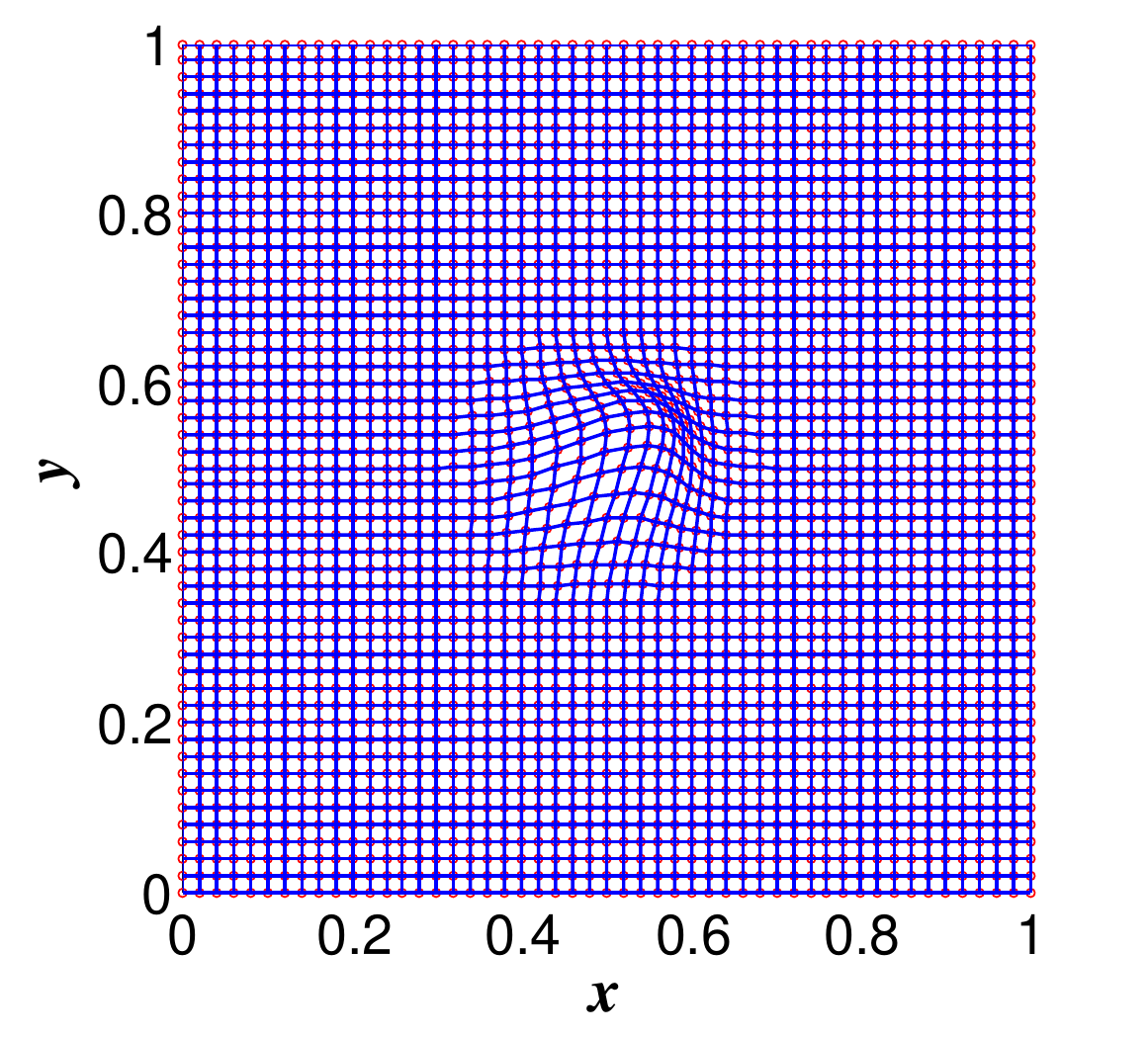}
}
\subfigure[$q(t_f)$]{
\includegraphics[width=0.3\linewidth]{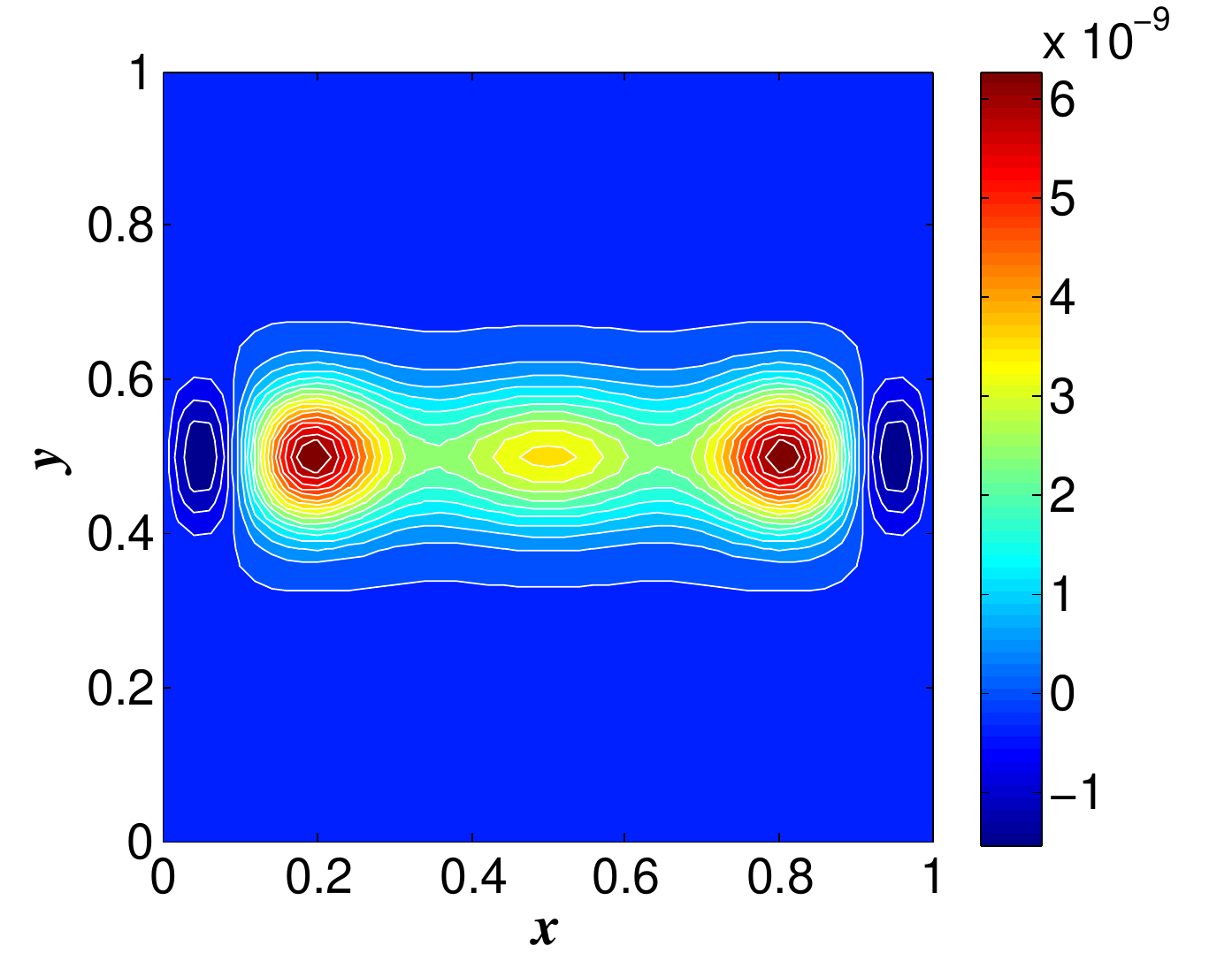}
}
\subfigure[$r(t_f)$]{
\includegraphics[width=0.3\linewidth]{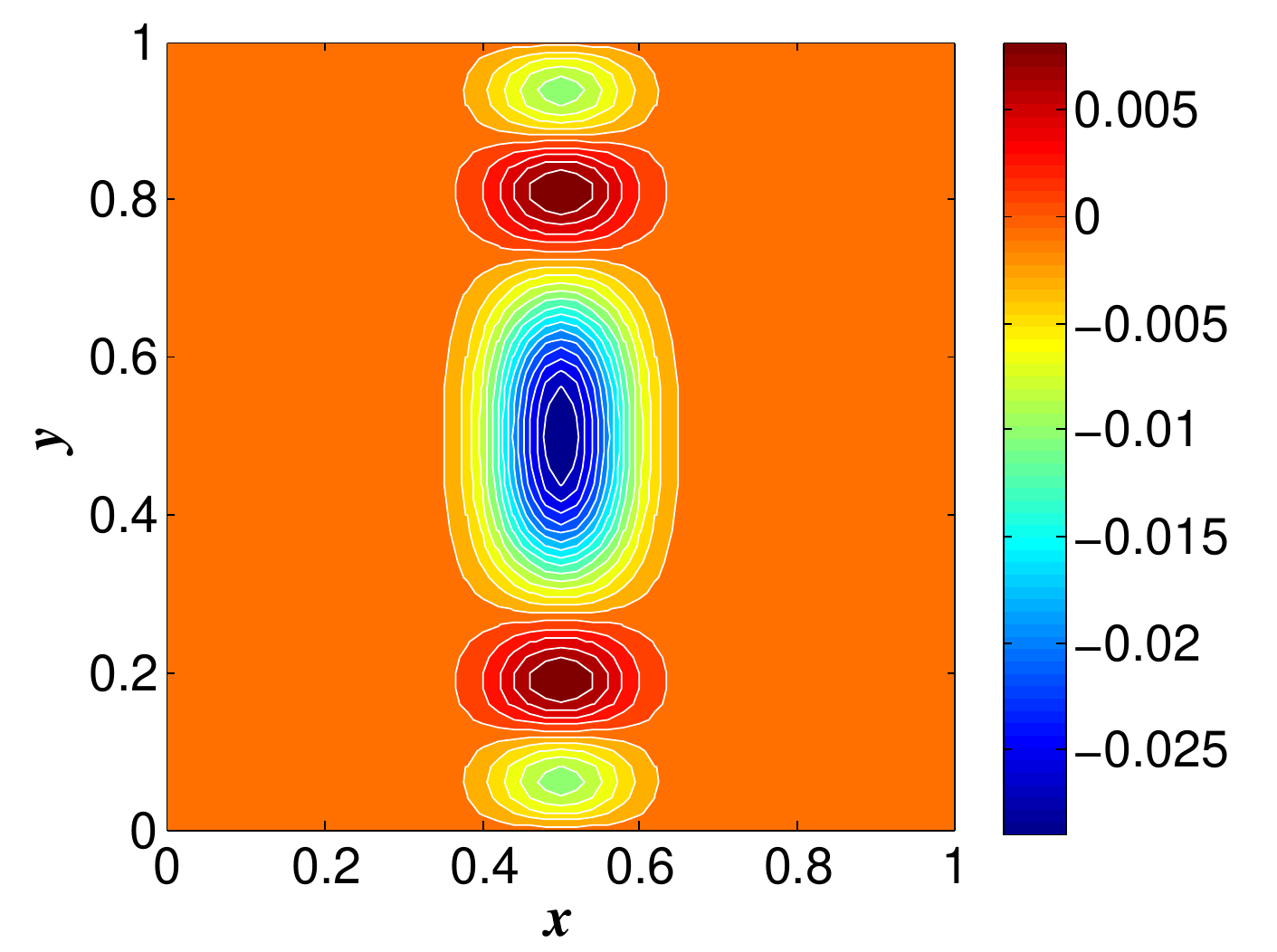}
}
\subfigure[Mass $(i,j)$ at $t=t_f$, $i,j=0, \ldots, 50$]{
\includegraphics[width=0.3\linewidth]{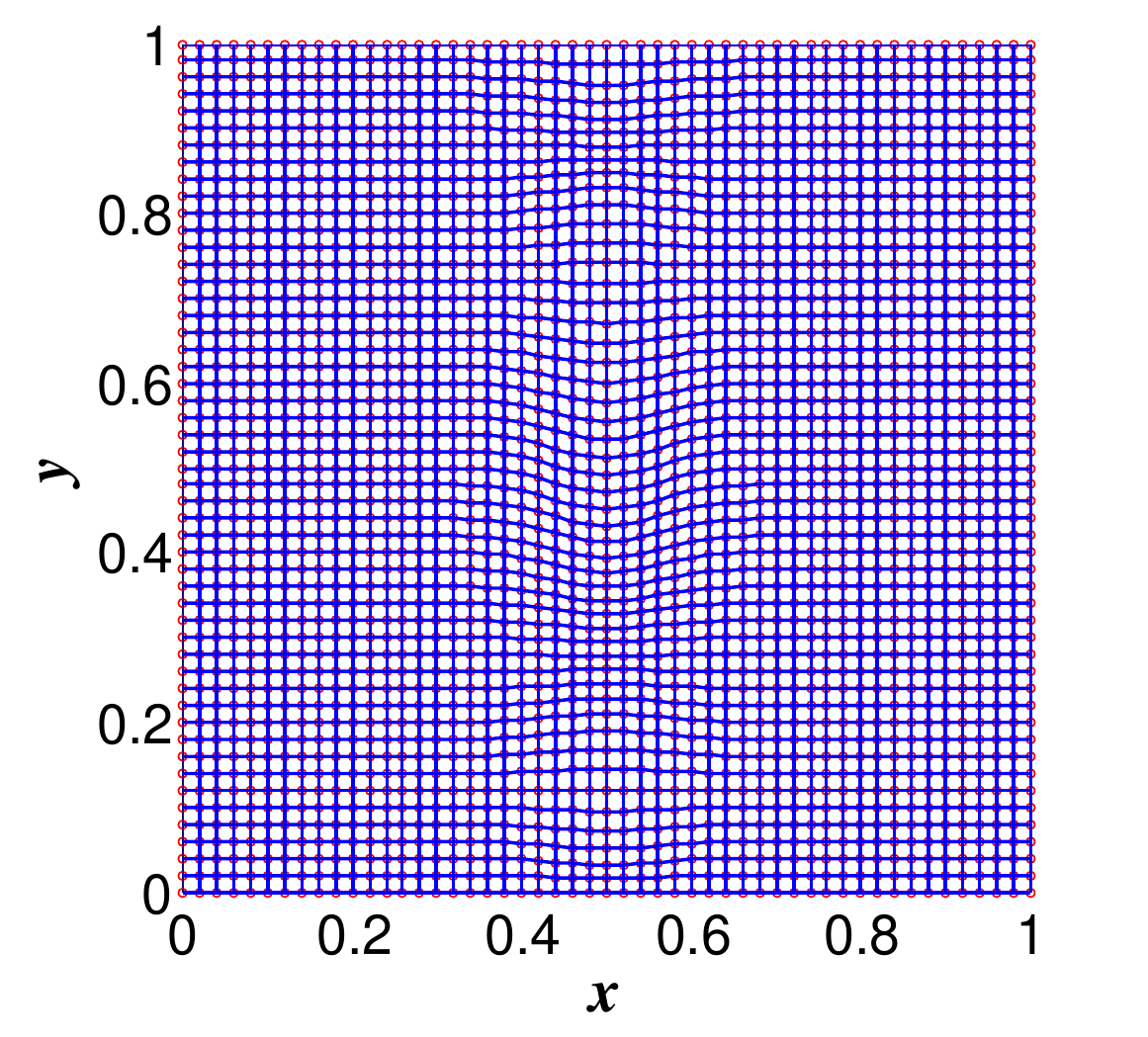}
}
 \caption{\textit{2D mass--spring example}. Initial condition and final state.  } \label{fig:snap_init_final}
\end{center}
\end{figure}

Let  $\alpha(x)=\frac{ |x-{l/2} |}{l/10}$ with $l=1$ the length of the spatial
interval in each direction and  $h(\alpha)$ be a cubic spline
\begin{equation*}
\begin{array}{l}
h(\alpha)=\left\{
\begin{array}{ccl}
 \vspace{3pt} 1-\frac{3}{2}\alpha^2+\frac{3}{4}\alpha^3&  \text{if }&  0\le \alpha \le 1, \\
 \frac{1}{4}(2-\alpha)^3&  \text{if }&  1<\alpha \le2 ,\\ 0&  \text{if }&  \alpha>2.
\end{array} \right. 
\end{array}
\end{equation*}
Let $x_i=il/\bar n$ and $y_j=jl/\bar n$.
For our numerical experiments,  the initial condition is provided by
\begin{equation}\label{initial_2D}
q_{i,j}(0)=   r_{i,j}(0) = h(\alpha(x_i)) h(\alpha(y_j)), \quad \quad \quad
p_{i,j}(0)=s_{i,j}(0)=0.
\end{equation}
We employ a time step of 
$\delta t=0.002$ and set the final time to $t_f = 15$ to compute the errors $\eta$ and $\eta_E$; Thus, $T=75001$. Figure
\ref{fig:snap_init_final} depicts the initial condition and final state computed by
the full-order model. For the purpose of constructing basis matrices, we collect $N= 101$ snapshots from the time domain $[0, 5]$ with snapshot interval $\Delta t=0.05$. 

\begin{figure}
\begin{center}
\subfigure[The evolution of the state-space error  $\|e(t)\|=\|x(t)-\hat x(t)\|$]{
\includegraphics[width=0.45\linewidth]{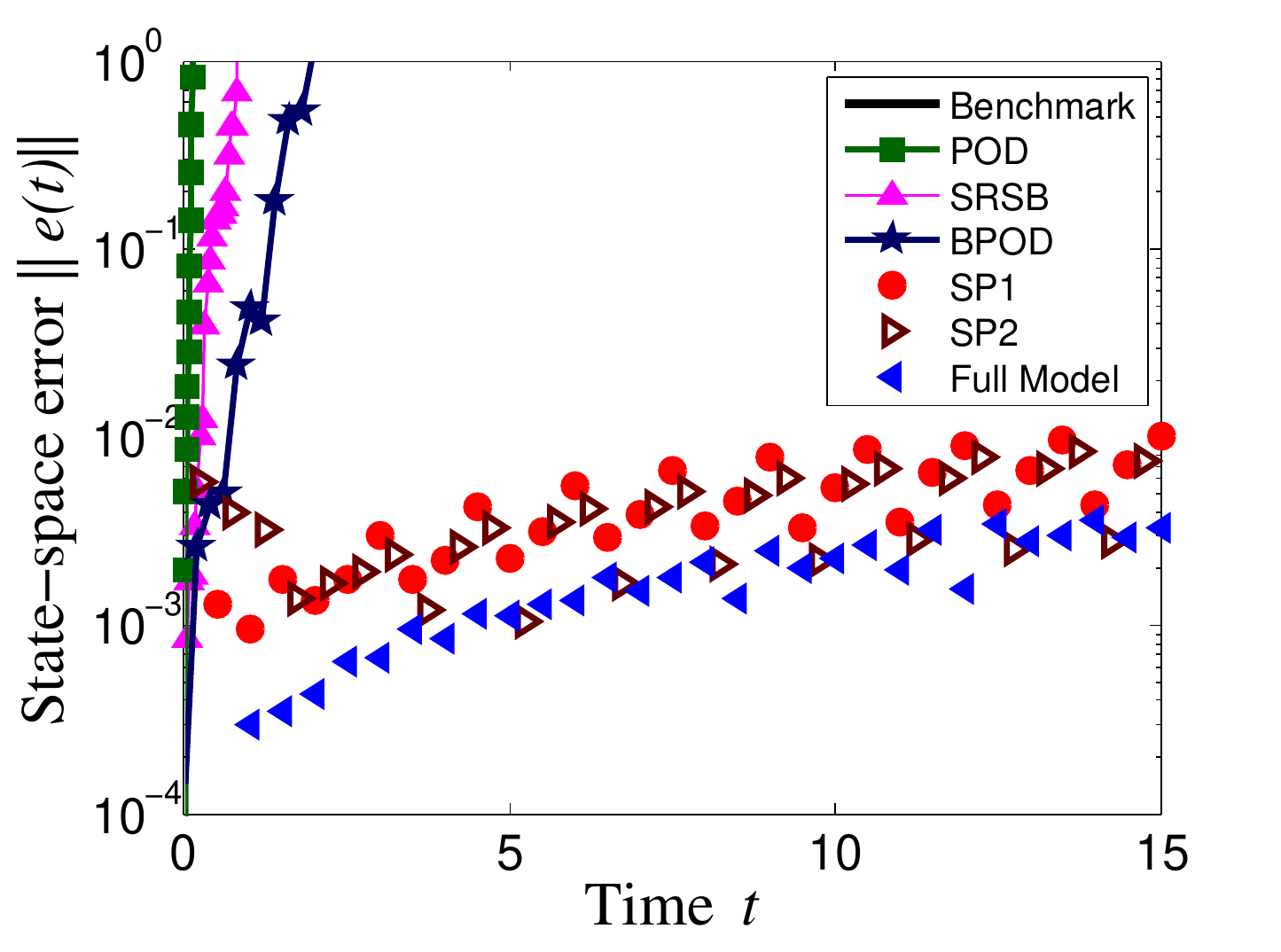}
}
\subfigure[The evolution of the system energy $E(t)$]{
\includegraphics[width=0.45\linewidth]{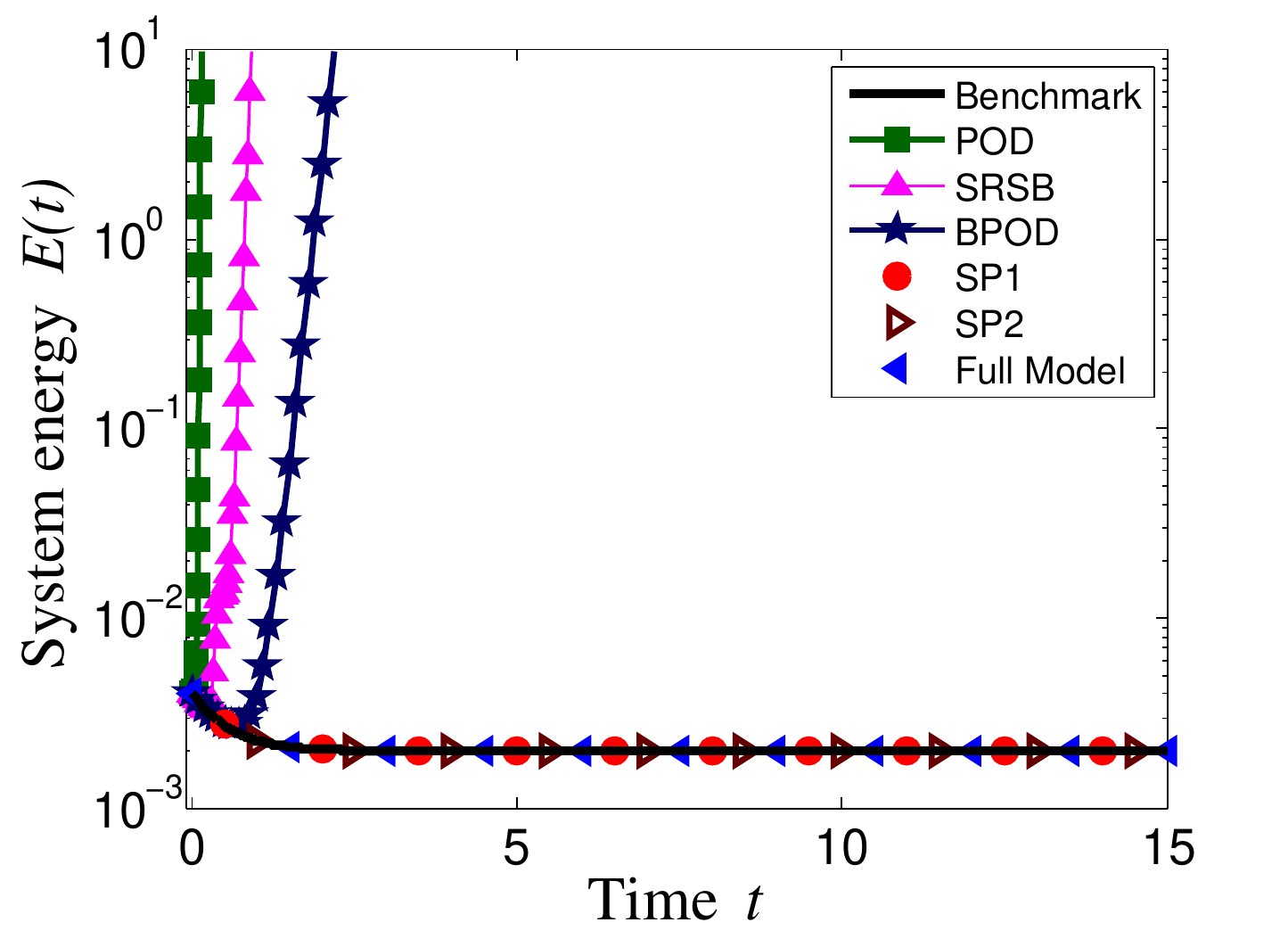}
}
 \caption{\textit{2D mass--spring example}. The evolution of the state-space error  $\|e(t)\|=\|x(t)-\hat x(t)\|$ and
 system energy $E(t)$  for all tested methods and reduced dimension $k=40$.} 
 \label{fig:err_2D}
\end{center}
\end{figure}

Table \ref{tab:test_eg2} compares the performance of  different
reduced-order models (all of dimension $k=40$), while Figure \ref{fig:err_2D}  plots the $\ell^2$-norm
of the state-space error $e(t)\defeq x(t)-\hat{x}(t)$ and the system energy
$E(t)$ for those reduced-order models as a function of time. 
Here, the system energy
is defined by the total Hamiltonian,  i.e.,
$E=H(q_{i,j}, p_{i,j}, r_{i,j}, s_{i,j})$, and its infinite-time value is computed
by eigenvalue analysis.


\begin{table}
\begin{center}
\caption{\textit{2D mass--spring example}. Comparison of different
model-reduction methods for reduced dimension $k=40$.}
\label{tab:test_eg2}
{%
\begin{tabular}{|c|c|c|c|c|c|c|}
\hline
& POD & SRSB &  BPOD  & SP1 &  SP2 & \begin{tabular}{@{}c@{}}  Full-order \\ model   \end{tabular} \\
\hline
\hline
\begin{tabular}{@{}c@{}}  Number of \\ unstable modes \end{tabular} & 8 & 16  & 18 & 0 & 0 &0 \\
 \hline
\begin{tabular}{@{}c@{}}  Instability margin \\ $\max ({\rm{Re}}(\lambda))$  \end{tabular}    & 50.480  & 10.586  & 3.695 & 0  & 0 &0 \\
 \hline
\begin{tabular}{@{}c@{}}  Marginal-stability\\ preservation \end{tabular}     &No & No & No & Yes & Yes   & Yes\\
\hline
\begin{tabular}{@{}c@{}}  Relative  state-space \\ error $\eta$   \end{tabular}      &  $+\infty$ & $+\infty$ & $+\infty$ &  0.11156 &  0.10214& 0.04358 \\
\hline
\begin{tabular}{@{}c@{}}  Relative  system-energy \\error $\eta_E$\end{tabular}   &$+\infty$ & $+\infty$ & $+\infty$ & $8.6868 \times 10^{-5}$ & $4.8843 \times 10^{-3}$ & $3.413 \times 10^{-5}$\\
\hline
\begin{tabular}{@{}c@{}}  Infinite-time\\ energy \end{tabular} & $+\infty$ & $+\infty$  & $+\infty$ & $1.9958\times 10^{-3}$ & $1.9959\times 10^{-3}$   & $1.9959\times 10^{-3}$ \\
\hline
\end{tabular}}
\end{center}
\end{table}

First, note that
among all the tested methods,  only the full-order model and the proposed SP
reduced-order models
preserve marginal stability and have finite errors $\eta$ and $\eta_E$.
Further, the SP methods ensure that the reduced-order model has a pure
marginally stable subsystem, and thus a finite infinite-time energy that is
nearly identical to that of the full-order model. Because POD,
SRSB, and BPOD have unstable modes, they yield unbounded infinite-time
energy. Further, due to their relatively
large instability margins, their errors and energy grow rapidly within the
considered time interval, leading to significant errors.

\begin{figure}
\begin{center}
\subfigure[Relative state-space error $\eta$ versus subspace dimension $k$]{
\includegraphics[width=0.45\linewidth]{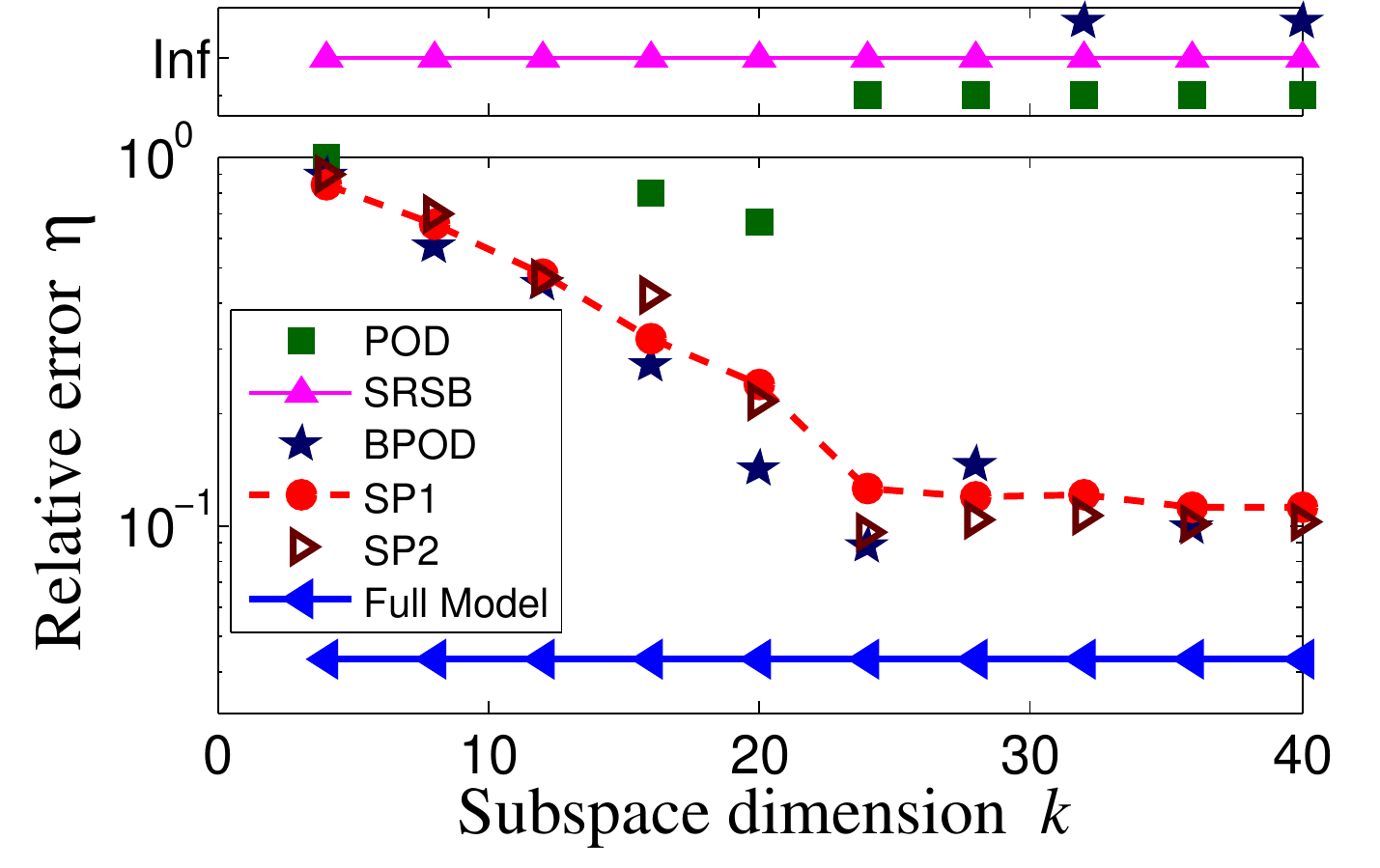}
}
\subfigure[Relative system-energy error $\eta_E$ versus subspace dimension $k$]{
\includegraphics[width=0.45\linewidth]{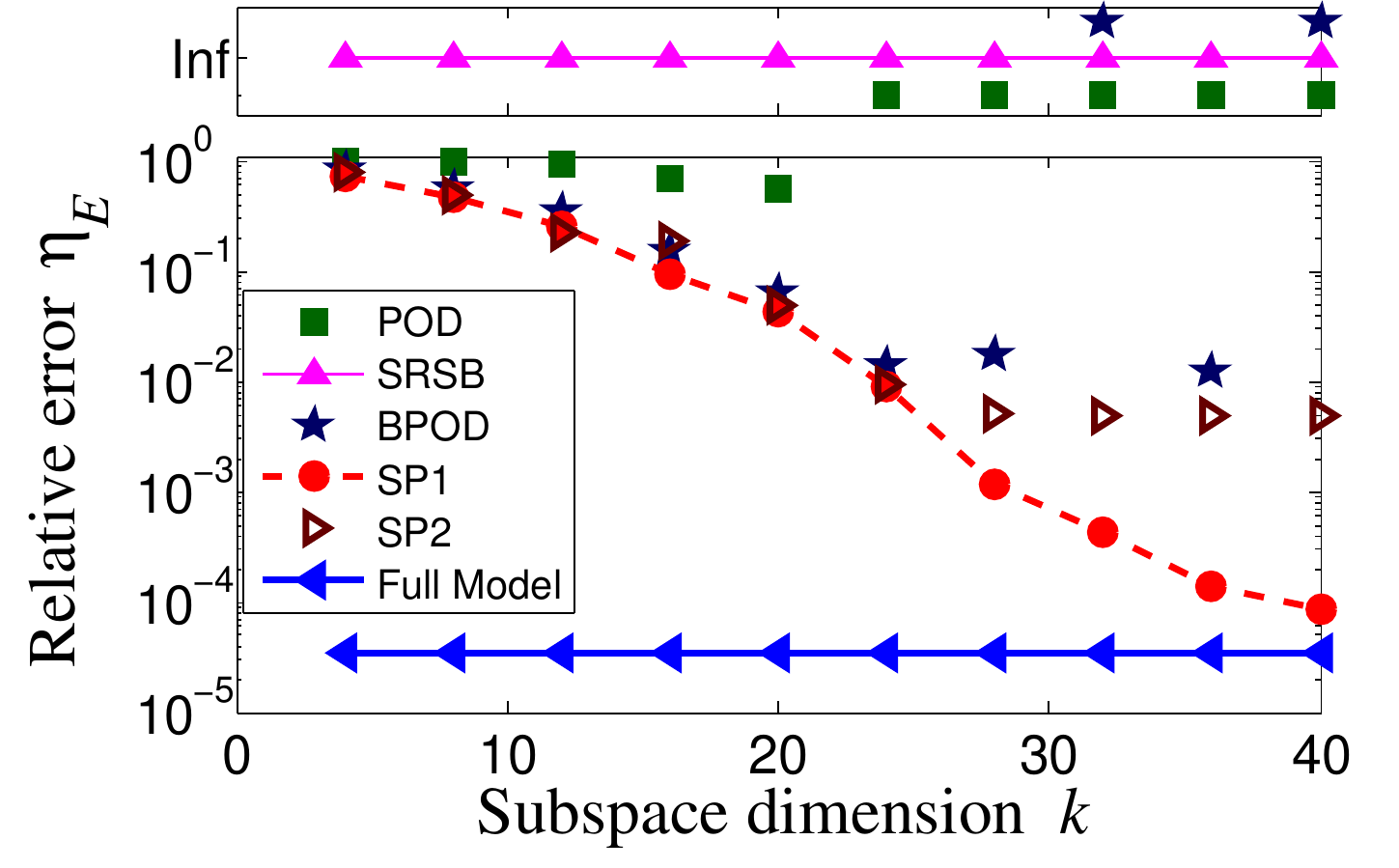}
}
 \caption{\textit{2D mass--spring example}. Method performance as a function
 of reduced dimension $k$.} 
\end{center}
\label{figure:dim_2D}
\end{figure}


Finally, we vary the reduced dimension between $k=4$ to $k=40$ to
assess the effect of subspace dimension on method performance. Figure 5.4
plots the relative state-space error $\eta$ of state variable and the
relative system-energy error $\eta_E$ as a function of $k$. Only the
full-order model
and the SP reduced-order models yield finite values of $\eta$ and $\eta_E$ for
all the tested values of subspace dimension $k$.

\section{Conclusions}\label{sec:conclusion}
This work proposed a model-reduction method that preserves marginal stability
for linear time-invariant (LTI) systems. 
	The method decomposes the LTI system into asymptotically stable and pure
	marginally stable subsystems, and subsequently performs structure-preserving
	model reduction on the subsystems separately. Advantages of the method
	include
\begin{itemize} 
\item its ability to preserve marginal stability,
\item its ability to ensure finite infinite-time energy,
\item its ability to balance primal and dual energy functionals for both
subsystems.
\end{itemize}
A geometric perspective enabled a unified comparison of the proposed
inner-product and symplectic projection methods.

Two numerical examples demonstrated the stability and accuracy of the proposed
method. In particular, the proposed method yielded a finite
infinite-time energy, while all other tested methods (i.e., POD--Galerkin,
shift-reduce-shift-back, and balanced POD) produced an infinite (unstable) or
zero (asymptotically stable) response.

\section*{Acknowledgments}
The authors thank Mohan Sarovar for his invaluable input and contributions to
this work. Sandia National Laboratories is a multi-program laboratory managed
and operated by Sandia Corporation, a wholly owned subsidiary of Lockheed
Martin Corporation, for the U.S.\ Department of Energy's National Nuclear
Security Administration under contract DE-AC04-94AL85000.

\appendix

\section{System decomposition in the general case}\label{sec:decompGeneral} We
now extend the decomposition method (in Section \ref{sec:decomp}) to a general case where the original system is unstable and $A$ is singular. 
Let $T\in \mathbb{R}^{n\times n}$ be    a nonsingular matrix such that a
 similarity transformation gives 
 \begin{equation}\label{eq:diagA_general}
 A= T
 \begin{bmatrix}
A_s & 0  & 0 & 0\\
 0 & A_m &0 &0 \\
 0 & 0 &A_u &0 \\
0 & 0 & 0 &0
 \end{bmatrix} T^{-1},
 \end{equation}
where all eigenvalues of $A_u \in \mathbb{R}^{n_u \times n_u}$ have a positive
real part.
Substituting $x = T[x_s^{\tau}\ x_m^{\tau}\ x_u^{\tau}\ x_f^{\tau}]^{\tau}$
into
	\eqref{eq:linear_control_sys} and premultiplying the first set of equations
	by $T^{-1}$ yields a decoupled LTI system
 \begin{align}\label{eq:general_relinsys}
 \begin{split}
  \frac{d}{dt}
     \begin{bmatrix}
    x_s \\ x_m \\x_u \\ x_f
     \end{bmatrix}
  &=
   \begin{bmatrix}
A_s & 0  & 0 & 0\\
 0 & A_m &0 &0 \\
 0 & 0 &A_u &0 \\
0 & 0 & 0 &0
 \end{bmatrix}
      \begin{bmatrix}
    x_s \\ x_m \\x_u \\ x_f
     \end{bmatrix}+
		\begin{bmatrix} 
		B_s\\
		B_m\\
		B_u\\
		B_f
		\end{bmatrix} 
		u\\
		y &=\begin{bmatrix}
C_s & C_m & C_u & C_f
		\end{bmatrix}
     \begin{bmatrix}
    x_s \\ x_m \\x_u \\ x_f
     \end{bmatrix},
 \end{split}
 \end{align}
where $T^{-1}B = \begin{bmatrix} B_s^{\tau}& B_m^{\tau} &B_u^{\tau}&B_f^{\tau}\end{bmatrix} ^{\tau}$ and $
		CT = 
		\begin{bmatrix}
C_s & C_m &C_u & C_f
		\end{bmatrix}
		$.
 Here, the subsystem associated
 with $ x_u $ is antistable, and the subsystem associated with $x_f$ has $0$ as system matrix and is
  marginally stable.  

In the general
case characterized by decomposition \eqref{eq:general_relinsys}, we can
perform this reduction by defining biorthogonal test and trial basis matrices for each
subsystem $\Psi_i\in\RRstar{n_i\times k_i}$, $\Phi_i\in\RRstar{n_i\times k_i}$, $i\in\{s,m,u,f\}$ . Applying Petrov--Galerkin projection to \eqref{eq:general_relinsys} with test basis  matrix $\diag(\Psi_i)$ and trial basis matrix $\diag(\Phi_i)$ yields a decoupled reduced LTI system
 \begin{align}\label{eq:general_relinsysRed}
 \begin{split}
  \frac{d}{dt}
     \begin{bmatrix}
    z_s \\ z_m \\z_u \\ z_f
     \end{bmatrix}
  &=
   \begin{bmatrix}
\tilde A_s & 0  & 0 & 0\\
 0 & \tilde A_m &0 &0 \\
 0 & 0 &\tilde A_u &0 \\
0 & 0 & 0 &0
 \end{bmatrix}
      \begin{bmatrix}
    z_s \\ z_m \\z_u \\ z_f
     \end{bmatrix}+
		\begin{bmatrix} 
		\tilde B_s\\
		\tilde B_m\\
		\tilde B_u\\
		\tilde B_f
		\end{bmatrix} 
		u\\
		y &=\begin{bmatrix}
\tilde C_s & \tilde C_m & \tilde C_u & \tilde C_f
		\end{bmatrix}
     \begin{bmatrix}
    z_s \\ z_m \\z_u \\ z_f
     \end{bmatrix},
 \end{split}
 \end{align}
where $\tilde A_i = \Psi_i^\tau A\Phi_i\in\RR{k_i\times k_i}$, $\tilde B_i =
\Psi_i^\tau B\in\RR{k_i\times p}$, $\tilde C_i =
C\Phi_i\in\RR{q\times k_i}$, $i\in\{s,m,u,f\}$.

The techniques proposed in this work can be employed to construct
$(\Psi_i,\Phi_i)$, $i\in\{s,m\}$, while bases $(\Psi_u,\Phi_u)$ can be computed to
preserve antistability in the associated reduced subsystem (e.g., via
techniques proposed in Refs.~\cite{MirnateghiN:13a,ZhouK:99a, PrakashR:90a}). Because $A_f=0$, $\tilde A_f= \Psi_f^\tau A_f \Phi_f =0$ holds for any $\Psi_f$ and $\Phi_f$. Thus, we can choose any existing method, (e.g., POD, balanced truncation, and balanced POD)   and the reduced subsystem  associated with $x_f$  always preserves pure marginal stability.

\section{Canonical form of Lyapunov equation}\label{sec:can_lya} 
We now prove a claim at the end of Section \ref{sec:asymp_stable_systems},
which  states that any Hurwitz matrix can be transformed by a similarity transformation to a matrix with
negative-definite symmetric part; in other words, a similarity transform
enables any Hurwitz matrix to satisfy the canonical Lyapunov inequality
$A_0^\tau + A_0 \prec 0$. This is in analogue to Lemma
\ref{lemma:hami_gen_can}, which shows that a generalized Hamiltonian matrix
can be transformed into a Hamiltonian matrix that satisfies the canonical
Hamiltonian property.

\begin{lemma}\label{lemma:spd_can}
A Hurwitz matrix $A$ can be transformed into a matrix $A_0$  with negative symmetric part by similarity transformation  with a real matrix $G$. Conversely, if the system matrix $A$ can be transformed into a matrix $A_0$ with negative symmetric part by a similarity transformation,  then $A$ is a Hurwitz matrix.
\end{lemma}
\begin{proof}
Because $A\in \hurwitz{n}$,  there exists  $\Theta \in \SPD{n}$ such that the Lyapunov inequality
\eqref{eq:equiLya4} holds. Choose $G\in \RRstar{n}$ (e.g., $G=\Theta^{-1/2}$) such that $G^\tau \Theta G =I_n$.
Left- and right-multiplying \eqref{eq:equiLya4} by $G^\tau$ and
$G$, respectively, yields
\begin{equation*}
(G^\tau A G^{-\tau}) (G^\tau \Theta G)   +
(G^\tau \Theta G) (G^{-1}A G)
\prec 0.
\end{equation*}
Let $A_0=G^{-1}AG$. Then,
the above equation implies $A_0^{\tau} + A_0 \prec 0$, i.e., $A_0$ has negative symmetric part. 

Conversely, suppose that $A_0 =G^{-1} A G$, where
$A_0$ has negative symmetric part and $G$ is nonsingular.  Substituting
$A_0=G^{-1}A G$ into $A_0^{\tau} + A_0 \prec 0$  yields
\begin{equation*} 
(G^{-1}A G)^{\tau} +  G^{-1}A G \prec 0.
 \end{equation*}
Left- and right-multiplying the above equation by $G^{-\tau}$
and $G^{-1}$, respectively, yields
\begin{equation*} 
A^{\tau} (G^{-\tau} G^{-1})+ (G^{-\tau}  G^{-1})A
\prec 0.
 \end{equation*}
Let  $\Theta=G^{-\tau} G^{-1}$. Then, $\Theta\in \SPD{n}$ and the above equation gives 
\eqref{eq:equiLya4}. Thus,  $A\in\hurwitz{n}$ by Lemma \ref{lemma:stabilityequiv}.
\hfill
\end{proof}

\section{Cotangent lift}\label{sec:cotangent_lift} 
The end of Section \ref{sec:PSD} mentions that there is no general way to
construct a trial basis matrix satisfying $\Phi \in \Sp(J_\Omega, J_\Pi)$.
This section briefly reviews the cotangent lift method, which is an SVD-based
method to construct $\Phi_0 \in \Sp(J_{2n}, J_{2k})$; from this matrix, a
trial basis matrix satisfying $\Phi \in \Sp(J_\Omega, J_\Pi)$ can then be
computed from Method 3 in Table \ref{tab:two_case_symplectic} using $\Phi_0$
as an input.

The cotangent lift method \cite{PengL:16a,PengL:16c} assumes that $\Phi_0$ has a block
diagonal form, i.e.,  $\Phi_0 = {\rm{diag}}(\bar \Phi, \bar \Phi)$  for some
$\bar \Phi \in \RRstar{n\times k}$. Then $\Phi_0^{\tau}J_{2n}\Phi_0=J_{2k}$
holds if and only if $\bar \Phi^{\tau}\bar \Phi=I_k$. Thus, $\bar \Phi$ is
orthonormal, i.e., $\bar\Phi\in O(I_n,I_k)$. 
Assume we have snapshots of a pure marginally stable system $\{x_i\}_{i=1}^N$;
then, we apply the inverse symplectic transformation to obtain the associated
snapshots in the canonical coordinates $\{y_i\}_{i=1}^N$ with $y_i =
G^{-1}x_i$, $i=1,\ldots,N$. Writing the decomposition $y_i = \begin{bmatrix} q_i^\tau \
p_i^\tau \end{bmatrix}^\tau\in\RR{2n}$ with $q_i,p_i\in\RR{n}$,  $\bar \Phi$ can be
computed by the SVD of an extended snapshot matrix $M_{\rm{cot}}=\begin{bmatrix} q_1
\cdots  q_N \ p_1 \cdots p_N \end{bmatrix}\in \mathbb{R}^{n\times 2N}$.

\begin{algorithm}
\caption{Cotangent lift} \label{alg:lift}
\begin{algorithmic}[1]
 \REQUIRE
Snapshots $\{x_i\}_{i=1}^N\subset\RR{2n}$ with and a
symplectic transformation matrix $G$ associated with a pure marginally stable
system.
\ENSURE A symplectic matrix $\Phi_0\in \Sp(J_{2n}, J_{2k})$ in block-diagonal form.
\STATE Apply inverse symplectic transformation to snapshots $y_i = G^{-1}x_i$,
$i=1,\ldots, N$.
\STATE Form the extended snapshot matrix
$M_{\rm{cot}}=\begin{bmatrix} q_1
\cdots\  q_N \  p_1 \cdots\  p_N\end{bmatrix}$, where $y_i = \begin{bmatrix} q_i^\tau\
p_i^\tau \end{bmatrix}^\tau$.
\STATE Compute the SVD of $M_{\rm{cot}}$; the basis matrix $\bar \Phi$
comprises the first $k$ left
singular vectors of $M_{\rm{cot}}$.
\STATE Construct the symplectic matrix $\Phi_0={\rm{diag}}(\bar \Phi, \bar \Phi)$.
\end{algorithmic}
\end{algorithm}

Algorithm \ref{alg:lift} lists the detailed procedure of the cotangent lift.
Although the cotangent lift method can only find a near optimal solution to
fit empirical data,  we can prove that  the projection error of cotangent lift
is no greater than the projection error of POD with a constant factor \cite{PengL:16c}.

\section{Generalized system energy}\label{sec:sysenergy}
If the original system is asymptotically stable, we can define a quadratic function as the system energy~\cite{RowleyCW:04a}. When it is marginally stable, we can extend the definition; the system energy is used in Section \ref{sec:1Dexample} to measure the performance of several model reduction methods. Suppose the  matrix $\Theta=M\in \SPD{n_s}$ satisfies the
Lyapunov equation \eqref{eq:semiequiLya}  for the asymptotically stable subsystem. Suppose $H: x_m \mapsto  {\frac{1} {2}}{x_m^{\tau}}Lx_m$ is the Hamiltonian function of the marginally stable subsystem with $L\in \SPD{n_m}$.
With $(x_s^{\tau}, x_m^{\tau})^{\tau}= T^{-1}x$, the system energy can be defined as
\begin{equation}\label{ref:energy}
E(t)= {\frac{1}{2}} {\begin{bmatrix} x_s(t)^{\tau} & x_m(t)^{\tau} \end{bmatrix}}   \begin{bmatrix} M & 0\\ 0& L \end{bmatrix} {\begin{bmatrix} x_s(t) \\ x_m(t) \end{bmatrix}}=\frac{1}{2}\|x_s(t)\|_M^2+ H(x_m(t)),
\end{equation}
The time evolution of the system energy is given by 
\begin{equation*}
\begin{aligned}
\frac{d}{dt} E(t) &= \frac{1}{2} (\dot x_s^{\tau} M x_s + x_s^{\tau} M \dot x_s) +\frac{1}{2} (\dot x_m^{\tau} L x_m + x_m^{\tau} L \dot x_m) 
= \frac{1}{2} x_s^{\tau}  ( A_s^{\tau} M  + M A_s )x_s +\frac{1}{2}x_m^{\tau} (  A_m^{\tau} L  + L A_m) x_m\\
&= -\frac{1}{2} x_s^{\tau} Q x_s +\frac{1}{2}x_m^{\tau} (-L  J  L  + L J  L) x_m 
= -\frac{1}{2} x_s^{\tau} Q x_s.
\end{aligned}
\end{equation*}
In the third equality, we use the Lyaponuv equation $A_s^{\tau} M  + M A_s =-Q$, the definition  $A_m=JL$,  and the fact that  $J$ is skew-symmetric. Because $Q \in \SPD{n_s}$, the system energy is strictly decreasing in time when $x_s\ne 0$.

\section{Review of existing model reduction methods}\label{sec:geoexisting}
In this section, we briefly  review a few existing model reduction methods,
including POD--Galerkin (POD), balanced truncation,  balanced POD, and
shift-reduce-shift-back (SRSB), as listed in Table \ref{tab:four_methods}.  Section \ref{sec:example} numerically compares the
performance of these methods with the proposed structure-preserving technique.
We show that each of these methods exhibits an inner-product structure (see
Table \ref{tab:comp_inner} in Section \ref{sec:existenceInner}); however, the
associated inner-product matrix $M$ does not associate with a Lyapunov matrix
in all cases, which precludes some methods from ensuring asymptotic-stability
preservation.

\begin{table}
\begin{center}
\caption{Algorithms for computing test and trial basis matrices using existing
model-reduction methods}
\label{tab:four_methods}
 \resizebox{\textwidth}{!}
 {
 \begin{tabular}{|c|c|c|c|c|}
\hline
& POD--Galerkin & Balanced truncation & Balanced POD & SRSB \\
\hline
\hline
 Input &
 \begin{tabular} {@{}l@{}}   
     Snapshots $X$ in \eqref{eq:POD_SVD}
     \end{tabular} 
    &  $(A, B, C)$ &  
     \begin{tabular} {@{}l@{}}  
     Primal snapshots $S$ in \eqref{eq:defU} \\
     Dual snapshots $R$ in  \eqref{eq:defV}
     \end{tabular} 
 &  
 \begin{tabular} {@{}l@{}}  
     $(A, B, C)$ \\
     Shift margin $\mu$
     \end{tabular}
  \\
\hline
 Output & $\Psi$, $\Phi \in O(I_n, I_k)$.  & 
 \begin{tabular} {@{}l@{}}   
 $\Phi \in O(W_o, \Sigma_1)$,\\ $\Psi \in O(W_c, \Sigma_1)$ 
\end{tabular} 
 & 
 \begin{tabular} {@{}l@{}}   
 $\Phi \in O(\hat W_o,  \Sigma_1)$; \\ 
 $\Psi \in O(\hat W_c,  \Sigma_1)$. 
 \end{tabular}
 & 
 \begin{tabular} {@{}l@{}}   
 $\Phi \in O(W_o^\mu, \Sigma_1)$,\\ $\Psi \in O(W_c^\mu, \Sigma_1)$
 \end{tabular}
 \\
\hline
Algorithm  
&
     \begin{tabular} {@{}l@{}}   
    1.  Compute SVD\\
		\quad$X = U\Sigma V^\tau$.\\
    2.  $\Psi=\Phi=U_1$.
     \end{tabular}
 
&     
     \begin{tabular} {@{}l@{}}   
    1.  Compute $W_o$ by \eqref{eq:lyapObserve1} \\
    2. Compute $W_c$ by  \eqref{eq:lyapControl1}.\\
    3.  Compute symmetric factorization\\
     \quad $W_c=S S^{\tau}$, $W_o=R R^{\tau}$. \\
    4. Compute SVD\\
		\quad $R^\tau S = U\Sigma V^\tau$. \\
    5. $\Phi = S V_1 \Sigma_1 ^{-1/2}$. \\
    6. $\Psi = R U_1\Sigma_1 ^{-1/2}$.
     \end{tabular}
 &
      \begin{tabular} {@{}l@{}}   
    1. Compute SVD \\
    \quad $R^{\tau} S =  U   \Sigma  V^{\tau}$\\
    2. $\Phi = S  V_1   \Sigma_1 ^{-1/2}$ \\
    3. $\Psi = R  U_1  \Sigma_1 ^{-1/2}$.
     \end{tabular}    
 &    
      \begin{tabular} {@{}l@{}}   
    1.  Compute $W_o^\mu$ by \eqref{eq:reWo} \\
    2. Compute $W_c^\mu$ by  \eqref{eq:reWc}.\\
    3.  Compute symmetric factorization\\
     \quad $W_c^\mu=S S^{\tau}$, $W_o^\mu=R R^{\tau}$. \\
    4. Compute SVD\\
		\quad $R^\tau S =  U \Sigma V^\tau$. \\
    5. $\Phi = S V_1 \Sigma_1 ^{-1/2}$. \\
    6. $\Psi = R U_1\Sigma_1 ^{-1/2}$.
     \end{tabular}
     \\
\hline
\end{tabular}}
\end{center}
\end{table}

\subsection{POD--Galerkin}\label{sec:POD}
POD  \cite{HolmesP:12a} computes a basis $\Phi$ that minimizes the mean-squared projection error
of a set of snapshots $\{x_i\}_{i=1}^N$, i.e., satisfies
optimality property \eqref{eq:PODoptimality_phi} with $M=I_n$ and $N = I_k$.
Algebraically, POD computes the singular value
decomposition (SVD) 
\begin{equation}\label{eq:POD_SVD}
 X = \begin{bmatrix} x_1 & \cdots & x_N \end{bmatrix} =U\Sigma  V^{\tau},
\end{equation} where $U \in O(I_n,I_r)$, $\Sigma =
\diag(\sigma_1,\ldots,\sigma_r)$ with singular values $\sigma _1 \ge ... \ge \sigma _r \ge
0$, $V \in O(I_N,I_r)$, and $r = \min (n, N)$.
Then, both the trial and test basis matrices are set to the first $k$ columns of $U$, which is equivalent to enforcing the
Galerkin orthogonality condition
(i.e., performing Galerkin projection).

As reported in Table \ref{tab:comp_inner}, it can be verified that
POD--Galerkin corresponds to an inner-product balancing with $\Xi = \Xi' =
XX^\tau$. Thus,
$\Psi,\Phi \in O(XX^T,\Sigma^2)$; however, note that we also have
$\Psi,\Phi\in O(I_n,I_k)$.

\subsection{Balanced truncation} \label{sec:BT}
Balanced truncation~\cite{MooreBC:81a} can be applied to LTI systems that are asymptotically
stable, controllable, and observable. 
In the present framework, balanced truncation corresponds to a specific type of inner-product balancing
with $\Xi=W_o$ and $\Xi'=W_c$, where $W_o$ and $W_c$ represent observability
and controllability Gramians that satisfy primal and dual
Lyapunov equations
\begin{align}
A^{\tau} W_o + W_o A  &= -C^{\tau} C,
\label{eq:lyapObserve1} 
\\
A W_c + W_c A ^{\tau} &= -B B^{\tau},
\label{eq:lyapControl1} 
\end{align}
respectively, which are defined for observable and controllable asymptotically
stable LTI systems.  While the present framework cannot prove that balanced
truncation preserves asymptotic stability (the right-hand-side matrices in the
Lyapunov equations \eqref{eq:lyapObserve1}--\eqref{eq:lyapControl1} are positive \textit{semidefinite}), it can be shown using
observable and controllable conditions that balanced truncation
does in fact preserve asymptotic stability~\cite[pp.~213--215]{AntoulasAC:05a}.

\subsection{Balanced POD}\label{sec:BPOD}
 Several techniques exist to solve the Lyapunov equations
 \eqref{eq:lyapObserve1} and \eqref{eq:lyapControl1} for the controllability and observability Gramians~\cite{AntoulasAC:05a};
however, they are  prohibitively expensive for large-scale systems. For this
reason, several methods have been developed that instead employ \textit{empirical
Gramians} \cite{lall2002sab} that approximate the analytical Gramians. One
particular method is balanced POD (BPOD) 
\cite{willcox2002bmr,RowleyCW:05a}, which relies on collecting primal snapshots
for $N$ time steps (one impulse response of the forward system per column in
$B$):
 \begin{equation} \label{eq:defU}
S = \begin{bmatrix}
B&  e^{A \Delta t}B&  \cdots & e^{A (N-1) \Delta t} B
\end{bmatrix}
  \end{equation} 
and dual snapshots for $N$ time steps (one impulse response of the dual system per row in $C $):
 \begin{equation}  \label{eq:defV}
R= \begin{bmatrix}
C^{\tau} & e^{A^{\tau} \Delta t}C^{\tau} & \cdots&  e^{ A^{\tau}  (N-1) \Delta t}C^{\tau}
\end{bmatrix};
\end{equation} 
the empirical observability and controllability Gramians are then set to
$\hat W_o = R R^{\tau}$ and
$\hat W_c= S S^{\tau}$, respectively.  Then, 
BPOD corresponds to inner-product balancing with $\Xi= \hat W_o$ and $\Xi' =
\hat W_c$. 

Critically, empirical Gramians $\hat W_o$ and $\hat W_c$  may not be Lyapunov
matrices, i.e., they may not satisfy \eqref{eq:equiLya4} and \eqref{eq:dualinequiLya}, respectively.
Consequently, unlike balanced truncation, balanced POD does not guarantee
asymptotic-stability preservation.

 \subsection{SRSB}\label{sec:SRSB} 
The shift-reduce-shift-back (SRSB) method aims to extend the applicability of
balanced truncation to marginally stable and unstable systems~\cite{YangJ:92a,
DoyleJC:02a,SantiagoJM:85a, ZilouchianA:91a, YangJ:93a}.  

By Lemma \ref{lemma:stabilityequiv},  the real parts of the eigenvalues of $A$
are less than the shift margin $\mu\in\RR{}$ if and only if given any
$Q\in\SPD{n}$, there exists $\Theta \in\SPD{n}$ that is the
unique solution to
\begin{equation*}
 (A- \mu I )^{\tau} \Theta + \Theta (A- \mu I ) = -Q.
\end{equation*}
Thus, even if $A$ is marginally stable or unstable, we can choose $\mu$ such that the shifted system corresponding to the matrix $A- \mu
I $ is asymptotically stable and the test and trial basis matrices can be
computed by performing balanced truncation on the shifted system. The basis matrices can be applied to
  the original (unshifted) system. Table 
\ref{tab:four_methods} provides the associated algorithm, which amounts to
computing an 
inner-product balancing with $\Xi = W_o^\mu$ and
$\Xi' = W_c^\mu$, where the shifted observability and controllability Gramians
satisfy
\begin{align}
(A- \mu I)^{\tau} W_o^\mu + W_o^\mu(A- \mu I)  &= -C^{\tau}C,
\label{eq:reWo}
\\
 (A-\mu I)W_c^\mu + W_c^\mu(A-\mu I)^{\tau} &= -BB^{\tau}.
\label{eq:reWc}
\end{align}


While SRSB ensures that the shifted reduced system  $(\Psi^{\tau} (A-\mu
I_n)\Phi, \Psi^{\tau} B,  C \Phi)=(\tilde A - \mu I_k, \tilde B, \tilde C)$
remains asymptotically stable, this guarantee does not extend to the reduced
system $(\tilde A, \tilde B, \tilde C)$ that is used in practice. In
particular, there is no assurance that the reduced system will retain the
asymptotic or marginal stability  that characterized the original system
$(A,B,C)$.

\bibliographystyle{siam} \bibliography{RefA1} \end{document}